\input amstex
\documentstyle{amsppt}
\magnification 1200
\NoBlackBoxes
\input epsf
\vcorrection{-9mm}

\topmatter

\title
       Diffusion orthogonal polynomials in 3-dimensional
       domains bounded by developable surfaces
\endtitle

\author
       S.~Yu.~Orevkov
\endauthor

\abstract The following problem is studied: describe the triplets $(\Omega,g,\mu)$,
$\mu=\rho\,dx$, where $g= (g^{ij}(x))$ is the (co)metric associated
with the symmetric second order differential operator
$\bold L(f) = \frac{1}{\rho}\sum_{ij} \partial_i (g^{ij} \rho\,\partial_j f)$ 
defined on a domain $\Omega$ of $\Bbb R^n$ and such that there exists an
orthonormal basis  of $\Cal L^2(\mu)$ made of polynomials which are eigenvectors
of $\bold L$, and the basis is compatible with the filtration of the space
of polynomials with respect to some weighted degree.

In a joint paper with D.~Bakry and
M.~Zani this problem was solved in dimension 2 for the usual degree.
In the author's subsequent paper this problem was
solved in dimension 2 for any weighted degree.
In the present paper this problem is solved in dimension 3 for the usual degree
under the condition that $\partial\Omega$ contains a piece of
a tangent developable surface. The proof is based on Pl\"ucker-like formulas
in the form given by Ragni Piene.
All the found solutions are generalized for any dimension.
\endabstract

\address
IMT, Univ.~Paul Sabatier, 118 roure de Narbonne, Toulouse, France
\endaddress
\vskip-8pt
\address
Steklov Math. Inst., ul. Gubkina 8, Moscow, Russia
\endaddress
\endtopmatter

\rightheadtext{Diffusion orthogonal polynomials in 3-dimensional domains}

\def\refBB      {1}
\def\refBGL     {2}
\def\refBOZ     {3}
\def\refBou     {4}  
\def\refCox     {5}  
\def\refMe      {6}  
\def\refOr      {7}
\def\refPie     {8}
\def\refPieMA   {9}
\def\refSYS     {10}
\def\refSouArx  {11}
\def\refSouAFST {12}

\def\eqLL       {1}
\def\eqLLrho    {2}
\def\eqParam    {3}

\def\eqPie      {4}
\def\eqPieA     {5}
\def\eqPieB     {6}

\def\eqParamGen {7}

\def\eqTCubic   {8}
\def\eqPhiMu    {9}
\def\eqTQuartic {10}

\def\eqIV       {11}
\def\eqV        {12}
\def\eqVI       {13}

\def\eqModelAnI   {14}
\def\eqModelAnII  {15}
\def\eqModelAnIII {16}
\def\eqModelAnIV  {17}
\def\eqModelBn    {18}
\def\eqModelBnDe  {19}
\def\eqModelAA    {20}

\def\ctCubic {$1^\circ$}
\def\ctCQuar {$2^\circ$}
\def\ctOFQua {$3^\circ$}
\def\ctTFQua {$4^\circ$}
\def\ctFQuin {$5^\circ$}
\def\ctGQuar {$6^\circ$}
\def\ctSQuin {$7^\circ$}
\def\ctCSext {$8^\circ$}


\def\sectTan {2}
\def\sectLoc {3}
\def\sectAlg {4}
\def\sectDOP {5}
\def\sectFig {5.3}
\def\sectQuo {6}
\def\sectGa {6.1}
\def\sectA  {6.2}
\def\sectB  {6.3}
\def\sectBb {6.4}
\def\sectD  {6.5}
\def\sectJS {6.6}
\def\sectAA {6.7}
\def\sectAB {6.8}
\def\sectCt {6.9}
\def\sectBDt{6.10}
\def\sectAt {6.11}
\def\sectCone {7}


\def\propPie  {2.1}
\def\lemCoNo  {2.2}
\def\lemCompl {2.2(a)}
\def\lemNodal {2.2(b)}

\def\lemBrFlexAff    {3.1}  
\def\lemBrFlatAff    {3.2}  
\def\lemBrGen        {3.3}  
\def\lemBrFlex       {3.4}  
\def\lemBrFlat       {3.5}  
\def\lemBrDblFlat    {3.6}  
\def\lemBrGenTan     {3.7}  
\def\lemBrFlatTan    {3.8}  
\def\lemBrDblFlatTan {3.9}  
\def\lemBrCusp       {3.10} 
\def\lemBrFlatCusp   {3.11} 
\def\lemBrCuspTan    {3.12} 

\def\propAlgC   {4.1}
\def\propAlgR   {4.2}
\def\lemCubicI  {4.3}
\def\lemCubicII {4.4}
\def\lemQuartic {4.5}
\def\lemQuiSext {4.6}

\def\thB       {5.1}
\def\thU       {5.2}
\def\remCubicI {5.3}

\def\thC   {7.1}
\def\remC  {7.2}
\def\propC {7.3}
\def\lemC  {7.4}

\def\tabPie   {1}

\def\figAB     {1}
\def\figAiAii  {2}
\def\figAaff   {3}
\def\figCaff   {4}
\def\figABD    {5}
\def\figIrrel  {6}
\def\figCone   {7}

\def\bR{\Bbb{R}}
\def\bC{\Bbb{C}}

\def\CP{\Bbb{CP}}
\def\bK{\Bbb{K}}
\def\bP{\Bbb{P}}

\def\bS{\Bbb{S}}

\def\ord{\operatorname{ord}}

\def\LL{\bold{L}}
\def\cL{\Cal{L}}

\def\cP{\Cal{P}}

\def\xx{\bold x}

\def\ii{\bold i}
\def\gg{\bold g}
\def\nn{\bold n}
\def\ff{\bold f}

\def\eps{\varepsilon}

\loadbold
\def\De{{\boldsymbol\Delta}}
\def\Ga{{\boldsymbol\Gamma}}

\document

\head 1. Introduction
\endhead

This paper continues the study of the diffusion orthogonal polynomials started
in [\refBOZ] (see also [\refBB], [\refOr], [\refSouArx], [\refSouAFST]).
It is devoted to the following problem posed by Dominique Bakry: describe
all triples $(\Omega,\LL,\mu)$ where $\Omega$
is a domain in $\bR^n$ such that $\Omega=\operatorname{Int}\overline{\Omega}$,
$\LL$ is an elliptic second order operator of the form
$$
      \LL(f) = \sum_{i,j} g^{ij}(x)\partial_{ij}f
             + \sum_ib^i(x)\partial_i f                                 \eqno(\eqLL)
$$
with $g^{ij}$ and $b^i$ continuous
in $\Omega$, and $\mu=\rho\,dx$ a probability measure on $\Omega$
with $\Cal C^1$-smooth density $\rho$,
and such that there exists a polynomial orthogonal basis of
$\cL^2(\Omega,\mu)$ formed by eigenvectors of $\LL$, which is also a basis
(in the algebraic sense) of $\bR[x]$, $x=(x_1,\dots,x_n)$,
and which is compatible with the filtration of $\bR[x]$ by the degree
(a variant: by a weighted degree, see [\refBB], [\refOr]). The latter
condition means that the space
$\cP_m$ of polynomials of degree $\le m$ is $\LL$-invariant for any $m$.
We say that such a triple $(\Omega,\LL,\mu)$
is a solution of the {\it Diffusion Orthogonal Polynomial problem}
(DOP problem for short).
If in addition
$\int_\Omega f_1\LL f_2\,d\mu = \int_\Omega f_2\LL f_1\,d\mu$
for any pair of compactly supported functions
(for bounded domains this condition follows from the other ones), we say that
$(\Omega,\LL,\mu)$ is a solution of the {\it strong DOP problem} (SDOP problem
for short). In this case
$$
    \LL(f) = \frac1\rho\sum_{i,j}\partial_i\Big(g^{ij}\rho\,\partial_j f\Big),
                                                                        \eqno(\eqLLrho)
$$
thus
$\LL$ is determined by $g=(g^{ij})$ and $\rho$, and we therefore
talk about $(\Omega,g,\rho)$ as a solution of the SDOP problem.
If $\rho=(\det g)^{-1/2}$, then $\LL$ given by (\eqLLrho) is the Laplace-Beltrami
operator for the metric $(g_{ij})=g^{-1}$.


As shown in [\refBOZ, Thm.~2.21], $(\Omega,g,\rho)$ is a solution of the SDOP problem
(and hence of the DOP problem when $\Omega$ is bounded) if and only if
there exists a squarefree polynomial $\Gamma$ such that:
\roster
\item"(A1)" $g^{ij}\in\cP_2$ for each $i,j=1,\dots,n$;
\item"(A2)" $\Gamma$ divides $\det g$;
\item"(A3)" $\Gamma$ divides $\sum_j g^{ij}\partial_i\Gamma$ for each $i=1,\dots,n$;
\item"(A4)" $\partial\Omega\subset\{\Gamma=0\}$
            and $g|_\Omega$ is positive definite;
\item"(A5)" $\sum_j g^{ij}\partial_i\log\rho\in\cP_1$ for each $i=1,\dots,n$;
\item"(A6)" polynomials are dense in $\cL^2(\Omega,\rho\,dx)$.
\endroster
Condition (A3) is equivalent to the fact that for any
germ $\xi:(\bR^{n-1},0)\to(\bR^n,x)$
such that $\Gamma\circ\xi=0$, one has
$$
   \xi^*(\omega_i)=0, \quad
   \omega_i =\sum_j (-1)^j g^{ij}\,
         dx_1\wedge\dots\wedge\widehat{dx}_j\wedge\dots\wedge dx_n,
       \quad i=1,\dots,n.
                                                      \eqno(\eqParam)
$$
Note that Conditions (A1)--(A3) are purely algebraic and they make sense for
polynomials with coefficient in any field $\bK$. If they are satisfied for a field $\bK$,
we say that $(g,\Gamma)$ is a solution of the
{\it algebraic counterpart of the DOP problem over} $\bK$
(AlgDOP/$\bK$ problem for short).

In dimension 2, all solutions of the DOP problem are found in [\refBOZ]
for the usual degree and in [\refOr]
for any weighted degree. In the present paper we attack the classification
of the solutions in dimension 3 for the usual degree.
By (A3), $\partial\Omega$ sits on an algebraic hypersurface
of degree at most $2n$, thus on a quartic curve when $n=2$.
The arguments in [\refBOZ] essentially rely on the Pl\"ucker-like formulas
relating the singularities of this curve and those of its projectively dual curve.
It seems that this approach can be also applied at least in dimension 3.
Here we take the first step in this direction.
Namely, we describe all
irreducible surfaces $\Sigma$ in $\bR^3$ whose projective dual has dimension $1$
(i.e., is a curve, which we denote by $\check C$)
and such that a relatively open piece of $\Sigma$ appears in $\partial\Omega$
for some solution $(\Omega,g,\rho)$ of the SDOP problem. Moreover, in the case
when $\check C$ is not contained in any plane (in this case $\Sigma$ is the tangent
developable of another curve $C$ called the dual curve of $\check C$), we
describe all such solutions $(\Omega,g,\rho)$ (Theorems~\thB\ and \thU).
If $\check C$ is contained in some plane, then $\Sigma$ is a cylinder
or a cone over some planar curve $A$. If $\Sigma$ is a cylinder, then it is easy to
show that a piece of $A$ occurs in the boundary of a two-dimensional solution.
If $\Sigma$ is a cone, we prove in Theorem~\thC\ that $\deg A=2$, thus $\Sigma$ is
a standard quadratic cone. In the conical case there indeed exist some solutions 
(see Remark~\remC) but we do not know
if this list is exhaustive.

To prove Theorems \thB\ and \thU, we follow the strategy similar to that
in [\refBOZ, \S3].
Condition (\eqParam) yields equations for the coefficients of the
polynomials $g^{ij}$ and those of local parametrizations of the curve $C$.
By solving them we obtain in \S\sectLoc\ rather strong restrictions on a priori possible
types of local branches (real and complex) of $C$. Then in \S\sectAlg,
using Pl\"ucker-like formulas due to Ragni Piene [\refPie] (introduced in \S\sectTan),
we find all solutions of the AlgDOP problem
over $\bC$, and then (in \S\sectDOP) we
find $\Omega$, $\rho$, and the real form of $g$ satisfying
the remaining conditions (A4)--(A6).

In \S\sectQuo, for each bounded domain in Theorem~\thB, we show that the
Laplace-Beltrami solution is the image of Euclidean or spherical Laplace operator
through on an appropriate realization of quotient of $\bR^3$ or $\bS^3$ by
a Coxeter group, and we generalize this construction to any dimension.
In \S\sectCone\ we prove the aforementioned result about conical surfaces.

\head\sectTan. Tangent developables.
\endhead

Let $C=\nu(\widetilde C)$ be an irreducible algebraic curve
in $\bP^3$ of genus $g$ which is
not contained in a plane (here $\widetilde C$ is a smooth compact Riemann surface
and $\nu:\widetilde C\to\bP^3$ an analytic mapping).
Let $\Sigma$ be the {\it tangent developable surface of $C$},
i.e., $\Sigma$ is the union of all
lines tangent to $C$.

Following [\refPie], we introduce the following notation.
For any point $p\in\widetilde C$ there exists a local affine chart of $\bP^3$
centered at $\nu(p)$ such that the corresponding local branch of $C$ is parametrized
by $t\mapsto(t^{m_0},t^{m_1},t^{m_2})$ with
$$
    (m_0,m_1,m_2)=(1+l_0,\,2+l_0+l_1,\,3+l_0+l_1+l_2), \qquad l_j=l_j(p)\ge 0.
$$
Then we say that $p$ is a point of type $(m_0,m_1,m_2)$ and
we set $k_j = k_j(C) = \sum_{p\in\widetilde C} l_j(p)$, $j=0,1,2$.
We also denote the {\it osculating} plane at $p$ by $O_p$. In the above coordinates,
this is the plane spanned by $(1,0,0)$ and $(0,1,0)$.

The curve $\check C$ in the dual projective space $\check\bP^3$ parametrized by
$\check\nu:\widetilde C\to\check\bP^3$, $p\mapsto O_p$ is called the {\it dual curve} of $C$.
Let $r_0$, $r_1$, $r_2$ be the degrees of $C$, $\Sigma$, $\check C$ respectively
(see [\refPie] for a more uniform definition). It is immediate to check that
the dual of $\check C$ is $C$ and
$k_j(\check C)=k_{2-j}(C)$, $r_j(\check C)=r_{2-j}(C)$, $j=0,1,2$.
The classical Pl\"ucker--Cayley equations in the form given by Ragni Piene
in [\refPie], [\refPieMA,~Eq.~(1)] read
as follows:
$$
  \matrix
    r_1=2r_0+2g-2-k_0,       \hskip41pt& r_1=2r_2+2g-2-k_2,       \hskip41pt\\
    r_2=3(r_0+2g-2)-2k_0-k_1,\hskip7pt & r_0=3(r_2+2g-2)-2k_2-k_1,\hskip7pt \\
    k_2=4(r_0+3g-3)-3k_0-2k_1,         & k_0=4(r_2+3g-3)-3k_2-2k_1.
  \endmatrix                                                        \eqno(\eqPie)
$$
Any three of these equations imply the others.

\proclaim{ Proposition \propPie }
If $r_1\le 6$, then one of the cases listed in Table~1 takes place.
\endproclaim

\demo{ Proof }
If $r_0\le 3$ (recall that $r_0=\deg C$), then the only non-planar curve
(up to automorphism of $\bP^3$) is
the rational cubic parametrized by $t\mapsto(1:t:t^2:t^3)$. It
corresponds to Case~$1^\circ$. Then assume that $r_0\ge 4$ and (by the duality)
$r_2\ge 4$.

We assume for simplicity that $l_0(p)+l_1(p)+l_2(p)\le 1$ for
each point $p$ of $C$, i.~e., each point of $C$ contributes at most $1$
to $k_0+k_1+k_2$. It is not difficult to adapt the proof for the general case.

The degree of the cuspidal edge of $\Sigma$ is $r_0+k_1$
(see [\refPieMA, p.~112, l.~15~ff]), hence the genus formula for a
generic plane section of $X$ yields
$$
       r_0+k_1 \le (r_1-1)(r_1-2)/2.                       \eqno(\eqPieA)
$$

If $k_2>0$, then we consider the plane projection of $C$
from one of its cusps. Its degree is $r_0-2$ and it has $k_2-1$ cusps.
Hence, by the genus formula,
$$
      k_2>0 \quad\Rightarrow\quad g+k_0-1 \le (r_0-3)(r_0-4)/2. \eqno(\eqPieB)
$$
One can easily check that in Table~\tabPie\ there listed all non-negative
integer solutions $(g,k_0,k_1,k_1,r_0,r_1,r_2)$ of the equations (\eqPie)
combined with the inequalities (\eqPieA), (\eqPieB),
$r_0\ge 4$, $r_2\ge 4$, and $3\le r_1\le 6$.
\qed\enddemo

\midinsert
\vbox{
\def\TL#1#2#3#4#5#6#7#8#9{%
   \noindent\hskip0pt
 \hbox to 15mm{\hskip7mm#1\hskip5pt\hfill}
 \hbox to 53mm{#2\hfill}
 \hbox to  9mm{#3\hfill}
 \hbox to  6mm{#4\hfill}
 \hbox to  6mm{#5\hfill}
 \hbox to  9mm{#6\hfill}
 \hbox to  6mm{#7\hfill}
 \hbox to  6mm{#8\hfill}
 \hbox to  6mm{#9\hfill}}
\hrule
\vskip3pt
\TL{no.}{}{$g$}{$k_0$}{$k_1$}{$k_2$}{$r_0$}{$r_1$}{$r_2$}
\vskip3pt
\hrule
\vskip4pt
\TL{\ctCubic}{Twisted cubic                   } 0  0 0 0  3 4 3 \par
\TL{\ctCQuar}{Cuspidal quartic                } 0  1 0 1  4 5 4 \par
\vskip4pt
\TL{\ctOFQua}{Once inflected quartic          } 0  0 1 2  4 6 5 \par
\TL{\ctTFQua}{Twice inflected quartic         } 0  0 2 0  4 6 4 \par
\TL{\ctFQuin}{Inflected bicuspidal quintic    } 0  2 1 0  5 6 4 \par
\vskip4pt
\TL{\ctGQuar}{Generic quartic                 } 0  0 0 4  4 6 6 \par
\TL{\ctSQuin}{Non-inflected bicuspidal quintic} 0  2 0 2  5 6 5 \par
\TL{\ctCSext}{Four-cuspidal sextic            } 0  4 0 0  6 6 4 \par
\vskip4pt
\hrule
}   
\botcaption{ Table \tabPie } Curves whose tangent
                             developables have degree 6
\endcaption
\endinsert

\proclaim{ Lemma \lemCoNo } Suppose that $C$ is rational.
Let $p,q\in\widetilde C$, $p\ne q$, be points of the types $(m_0,m_1,m_2)$
and $(m'_0,m'_1,m'_2)$. Recall that $r_0=\deg C$.

\smallskip
(a). If
$ m'_2 \le m'_1+m_0 = m'_0+m_1 = m_2 = r_0 $,
then $C$ has parametrization
$t\mapsto\big(\sum_{j=0}^{r_0-m'_2}a_jt^j:t^{m_0}:t^{m_1}:t^{m_2}\big)$,
$a_0a_{r_0-m'_2}\ne 0$, in some homogeneous coordinates.

\smallskip
(b). If $r_0=4$ and $\nu(p)=\nu(q)$, then $C$ has parametrization
$t\mapsto\big(1+t^4:t:t^2:t^3\big)$ in some homogeneous coordinates.
\endproclaim

\demo{ Proof }
We have $\widetilde C=\bP^1$ and we may assume that $p=(0:1)$, $q=(1:0)$,
and the mapping $\nu$ is given by
$(t:s)\mapsto(f_0(t,s):\dots:f_3(t,s))$,
where $f_j$ are homogeneous polynomials of degree $r_0$
and $\ord_t(f_0,\dots,f_3)=(0,m_0,m_1,m_2)$.

\smallskip
(a). The condition $m_2=r_0$ implies that $f_3=t^{m_2}$ up to rescaling.
By a coordinate change $f_j\to f_j-c_j f_3$, $j=0,1,2$, we may attain
$\ord_s f_j\ge m'_0$ for $j\ge 2$. Then the condition $m'_0+m_1=r_0$ implies
that $f_2=t^{m_1}s^{m'_0}$ up to rescaling.
Proceeding in this way we arrive to the required parametrization.

\smallskip
(b). The condition $\nu(p)=\nu(q)$ implies $f_1(q)=f_2(q)=f_3(q)=0$, i.e.,
$\deg_t f_j\le 3$ for $j=1,2,3$. Then by coordinate changes $f_i\to f_i-c_{ij}f_j$,
$i<j$, we arrive to the required parametrization.
\qed\enddemo


\head\sectLoc.     Restrictions on local branches
\endhead

Let the notation be as in \S\sectTan\ but we fix an affine chart in $\bP^3$
with coordinates $(x,y,z)$. We denote the plane at infinity by $P_\infty$.
Let $\Gamma(x,y,z)=0$ be the equation of $\Sigma$.
Suppose that there exists a cometric $g=(g^{ij})$ such that $(g,\Gamma)$ is
a solution of the SDOP problem. We denote the coefficient of $x^k y^l z^m$ in $g^{ij}$
by $g^{ij}_{klm}$.
Let 
$$
   t\mapsto\gamma(t)=(\xi_1(t),\xi_2(t),\xi_3(t)) 
$$
be a local meromorphic branch of $C$ at a finite or infinite point.
Then $\Sigma$ admits parametrization
$$
    (t,u)\mapsto\big( \hat\xi_1(t,u),\hat\xi_2(t,u),\hat\xi_3(t,u) \big),
    \qquad  \hat\xi_j = \xi_j + u\dot\xi_j,
$$
at a neighbourhood of the line tangent to $C$ at $\gamma(0)$.
Then the equations (\eqParam) take the form $E_1=E_2=E_3=0$ where
$$
   E_i
    =\sum_{j=1}^3
    \frac{\partial(\hat\xi_{j+1},\hat\xi_{j-1})}{\partial(t,u)}
    g^{ij}(\hat\xi_1,\hat\xi_2,\hat\xi_3)
    =u\sum_{j=1}^3
      \big(\ddot\xi_{j+1}\dot\xi_{j-1}
         - \ddot\xi_{j-1}\dot\xi_{j+1} \big)
    g^{ij}(\hat\xi_1,\hat\xi_2,\hat\xi_3)
$$
(here the indices are considered mod 3).
We have
$E_i = \sum_{\alpha=\alpha_0}^{\infty} t^\alpha
           \sum_{\beta=1}^3 E_{\alpha,\beta,i} u^\beta$
where the $E_{\alpha,\beta,i}$ are linear forms in $g_{klm}^{ij}$ whose
coefficients are polynomial functions of the coefficients of the $\xi_i$'s.

In the following Lemmas~\lemBrFlexAff--\lemBrCuspTan, for several
a priori possible values of $\ord_t(\gamma)$, we either exclude them
or (in Lemma~\lemBrGen) show that the given value implies a certain explicit form of $C$.
In all the proofs (except those of Lemmas~\lemBrFlexAff--\lemBrFlatAff)
we assume that $\gamma$ is parametrized by $t\mapsto(x,y,z)$,
$$
   x=t^{j_1},\qquad y=t^{j_2}+\sum_{j>j_2}b_jt^j,\qquad
                    z=t^{j_3}+\sum_{j>j_3}c_jt^j,\qquad j_1>j_2>j_3,
$$
and, moreover, $b_0=b_{j_1}=c_0=c_{j_1}=c_{j_2}=0$. The latter condition
can be easily achieved by the change of variables $(y,z)\to(y_1,z_1)$,
$y_1=y-b_0-b_{j_1}x$, $z_1=z-c_0-c_{j_2}y_1-c_{j_1}x$.
Then we solve a system of some $n$ linear equations $E_{\alpha,\beta,i}=0$
for some $n$ unknowns $g^{ij}_{klm}$ whose determinant is a nonzero constant.
The number $n$ and the choice of the equations and unknowns is indicated in each proof.
In most cases the solution plugged into $g$ implies that $x^2$
divides $\det g$ which contradicts the condition $\deg\Gamma\ge 5$.
In other cases we then solve some few additional equations.

\proclaim{ Lemma \lemBrFlexAff }
If $\deg\Gamma\ge 5$, then $\ord_t(\gamma)\ne(1,3,4)$.
\endproclaim

\demo{ Proof }
We choose a parametrization of the form
$x=t$, $y=t^3+\sum_{\nu\ge 4}b_\nu t^\nu$,  $z=t^4+\sum_{\nu\ge 5}c_\nu t^\nu$.
By the change of variables $y\to y-b_4 z$ we make $b_4=0$.

All variables $g^{ij}_{klm}$ except
$g^{11}_{klm}$ with $k+l+m=2$,
$g^{12}_{0lm}$ with $l+m=2$,
$g^{13}_{002}$, and $g^{13}_{101}$ (thus $49$ variables)
can be found 
by solving the following system of $49$ equations:
$E_{1,\beta,i}$  (all $\beta$, $i$);
$E_{3,2,i}$, $E_{3,3,i}$, $E_{4,2,i}$, $E_{4,3,i}$, $E_{5,3,i}$ ($i=1,2,3$);
$E_{2,1,i}$, $E_{2,2,i}$, $E_{2,3,i}$,
$E_{5,1,i}$, $E_{5,2,i}$, $E_{6,1,i}$, $E_{6,2,i}$, $E_{7,1,i}$ ($i=1,2$);
$E_{6,3,i}$, $E_{7,3,i}$ ($i=2,3$);
$E_{4,1,1}$, $E_{7,2,2}$, $E_{8,2,2}$, $E_{8,3,2}$, $E_{9,3,2}$. 
\if01{
$E_{1,\beta,i}$  (all $\beta$, $i$);
$E_{3,2,i}$, $E_{3,3,i}$, $E_{4,2,i}$, $E_{4,3,i}$ ($i=1,2,3$);
$E_{2,\beta,i}$, 
$E_{5,\beta,i}$,  $E_{6,1,i}$, $E_{6,2,i}$, $E_{7,1,i}$ ($\beta=1,2,3$ and $i=1,2$);
$E_{6,3,i}$, $E_{7,3,i}$ ($i=2,3$);
$E_{4,1,1}$, $E_{5,3,3}$, $E_{7,2,2}$, $E_{8,2,2}$, $E_{8,3,2}$, $E_{9,3,2}$. 
$E_{8,2,2}$, $E_{8,3,2}$, $E_{9,3,2}$, and all the equations of the form
$E_{\alpha,\beta,i}$ with $1\le\alpha\le 7$ except
$E_{2,1,3}$, $E_{2,2,3}$, $E_{2,3,3}$, 
$E_{3,1,1}$, $E_{3,1,2}$, $E_{3,1,3}$,
$E_{4,1,2}$, $E_{4,1,3}$,
$E_{6,3,1}$, $E_{7,2,1}$, $E_{7,3,1}$, and
$(E_{\alpha,\beta,3})_{\beta=1,2}^{\alpha=5,6,7}$.
}\fi
The determinant of this system is a non-zero constant.
By plugging the solution to $E_{8,3,3}$ we obtain the equation
$36 c_5 g^{13}_{101}=0$. This equation implies that $z^2$ divides $\det g$
which contradicts the condition $\deg\Gamma\ge 5$.
\qed\enddemo

\proclaim{ Lemma \lemBrFlatAff }
If $\deg\Gamma\ge 5$, then $\ord_t(\gamma)\ne(1,2,4)$.
\endproclaim

\demo{ Proof }
We choose a parametrization of the form
$x=t$, $y=t^2+\sum_{\nu\ge 3}b_\nu t^\nu$,  $z=t^4+\sum_{\nu\ge 5}c_\nu t^\nu$.
By the change of variables $y\to y-b_4 z$ we make $b_4=0$.
We solve $40$ equations for $40$ unknowns.
The equations:
$E_{0,\beta,i}$ (all $\beta$, $i$);                
$E_{1,2,i}$, $E_{1,3,i}$, $E_{2,3,i}$,
$E_{3,2,i}$, $E_{3,3,i}$, $E_{4,3,i}$ ($i=1,2,3$); 
$E_{2,1,i}$, $E_{2,2,i}$,
$E_{4,1,i}$, $E_{4,2,i}$, $E_{5,1,i}$ ($i=1,2,3$); 
$E_{3,1,1}$, $E_{5,2,3}$, $E_{5,3,3}$. 
The unknowns:
$g^{ij}_{kl0}$ ($1\le i\le j\le 2$, $k,l\ne 2$),                       
$g^{12}_{200}$, $g^{22}_{200}$,                                        
and all the $g^{i3}_{klm}$
except $g^{13}_{011}$, $g^{13}_{002}$, $g^{23}_{002}$, $g^{33}_{002}$. 
The determinant is a nonzero constant.
Plugging the solution to $g$, we obtain that $z^2$ divides $\det g$
which contradicts $\deg\Gamma\ge 5$.
\qed\enddemo

\proclaim{ Lemma \lemBrGen }
If $\deg\Gamma\ge 5$ and $\ord_t(\gamma)=(-1,1,2)$
(i.e., $\gamma$ is a generic branch transverse to $P_\infty$),
 then
$\Gamma$ is parametrized by
$$
   t\mapsto(t^{-1} + t,\, 3t - t^3,\, 2t^2 - t^4)  \eqno(\eqParamGen)
$$
in some affine coordinates.
\endproclaim

\demo{ Proof }
$n=58$.
The equations:
$E_{-3,\beta,i}$, $E_{-1,\beta,i}$, $E_{0,\beta,i}$  (all $\beta$, $i$);
$E_{-7,3,i}$, $E_{-5,2,i}$, $E_{-5,3,i}$, $E_{-4,3,i}$, $E_{-2,3,i}$, $E_{1,2,i}$
($i=1,2,3$);
$E_{-6,3,i}$, $E_{-4,2,i}$, $E_{-2,1,i}$, $E_{1,1,i}$, $E_{1,3,i}$, $E_{2,1,i}$ ($i=2,3$);
$E_{2,3,3}$.
The unknowns are all the
$g^{ij}_{klm}$ except $g^{33}_{002}$ and $g^{22}_{011}$
which we denote by $h$ and $h_1$ respectively.
Plugging the solution into $g$, we see that $x^2$ divides $\det g$ when $h_1=0$.
This contradicts $\deg\Gamma\ge 5$,
hence we may set $h_1=1$. Then $E_{-2,2,3}$ yields $b_4=0$. Putting this into
$E_{1,3,1}$, $E_{2,2,1}$, $E_{2,2,2}$, $E_{2,3,1}$ we obtain a linear system with
constant coefficients for $c_3$, $c_4$, $b_5$, $b_6$ which yields
$$
  c_3=16 b_3 h-\tfrac32 b_3^2, \qquad
  b_5=-\tfrac83 b_3 h-\tfrac12 b_3^2,\qquad c_4=b_6=0.
$$
Plugging this solution into $E_{2,3,2}$ and $E_{3,2,2}$, we obtain a linear system
with constant coefficients for the unknowns $c_5$ and $b_7$ which yields
$$
  c_5 = b_3^3 - \tfrac{728}5 b_3^2 h + \tfrac{1648}{45}b_3h^2, \qquad
  b_7 = \tfrac12 b_3^3 + 40 b_3^2 h - \tfrac{80}9 b_3h^2.
$$
Putting this into $E_{3,2,3}$, we obtain the equation
$b_3h^2(3 b_3 - 2 h)=0$.
If $h=0$, then $\det g=0$. If $b_3=0$, then $x^2$ divides $\det g$.
Hence $b_3\ne 0$ and we may set $b_3=2$ by rescaling the parameter $t$.
Then $h=3$ and this gives us all coefficients of $g$.

Thus the curve $C$ is uniquely determined up to an affine linear change of variables.
It remains to observe that (\eqParamGen) has the required branch at $t=0$, and
to check that (\eqParamGen) gives a solution of the AlgDOP problem.
\qed\enddemo

\proclaim{ Lemma \lemBrFlex }
If $\deg\Gamma\ge 5$, then $\ord_t(\gamma)\ne(-1,2,3)$,
i.e., $\gamma$ cannot be of type $(1,3,4)$ (flex) with $\gamma\cdot P_\infty=1$.
\endproclaim

\demo{ Proof }
$n=45$.
The equations:
$E_{-6,3,i}$, $E_{-4,2,i}$, $E_{-3,3,i}$, $E_{-2,1,i}$, $E_{-2,3,i}$, 
$E_{-1,2,i}$, $E_{-1,3,i}$, $E_{0,2,i}$, $E_{0,3,i}$ ($i=1,2,3$);
$E_{-5,3,i}$, $E_{-3,2,i}$, $E_{-1,1,i}$, $E_{1,3,i}$, $E_{2,3,i}$, $E_{3,3,i}$ 
($i=2,3$); $E_{0,1,3}$, $E_{1,1,2}$, $E_{2,2,2}$, $E_{3,2,2}$.
The unknowns: $g^{ij}_{klm}$ ($1\le i\le j\le 2$),
$g^{13}_{klm}$ with $(k,l,m)\not\in\{200,110\}$, and
$g^{23}_{klm}$ with $(k,l,m)\not\in\{020,110,200\}$.
Plugging the solution into $E_{4,2,2}$ and $E_{4,3,1}$, we obtain the equations
$(58 b_4 + 23 c_1)g^{33}_{002} = 0$ and $(8 b_4 + c_1)g^{33}_{002} = 0$.
If $g^{33}_{002}=0$ or $b_4=c_1=0$, then $x^2$ divides $\det g$.
\qed\enddemo

\proclaim{ Lemma \lemBrFlat }
If $\deg\Gamma\ge 5$, then $\ord_t(\gamma)\ne(-1,1,3)$,
i.e., $\gamma$ cannot be of type $(1,2,4)$ (flat branch) with $\gamma\cdot P_\infty=1$.
\endproclaim

\demo{ Proof }
$n=42$.
The equations:
$E_{-3,\beta,i}$ (all $\beta$, $i$);
$E_{-7,3,i}$, $E_{-5,2,i}$, $E_{-5,3,i}$, $E_{-2,2,i}$, $E_{-1,2,i}$ ($i=1,2,3$);
$E_{-4,2,i}$, $E_{-2,1,i}$, $E_{-1,1,i}$, $E_{0,1,i}$ ($i=2,3$);
$E_{-1,3,i}$, $E_{0,2,i}$, $E_{1,3,i}$ ($i=1,2$);
$E_{1,1,2}$, $E_{1,2,2}$, $E_{2,1,2}$, $E_{2,2,1}$.
The unknowns: $g^{33}_{002}$,
$g^{33}_{001}$, and
$g^{ij}_{klm}$ with $(i,j)\ne(3,3)$ except
$g^{12}_{200}$, 
$g^{13}_{200}$, 
$g^{13}_{110}$,
$g^{22}_{200}$,
$g^{22}_{110}$,
$g^{23}_{200}$,
$g^{23}_{110}$,
$g^{23}_{011}$,
$g^{23}_{101}$,
$g^{23}_{100}$.
The solution implies that $x^2$ divides $\det g$.
\qed\enddemo

\proclaim{ Lemma \lemBrDblFlat }
If $\deg\Gamma\ge 5$, then $\ord_t(\gamma)\ne(-1,1,4)$,
i.e., $\gamma$ cannot be of type $(1,2,5)$ (doubly flat branch)
with $\gamma\cdot P_\infty=1$.
\endproclaim

\demo{ Proof }
$n=31$.
The equations:
$E_{-3,\beta,i}$ (all $\beta$, $i$);
$E_{-7,3,i}$, $E_{-5,2,i}$, $E_{-5,3,i}$, $E_{-2,3,i}$ ($i=1,2,3$);
$E_{-4,3,i}$, $E_{-2,2,i}$, $E_{-1,2,i}$, $E_{0,1,i}$ ($i=2,3$);
$E_{0,2,1}$, $E_{0,3,1}$.
The unknowns:
$g^{11}_{100}$, $g^{11}_{110}$,                                           
$g^{1j}_{101}$ ($j=1,2,3$),                                               
$g^{2j}_{000}$, $g^{2j}_{001}$, $g^{2j}_{011}$, $g^{2j}_{002}$ ($j=2,3$), 
and all the $g^{1j}_{0lm}$. 
The solution implies that $x^2$ divides $\det g$.
\qed\enddemo

\proclaim{ Lemma \lemBrGenTan }
If $\deg\Gamma\ge 5$, then $\ord_t(\gamma)\ne(-2,-1,1)$,
i.e., $\gamma$ cannot be of type $(1,2,3)$ (generic branch)
with $\gamma\cdot P_\infty=2$.
\endproclaim

\demo{ Proof }
$n=43$.
The equations:
$E_{-12,3,i}$, $E_{-11,3,i}$, $E_{-10,3,i}$,
$E_{-9,2,i}$, $E_{-9,3,i}$,
$E_{-8,2,i}$, $E_{-8,3,i}$,                                 
$E_{-6,1,i}$, $E_{-6,2,i}$, $E_{-6,3,i}$ ($i=1,2,3$),                           
$E_{-8,1,i}$, $E_{-7,1,i}$, $E_{-7,2,i}$, $E_{-7,3,i}$, $E_{-5,1,i}$ ($i=2,3$), 
$E_{-5,2,1}$, $E_{-5,3,1}$, $E_{-4,2,1}$.                                       
The unknowns:
 $g^{1j}_{klm}$ ($j=1,2,3$),
 $g^{j3}_{001}$, 
 $g^{j3}_{011}$, 
 $g^{j3}_{002}$, 
 $g^{j3}_{020}$  
 ($j=2,3$),
 $g^{22}_{001}$, 
 $g^{22}_{010}$, 
 $g^{22}_{101}$, 
 $g^{22}_{011}$, 
 $g^{22}_{002}$. 
We obtain that $x^2$ divides $\det g$.
\qed\enddemo

\proclaim{ Lemma \lemBrFlatTan }
If $\deg\Gamma\ge 5$, then $\ord_t(\gamma)\ne(-2,-1,2)$,
i.e., $\gamma$ cannot be of type $(1,2,4)$ (flat branch)
with $\gamma\cdot P_\infty=2$.
\endproclaim

\demo{ Proof }
$n=29$.
The equations:
$E_{-12,3,i}$, $E_{-11,3,i}$, $E_{-10,3,i}$, $E_{-9,2,i}$, $E_{-8,3,i}$, $E_{-8,2,i}$,
$E_{-7,3,i}$, $E_{-6,1,i}$ ($i=1,2,3$);
$E_{-9,3,i}$, $E_{-7,2,i}$ ($i=2,3$);
$E_{-5,2,1}$.
The unknowns:
$g^{1j}_{101}$, $g^{1j}_{110}$ ($j=1,2,3$),
$g^{2j}_{002}$, $g^{2j}_{011}$ ($j=2,3$),
$g^{11}_{100}$, and all the $g^{1j}_{0lm}$.
The solution implies that $x^2$ divides $\det g$.
\qed\enddemo

\proclaim{ Lemma \lemBrDblFlatTan }
If $\deg\Gamma\ge 5$, then $\ord_t(\gamma)\ne(-2,-1,3)$,
i.e., $\gamma$ cannot be of type $(1,2,5)$ (doubly flat branch)
with $\gamma\cdot P_\infty=2$.
\endproclaim

\demo{ Proof }
$n=43$.
The equations:
$E_{-12,3,i}$, $E_{-11,3,i}$, $E_{-10,3,i}$, $E_{-9,2,i}$, $E_{-8,2,i}$,
$E_{-7,3,i}$, $E_{-6,1,i}$, $E_{-6,3,i}$, $E_{-4,2,i}$ ($i=1,2,3$);
$E_{-8,3,i}$, $E_{-6,2,i}$, $E_{-5,2,i}$, $E_{-5,3,i}$, $E_{-4,1,i}$,
$E_{-3,\beta,i}$ ($i=2,3$).
The unknowns:
$g^{ij}_{001}$, $g^{ij}_{010}$, $g^{ij}_{011}$, $g^{ij}_{002}$ ($2\le i\le j\le 3$); 
$g^{i3}_{020}$, $g^{2i}_{101}$ ($i=2,3$); 
and all the $g^{1j}_{klm}$ with $k\ne2$.  
We obtain that $x^2$ divides $\det g$.
\qed\enddemo

\proclaim{ Lemma \lemBrCusp }
If $\deg\Gamma\ge 5$, then $\ord_t(\gamma)\ne(-2,1,2)$,
i.e., $\gamma$ cannot be of type $(2,3,4)$ (cusp)
with $\gamma\cdot P_\infty=2$.
\endproclaim

\demo{ Proof }
$n=43$.
The equations:
$E_{-4,\beta,i}$ (all $\beta$, $i$),
$E_{-10,3,i}$, $E_{-7,2,i}$, $E_{-7,3,i}$, $E_{-6,3,i}$ ($i=1,2,3$),
$E_{-9,3,i}$, $E_{-6,2,i}$, $E_{-5,2,i}$, $E_{-5,3,i}$, $E_{-3,1,i}$ ($i=2,3$), 
$E_{-3,2,i}$, $E_{-3,3,i}$ ($i=1,2$),
$E_{-2,\beta,2}$, $E_{-1,\beta,2}$ ($\beta=1,2,3$),
$E_{-2,3,1}$, $E_{0,2,2}$.
The unknowns: $g^{ij}_{klm}$ ($1\le i\le j\le 2$),
$g^{13}_{klm}$ with $(k,l,m)\not\in\{001,011,002\}$, and
$g^{23}_{klm}$ with $(k,l,m)\not\in\{100,110,020,200\}$.
Plugging the solution into $E_{-1,3,1}$, we obtain the equation
$c_{-1}g^{33}_{002}=0$. If $g^{33}_{002}=0$, then $\det g=0$.
Hence $c_{-1}=0$ and we may set $g^{33}_{002}=1$.
Putting this into $E_{0,3,1}$, we obtain $b_3=0$. Then
$x^2$ divides $\det g$.
\qed\enddemo

\proclaim{ Lemma \lemBrFlatCusp }
If $\deg\Gamma\ge 5$, then $\ord_t(\gamma)\ne(-2,1,3)$,
i.e., $\gamma$ cannot be of type $(2,3,5)$ (flat cusp)
with $\gamma\cdot P_\infty=2$.
\endproclaim

\demo{ Proof }
$n=36$.
The equations:
$E_{-4,\beta,i}$ (all $\beta$, $i$),                                  
$E_{-10,3,i}$, $E_{-7,2,i}$, $E_{-7,3,i}$, $E_{-5,3,i}$  ($i=1,2,3$); 
$E_{-8,3,i}$, $E_{-5,2,i}$,  $E_{-3,1,i}$, $E_{-2,1,i}$  ($i=2,3$);   
$E_{-2,2,i}$, $E_{-2,3,i}$                               ($i=1,2$);   
$E_{-3,3,2}$, $E_{-1,1,2}$, $E_{-1,2,2}$.                             
The unknowns:
$g^{1j}_{101}$                                               ($j=1,2,3$), 
$g^{1j}_{100}$, $g^{1j}_{110}$                                 ($j=1,2$), 
$g^{2j}_{000}$, $g^{2j}_{001}$, $g^{2j}_{011}$, $g^{2j}_{002}$ ($j=2,3$), 
$g^{22}_{010}$, $g^{22}_{020}$, $g^{22}_{101}$,                           
and all the $g^{1j}_{0lm}$.                                               
Plugging the solution into $E_{-3,2,3}$ and $E_{-1,3,2}$, we obtain the equations
$(40 b_2 - 7 c_{-1})g^{33}_{002} = 0$ and $(10 b_2 - c_{-1})g^{33}_{002} = 0$.
If $g^{33}_{002}=0$ or $b_4=c_{-1}=0$,
then $x^2$ divides $\det g$.
\qed\enddemo

\proclaim{ Lemma \lemBrCuspTan }
If $\deg\Gamma\ge 5$, then $\ord_t(\gamma)\ne(-3,-1,1)$,
i.e., $\gamma$ cannot be of type $(2,3,4)$ (cusp)
with $\gamma\cdot P_\infty=3$.
\endproclaim

\demo{ Proof }
$n=42$.
The equations:
$E_{-15,3,i}$, $E_{-13,3,i}$, $E_{-11,2,i}$, $E_{-11,3,i}$, $E_{-10,2,i}$, $E_{-9,2,i}$,
$E_{-9,3,i}$, $E_{-7,1,i}$ ($i=1,2,3$);
$E_{-12,2,i}$, $E_{-9,1,i}$, $E_{-8,1,i}$, $E_{-8,2,i}$ ($i=2,3$);
$E_{-7,2,i}$, $E_{-7,3,i}$, $E_{-6,2,i}$ ($i=1,2$);
$E_{-6,1,2}$, $E_{-5,1,2}$, $E_{-6,3,1}$, $E_{-4,2,1}$.
The unknowns:
$g^{ij}_{klm}$ with $1\le i\le j\le 2$ (except $g^{22}_{100}$ and $g^{22}_{200}$);
$g^{13}_{klm}$ (except $g^{13}_{100}$ and $g^{13}_{200}$);
$g^{23}_{001}$, $g^{23}_{002}$, $g^{23}_{020}$, $g^{23}_{011}$, $g^{23}_{101}$, 
$g^{33}_{011}$.
Plugging the solution into $E_{-6,3,1}$, we obtain
$c_{-2}g^{33}_{002}=0$. Then
$x^2$ divides $\det g$.
\qed\enddemo


\head\sectAlg. Tangent developables which admit solutions of the AlgDOP problem
\endhead

Let the notation be as in \S\sectLoc.
Thus $C$ is an irreducible curve in $\bC^3$ not lying in any plane,
and $\Gamma(x,y,z)=0$ is the equation of its tangent developable.

\proclaim{ Proposition \propAlgC } 
Suppose that there exists $g=(g^{ij})$ such that $(g,\Gamma)$
is a solution of the AlgDOP problem over $\bC$. Then $C$ admits one of the
following parametrizations in some affine coordinates in $\bC^3$:
\roster
\item"(i)"
       $t\mapsto(t,t^2,t^3)$;
\item"(ii)"
       $t\mapsto(t^{-1},t,t^2)$;
\item"(iii)"
       \hbox to 76mm{%
       $t\mapsto(t^2,\, 2t^3,\, 3t^4)$;\hfill} cusp at $t=0$;
\item"(iv)"
       \hbox to 76mm{%
       $t\mapsto(t^{-1} + t,\, 3t - t^3,\, 2t^2 - t^4)$ (cf.~Lemma~\lemBrGen);\hfill}
       cusps at $t=\pm1$;
\item"(v)"
       \hbox to 76mm{%
       $t\mapsto(3t-t^3,\,4t^2-2t^4,\,5t^3-3t^5)$;\hfill}
       cusps at $t=\pm1$;
%
\item"(vi)"
       \hbox to 76mm{%
       $\theta\mapsto(3\cos\theta+\cos 3\theta,\,
                      3\sin\theta-\sin 3\theta,\,6\cos 2\theta)$;\hfill}
       cusps at $\theta=0,\pi,\pm\pi/2$.
\endroster
They correspond respectively to Cases
\ctCubic\!, \ctCubic\!, \ctCQuar\!, \ctFQuin\!, \ctSQuin\!, \ctCSext\ of Table~\tabPie.
\endproclaim 

\demo{ Proof }
By Proposition~\propPie\ one of the cases in Table~\tabPie\ takes place.

If $\deg C=3$, then $C$ is the rational normal curve, i.e., it admits
a parametrization $t\mapsto(1:t:t^2:t^3)$ in some projective coordinates.
In these coordinates, $H_\infty$ is uniquely determined by the divisor $D$
which it cuts on $C$. Thus there are only three possibilities:
$D=3p_1$ (then Case (i) occurs);
$D=p_1+2p_2$ (then Case (ii) occurs);
$D=p_1+p_2+p_3$ and then there is no solution
(we compute $g$ by solving a linear system (\eqParam), and see that $\det g=0$).

Let $\deg C\ge 4$ and then $\deg\Gamma\ge 5$ (see Table~\tabPie).
We assume also that Case (iv) does not occur.
We see in Table~\tabPie\ that either $C$ or $\check C$ has degree at most $5$.
Hence each branch of $C$ has type $(m_1,m_2,m_3)$ with $m_3\le 5$
(i.e. contributes at most $2$ into $k_0+k_1+k_2$)
and $m_3=5$ (i.e. the contribution $2$) is possible in Case~\ctSQuin\ only.
Then Lemmas~\lemBrGen--\lemBrFlatCusp\ imply that $C$ does not have
any branch $\gamma$ such that $\gamma\cdot P_\infty=1$ or $2$.
Thus one of the following three cases occurs.

\smallskip
{\it Case 1.}
{\sl $\deg C=4$ and $C$ has a branch $\gamma$ such that $\gamma\cdot P_\infty=4$.}
Cases~\ctGQuar, \ctOFQua, and \ctTFQua\ of Table~\tabPie\ are impossible because
$C$ cannot have branches of type $(1,3,4)$ or $(1,2,4)$ (flex or flat branch)
in $\bC^3$ by Lemmas~\lemBrFlexAff\ and \lemBrFlatAff. In Case~\ctCQuar\ we
obtain (iii) by Lemma~\lemCompl.

\smallskip
{\it Case 2.}
{\sl $\deg C=5$ and $C$ has a branch $\gamma$ such that $\gamma\cdot P_\infty=5$.}
As above, Case~\ctFQuin\ of Table~\tabPie\ is impossible by 
Lemmas~\lemBrFlexAff\ and \lemBrFlatAff, thus Case~\ctSQuin\ takes place.
Then $\gamma$ is of the type $(3,4,5)$, $(2,3,5)$, or $(1,2,5)$
which corresponds to $(l_0,l_2)=(2,0)$, $(1,1)$, or $(0,2)$ respectively.
If there is an affine branch with $l_0+l_2=2$, i.e., of type
$(m_1,m_2,m_3)=(1,2,5)$, $(2,3,5)$, or $(3,4,5)$, then Lemma~\lemCompl\ implies that
$C$ is parametrized by
$t\mapsto(t^{m_1},t^{m_2},t^{m_3})$.
Solving the corresponding linear systems (\eqParam), we obtain that $\det g=0$
in all the three cases.

Thus $C$ has two affine branches $\gamma_1$ and $\gamma_2$, each contributing
$1$ to $k_0+k_2$.
By Lemma~\lemBrFlatAff, $C$ does not have any ordinary flat branch in $\bC^3$,
hence $\gamma$ is of type $(1,2,5)$, and $\gamma_1$, $\gamma_2$ are of type $(2,3,4)$.
Then the dual branches $\check\gamma$, $\check\gamma_1$, $\check\gamma_2$
of $\check C$ are of type $(3,4,5)$, $(1,2,4)$, $(1,2,4)$.
By Lemma~\lemCompl\ applied to $\check\gamma$ and $\check\gamma_1$,
the dual curve $\check C$ is parametrized by
$t\mapsto(1+t:t^3:t^4:t^5)$ in some homogeneous coordinates.
Thus $\check C$ (and hence $C$ as well) is
uniquely determined
up to automorphism of $\bP^3$.
The choice of $P_\infty$ is also unique because
it is the osculating plane at $\gamma$. It remains to check that the curve in (v)
gives a solution of the AlgDOP problem and its
non-generic branches are of the required types.

\smallskip
{\it Case 3.}
{\sl $\deg C=6$ and $C$ has branches $\gamma_1$, $\gamma_2$
such that $\gamma_1\cdot P_\infty=\gamma_2\cdot P_\infty=3$.}
Then $C$ is a 4-cuspidal sextic (see.~ Table~\tabPie).
Since $\deg\check C=4$, all cusps are ordinary. Then Lemma~\lemBrCuspTan\
implies that $\gamma_1$ and $\gamma_2$ are of type $(1,2,3)$ 
and $P_\infty$ is the the osculating plane for each of them.
Hence $\check C$ has a point with two local branches. Then Lemma~\lemNodal\
implies that $\check C$ is uniquely determined,
and we obtain (vi) by the same argument as in Case 2.
\qed\enddemo

\proclaim{ Proposition \propAlgR }
Suppose that $C$ is real and there exists $g=(g^{ij})$ such that $(g,\Gamma)$
is a solution of the AlgDOP problem over $\bR$. Then $C$ admits one of the
following parametrizations in some affine coordinates in $\bR^3$:
\rm{(i)--(vi)} of Proposition~\propAlgC\ or
\roster
\item"(iv${}'$)"
       \hbox to 60mm{%
       $t\mapsto(t^{-1} - t,\, 3t + t^3,\, 2t^2 + t^4)$;\hfill}
       cusps at $t=\pm i$;
\item"(v${}'$)"
       \hbox to 60mm{%
       $t\mapsto(3t+t^3,\,4t^2+2t^4,\,5t^3+3t^5)$;\hfill}
       cusps at $t=\pm i$;
\item"(vi${}'$)"
       \hbox to 60mm{%
       $t\mapsto(3t^{-1}+t^3,\,3t^{-2}+3t^2,\, t^{-3}+3t)$;\hfill}
       cusps at $t=\pm1,\pm i$;
\item"(vi${}''$)"
       \hbox to 60mm{%
       $t\mapsto(3t^{-1}-t^3,\,3t^{-2}-3t^2,\, t^{-3}-3t)$;\hfill}
       cusps at the roots of $t^4+1$.
\endroster
\endproclaim

\demo{ Proof }
It is easy to see that the real form of $C$ is determined by the involution
of complex conjugation on $\widetilde C$. The latter must preserve
the set of points of each type and the divisor $\nu^*(P_\infty)$ on $\widetilde C$.
This implies that the list is exhaustive.
A computation shows that all the cases are realizable.
\qed\enddemo

The tangent developable of the twisted cubic
$t\mapsto(t,t^2,t^3)$ is given by the equation $\Gamma_4=0$
where $\Gamma_4$ is the discriminant of
$P(u)=u^3+3xu^2+3yu+z$, i.e.,
$$
  \Gamma_4=3x^2y^2 - 4y^3 - 4x^3z + 6xyz - z^2.             \eqno(\eqTCubic)
$$

\proclaim{ Lemma \lemCubicI } {\rm(cf.~Prop.~\propAlgC(i))}
Let $(g,\Gamma)$ be a solution of the
AlgDOP problem over $\bR$ such that the surface $\Gamma=0$ contains the
tangent developable of the curve $C$ parametrized by
$t\mapsto(t,t^2,t^3)$. Then
$$
\split
  g=\;&a\left(\smallmatrix
     0         & 0        & 0       \\
     *         & 2(x^2-y) & 3(xy-z) \\ 
     *         & *        & 18(y^2-xz)
  \endsmallmatrix\right)
  +b\left(\smallmatrix
     1         & 2x    & 3y    \\
     *         & 4x^2  & 6xy   \\ 
     *         & *     & 9y^2
  \endsmallmatrix\right)
  +c\left(\smallmatrix
     x   &  2x^2   &  3xy    \\
     *   &  5xy-z  &  6y^2   \\ 
     *   &  *      &  9yz
  \endsmallmatrix\right)
\\
  &+d\left(\smallmatrix
     x  &  2y     &  3z   \\
     *  &  3xy+z  &  6xz  \\ 
     *  &  *      &  9yz
  \endsmallmatrix\right)
  +e\left(\smallmatrix
     x^2  &  2xy   &  3xz    \\
     *    &  4y^2  &  6yz   \\ 
     *    &  *     &  9z^2
  \endsmallmatrix\right)
  +f\left(\smallmatrix
     2(x^2-y)  &  xy-z       &  0  \\
     *         &  2(y^2-xz)  &  0  \\ 
     *         &  *          &  0
  \endsmallmatrix\right),
\endsplit
$$
$\det g=\Gamma_4\Gamma_2$ where $\Gamma_4$ is as in (\eqTCubic) and
$$
\smallmatrix
\Gamma_2=
a^2b + a^2cx - 2abcx + a^2dx + 2abdx - 2ac^2x^2 + a^2ex^2 + 
 2abex^2 + 2a^2fx^2 + 4abfx^2 + bc^2y - 2bcdy + 2ad^2y \\ + 
 bd^2y - 2abey - 2a^2fy - 2abfy + c^3xy - c^2dxy - 
 2acexy - bcexy + 2adexy + bdexy + 2acfxy - 
 2bcfxy \\ + 4adfxy + 2bdfxy - 2c^2fy^2 + 4aefy^2 + 
 2befy^2 + 2af^2y^2 + bf^2y^2 - cd^2z + d^3z + bcez - bdez - 
 2adfz + c^2exz \\ - 2cdexz + d^2exz + 2d^2fxz - 2aefxz - 
 2befxz - 2af^2xz - 2cefyz + 2defyz + cf^2yz + 
 df^2yz + ef^2z^2
\endsmallmatrix
$$
and one of the following cases occurs up to
affine linear change of coordinates:

\roster
\item"(i${}_1$)"
  $\Gamma=\Gamma_4$ and
  $(a,c-d,f)\ne(0,0,0)$, 
  $(b,c+d,e,af-d^2)\ne(0,0,0,0)$;
  in this case $\Gamma_2$ is a non-zero constant if and only if
  $c=d=e=f=0$ and $ab\ne 0$;
\smallskip
\item"(i${}_2$)"
  $b=d=f=0$, $(c,ae)\ne(0,0)$, and $\Gamma=x\Gamma_4$, then we have
  $\Gamma_2=x\Gamma_1$ where $\Gamma_1=a^2c + a(ae-2c^2)x + c(c^2-2ae)y + c^2ez$;
  in this case $\Gamma_1$ cannot be a nonzero constant, and we have $\Gamma_1=x$
  if and only if $c=0$ and $ae\ne 0$;
  
%
\smallskip
\item"(i${}_3$)"
 $a=b=c=0$, $(d,ef)\ne(0,0)$, and $\Gamma=z\Gamma_4$, then we have
  $\Gamma_2=z\Gamma_1$ where $\Gamma_1=d^3 + d^2(e + 2f)x + df(2e + f)y + ef^2z$;
  in this case $\Gamma_1$ is a nonzero constant if and only if
  $e=f=0$ and $d\ne 0$; we have $\Gamma_1=z$ if and only if $d=0$ and $ef\ne 0$;
%
%
%
%
\item"(i${}_4$)"
  $(a,\dots,f)=(0,\,0,\,0,\,1,\,-1,\,0)$,  $\,\Gamma=(x-1)z\Gamma_4$;  
%
\item"(i${}_5$)"
  $(a,\dots,f)=(1,1,0,0,-1,-1)$,
  $\Gamma=P(1)P(-1)\Gamma_4$,
  therefore $\{\Gamma_2=0\}$ is the union of two osculating planes of $C$;
  recall that $\Gamma_4=\text{\rm discr}_u P(u)$;
\item"(i${}_6$)"
   $(a,\dots,f)=(2\alpha,1,0,0,\pm1,0)$, $\alpha\ne0$,
   $\Gamma=(\alpha+1)x^2-y\pm\alpha$.
%
\item"(i${}_7$)"
  $(a,\dots,f)=(1,\;0,\;1,\;1,\;0,\;0)$, $\Gamma=(x-x^2+y)\Gamma_4$;
\endroster
\endproclaim

\demo{ Proof }
{\it Step 1.}
We find $g$ by solving the system of linear equations (\eqParam).
If $\Gamma=\Gamma_4$, we arrive to (i${}_1$) where the indicated condition
on $(a,\dots,f)$ is equivalent to $\det g\ne 0$.
We have
$$
   \Gamma_2(t,t^2,t^3)=(b + ct + dt + et^2)(a - ct + dt + ft^2)^2.
$$ 
Hence $\{\Gamma_2=0\}$ is disjoint from the curve $C$ (in $\bC^3$) if and only if
$c=d=e=f=0$, i.e., if and only if $\Gamma_2$ is a non-zero constant.

\medskip
{\it Step 2.}
The variable changes
$\varphi_\mu:(x,y,z)\mapsto(x,\,y+2\mu x,\,z+3\mu y+3\mu^2 x)$ and
$\psi_\lambda:(x,y,z)\mapsto(\lambda x,\lambda^2 y,\lambda^3 z)$
preserve $\{\Gamma_4=0\}$
and replace $(a,\dots,f)$ with
$$
  \big(\,a+\mu(c-d)+\mu^2 f,\, b-\mu(c+d)+\mu^2 e,\,
           c+\mu(f-e),\, d-\mu(f+e),\, e,\, f\,\big)            \eqno(\eqPhiMu)
$$
and $(\lambda^2 a,\lambda^2 b,\lambda c,\lambda d,e,f)$ respectively.
Thus, if $f\ne 0$ or $c-d\ne 0$, we may assume that $a=0$;
if $e\ne 0$ or $c+d\ne 0$, we may assume that $b=0$.

\medskip
{\it Step 3.} Here we suppose that 
$\Gamma=\Gamma_4\Gamma_1$ with $\deg\Gamma_1=1$. Any affine plane cuts the
curve $C$. Hence, up to affine change of coordinates, we may assume that the plane
$P=\{\Gamma_1=0\}$ passes through the origin.

\smallskip
{\it Case 3.1.} $P$ is transverse to $C$ at the origin, i.e., $P$ is parametrized by
$(t,u)\mapsto(x,y,z)=(At+Bu,t,u)$. Then (\eqParam) has three solutions:
\roster
\item"$\bullet$" $A=B=b=d=f=0$ (this is (i${}_2$));
\item"$\bullet$" $b=e=0$, $c=-d=aA$, $f=aA^2$, $B=-\tfrac13A^2$ (then $\det g=0$);
\item"$\bullet$" $b=f=B=0$, $c=d=-aA$, $e=2aA^2$.
\endroster
In the latter case we have $\Gamma_2=2a^3A(x-Ay)$, thus $A\ne 0$, and
the variable change $\varphi_\mu$ followed by $\psi_\lambda$, $\lambda=\mu=-1/(2A)$
(see Step~2) gives (i${}_6$) with $\alpha=-1$.

\smallskip
{\it Case 3.2.} $P$ has an ordinary tangency with $C$ at the origin, i.e.,
up to rescaling of the coordinates, $\Gamma_1=z-y$.
Then (\eqParam) does not have any non-zero solution.

\smallskip
{\it Case 3.3.} $P$ is the osculating plane of $C$ at the origin, i.e., $\Gamma_1=z$.
Then the only non-zero solution of (\eqParam) is (i${}_3$).

\medskip
{\it Step 4.} Suppose that $\deg\Gamma=6$ and $\Gamma_2=\Gamma_1\tilde\Gamma_1$ with
$\deg\Gamma_1=\deg\tilde\Gamma_1=1$.
According to the result of Step~3, each of the planes
$\{\Gamma_1=0\}$, $\{\tilde\Gamma_1=0\}$ is either an osculating plane for $C$
(Case (i${}_3$)) or a plane of the form $\{x=x_0\}$ (Case (i${}_2$)).
Any two distinct points on $C$ can be mapped to any fixed positions by an affine
linear automorphism which preserves $C$ (see Step 2). Thus $\Gamma_2$ is
as in (i${}_4$) or (i${}_5$) unless $\Gamma_2=x^2-1$ or $\Gamma_2=xz$.
In the latter two cases the system (\eqParam) does not have any nonzero solution.
Notice that this fact can be checked without any computations. Indeed, $\Gamma_2=x^2-1$
would imply that the $g^{i1}$ are divisible by $x^2-1$ and $\Gamma_2=xz$
would imply that the $g^{i1}$ (resp. $g^{i3}$) are divisible by $x$ (resp. by $z$).
It is immediately seen from the form of $g$ that this is impossible.
 
\medskip
{\it Step 5.} Suppose that $\deg\Gamma_2=2$ and $\Gamma_2$ is irreducible.

\smallskip
{\it Case 5.1.} $ef\ne0$. Then $\deg_z\Gamma_2=2$ and its coefficient of $z^2$
is $ef^2$ (a nonzero constant). By the result of Step~2 we may assume that $a=0$.
Then we compute the remainders of the division of
$g^{11}\partial_x\Gamma_1+g^{12}\partial_y\Gamma_1+g^{13}\partial_z\Gamma_1$
(viewed as a polynomial in $z$) by $\Gamma_2$ and equate its coefficients to zero
(see (A3) in \S1). The obtained system of equations has only two solutions with $ef\ne0$.
These are: (S1) $b=c=0$, $e=f$ and (S2) $b=c=d=0$. In both cases $\Gamma_2$
reducible.

\smallskip
{\it Case 5.2.} $ef=a=0$. Then $\Gamma_2=p_0(y,z) + p_1(y,z)x$.

\smallskip
{\it Case 5.2.1.} 
$p_1(y,z)=0$.
If $f=0$, then
$p_1=((d-c)((c^2 - be)y + (c - d)ez))$. If it is zero, then
$d=c$ (then $\Gamma_2=0$) or $c=e=0$ (then $\deg\Gamma_2<2$).
If $e=0$ and $f\ne0$, then
$p_1=(d-c)(c^2 - 2bf)y - 2d^2fz$. If it is zero, then $cd=0$ which implies that
$\Gamma_2=\Gamma_1^{(1)}\Gamma_1^{(2)}$, $\deg\Gamma_1^{(k)}\le 1$.

\smallskip
{\it Case 5.2.2.} $p_1(y,z)\ne 0$. Then we solve the system (\eqParam) for
the parametrization $(t,u)\mapsto\big(p_0(t,u)/p_1(t,u),\,t,\,u\big)$
of $\{\Gamma_2=0\}$. If $f=0$, the solutions are:
(S1) $d=c$; (S2) $b\ne0$, $e=cd/b$; (S3) $b=c=0$; (S4) $b=d=0$.
If $e=0$, the solutions are:
(S5) $c=d=0$; (S6) $c=f=0$; (S7) $d=0$, $f\ne0$, $b=c^2/(2f)$;
(S8) $f=d-c=0$; (S9) $d=f=0$.
In all these cases we have
$\Gamma_2=\Gamma_1^{(1)}\Gamma_1^{(2)}$, $\deg\Gamma_1^{(k)}\le 1$.

{\it Case 5.3.} $ef=0$ and $a\ne 0$.
By the result of Step~1 we know that $C$ and $\{\Gamma_2=0\}$ have a common point
in $\bC^3$. Suppose first that there is a real common point.
By an affine linear change of coordinates in $\bR^3$
we can achieve that this is the origin.
Since $\Gamma_2(0,0,0)=a^2b$, we then have $b=0$.
By the result of Step~2 we may assume that $f=0$ and $d=c$ (otherwise we reduce
to Case 5.2). Then $c\ne0$ because otherwise $\Gamma_2=x^2$.
Thus we have $b=f=0$ and $d=c\ne0$. One can check that this is a
solution of AlgDOP problem. If $e=0$ we obtain (i${}_7$) by the coordinate change
$\psi_\lambda$ (see Step~2) with $\lambda=c/a$.
If $e\ne 0$, the coordinate change $\varphi_\mu$ with $\mu=c/e$ followed by
$\psi_\lambda$ with $\lambda=e/c$, we obtain (i${}_6$) with
$(a,\dots,f)=(-ae/c^2,1,0,0,-1,0)$.
In the case when $C$ and $\{\Gamma_2=0\}$ do not have common real points,
one can show that $(a,\dots,f)=(2\alpha,1,0,0,1,0)$, $\alpha\in\bR$ (we omit the details).
In both cases we have $\alpha\ne0$ (otherwise $\Gamma_2=0$).
\qed\enddemo

The following lemma is a direct computation.

\proclaim{ Lemma \lemCubicII } {\rm(cf.~Prop.~\propAlgC(ii))}
Let $(g,\Gamma)$ be a solution of the
AlgDOP problem over $\bR$ such that the surface $\Gamma=0$ contains the
tangent developable of the curve $t\mapsto(t^{-1},t,t^2)$.
Then
$$
  g=a\left(\smallmatrix
  0  &  0        &   0          \\
  *  &  2(1-xy)  &  3(y - xz)   \\
  *  &  *        &  18(y^2 - z)
  \endsmallmatrix\right)
  +b\left(\smallmatrix
   x^2  & -xy   & -2xz \\
    *   &  y^2  &  2yz \\
    *   &   *   &  4z^2
  \endsmallmatrix\right), \qquad ab\ne0,
$$
$\det g=9 a^2 b\, x^2 ( 3 y^2 - 4 x y^3 - 4 z + 6 x y z - x^2 z^2)$.
The coordinate change $(x,y,z)\mapsto(\lambda^{-1}x,\lambda y,\lambda^2 z)$
preserves $\det g$ and transforms $(a,b)$ into $(a,\lambda^2 b)$.
Thus we can reduce to $(a,b)=(1,\pm1)$. \qed
\endproclaim

The tangent developable of the cuspidal quartic curve
$t\mapsto(t^2,2t^3,3t^4)$ is given by the equation $\Gamma_5=0$
where $\Gamma_5$ is the discriminant of $P(u)=u^4 - 6xu^2 - 4yu - z$, i.e.,
$$
  \Gamma_5=-54x^3y^2 + 27y^4 + 81x^4z - 54xy^2z + 18x^2z^2 + z^3.  \eqno(\eqTQuartic)
$$

\proclaim{ Lemma \lemQuartic } {\rm(cf.~Prop.~\propAlgC(iii))}
Let $(g,\Gamma)$ be a solution of the
AlgDOP problem over $\bR$ such that the surface $\Gamma=0$ contains the
tangent developable of the curve $t\mapsto(t^2,2t^3,3t^4)$,
i.e.,
$\Gamma_5$ a factor of $\Gamma$.
Then
$$
  g=a\left(\smallmatrix
    4x &  6y        &  8z  \\
     * &  3(9x^2+z) &  36xy \\ 
     * &  *         & 144(y^2-xz)
  \endsmallmatrix\right)
  +b\left(\smallmatrix
     4y        & 2(9x^2+z) & 24xy  \\
     *         & 36xy      & 36y^2 \\ 
     *         & *         & 48yz
  \endsmallmatrix\right)
  +c\left(\smallmatrix
    4x^2 & 6xy  & 8xz  \\
     *   & 9y^2 & 12yz \\
     *   &  *   & 16z^2
  \endsmallmatrix\right),
$$
$\det g = \Gamma_5\Gamma_1$ where
$\Gamma_1=3a^3 + 3(a^2c - 3ab^2)x + 3(b^3 - abc)y + b^2cz$,
and one of the following cases occurs up to rescaling of the coordinates:
\roster
\item"(iii${}_1$)"
        $(a,b)\ne(0,0)$, $\Gamma=\Gamma_5$, in this case $\Gamma_1$ is a nonzero
        constant if and only if $a\ne0$ and $b=c=0$;
\item"(iii${}_2$)"
        $(a,b,c)=(3,1,-1)$, $\Gamma=\Gamma_1\Gamma_5$,
        in this case $\{\Gamma_1=0\}$ is the osculating plane at $t=3$,
        and we have $\Gamma_1=P(3)$ (recall that $P(u)=u^4-6xu^2-4yu-z$);
\item"(iii${}_3$)"
        $(a,b,c)=(1,0,\pm1)$, $\Gamma=(x\pm1)\Gamma_5$.
\endroster
\endproclaim

\demo{ Proof }
We find $g$ by solving the linear system of equations (\eqParam). Then $\det g$ is
as stated. It vanishes identically if and only if $a=b=0$.
Since $\Gamma_5$ divides $\Gamma$ and $\Gamma$ divides $\det g$, we have
either $\Gamma=\Gamma_5$ or $\Gamma=\det g=\Gamma_5\Gamma_1$.
In the former case everything is done. So, we suppose that $\Gamma=\Gamma_5\Gamma_1$.

The change $(x,y,z)\mapsto(\lambda^2 x,\lambda^3 y,\lambda^4 z)$
transforms $(a,b,c)$ to $(a,\lambda b,\lambda^2 c)$.
Thus, if $abc\ne0$, we may assume that $(a,b)=(3,1)$.
Then the remainder of the division of
$g^{11}\partial_x\Gamma_1+g^{12}\partial_y\Gamma_1+g^{13}\partial_z\Gamma_1$
(viewed as a polynomial in $z$) by $\Gamma_1$ is equal to
$(c+1)q(x,y)$ where $q(x,y)$ is a polynomial in $x,y$ such that
$q(0,0)\ne 0$, and we arrive to solution (iii${}_2$).

If $abc=0$, we may rescale the coordinates so that each of
$a,b,c$ is $0$ or $\pm1$. In each case we check if (\eqParam) is satisfied.
\qed\enddemo

The following lemma is a direct computation based on Proposition~\propAlgR.
In Cases (v), (v${}'$) instead of $C$ we consider its image under
$(x,y,z)\mapsto(x,\frac32y,2z)$.

\proclaim{ Lemma \lemQuiSext } 
Let $(g,\Gamma)$ be a solution of the
AlgDOP problem over $\bR$ such that the surface $\Gamma=0$ contains the
tangent developable of the curves in Props.~\propAlgC(iv)--(vi) and
\propAlgR(iv${}'$)--(vi${}''$).
Then $\Gamma=\det g$ and one of the following cases takes place:
\roster
\item"(iv)" $g$ is given by (\eqIV) with $\eps=1$:
$$
\left(\smallmatrix
   x^2  &    3xy - 12    & 4xz-4y \\       
   *    & \;9y^2 - 12z\; & 12yz \\ 
   *    &      *         & 16z^2
\endsmallmatrix\right)
+\eps
\left(\smallmatrix
  -4 &   0            & 0           \\             
   * & \;72 - 24xy\;  & 24y - 36xz  \\ 
   * &   *            & 32y^2 - 144z
\endsmallmatrix\right);                       \eqno(\eqIV)
$$
\item"(iv${}'$)" $g$ is given by (\eqIV) with $\eps=-1$;
\item"(v)" $g$ is given by (\eqV) with $\eps=1$:
$$
 \left(\smallmatrix
  4y - 9x^2   & \;2z - 12xy & \;\;-15xz \\ 
  *           & \;-16y^2    & \;\;-20yz \\
  *           & \;*         & \;\;-25z^2
 \endsmallmatrix\right)
 + \eps
 \left(\smallmatrix
  24   &   32 x              &  40 y        \\ 
  *    & \;16(6 x^2 - 5 y)\; & 120(xy  - z) \\
  *    &   *                 & 400(y^2 - xz)
 \endsmallmatrix\right),                           \eqno(\eqV)
$$
$\frac14\det(5g)$ is the discriminant of $u^5 - 10u^3 - 10xu^2 - 5yu - z$;
\item"(v${}'$)" $g$ is given by (\eqV) with $\eps=-1$;
\item"(vi)"
$$
g=\left(\smallmatrix
  8 + y^2 + 4z - 2x^2 & -3xy &  12x - 2xz \\ 
  *                   & 8 -4z + x^2 - 2y^2 & -12y - 2yz  \\ 
  *                   & *                   & 16 + 8x^2 + 8y^2 - 4z^2
\endsmallmatrix\right),
$$
$\tfrac14\det g$ is the discriminant of $u^4 - su^3 + zu^2 - \bar su + 1$,
$s=x+iy$;
\item"(vi${}'$)" $g$ is given by (\eqVI) with $\eps=-1$:
$$
\left(\smallmatrix
  3x^2 - 8y   &   2xy - 12z    &    xz \\ 
   *          &   4y^2 - 8xz   &   2yz \\
   *          &    *           &   3z^2
\endsmallmatrix\right)
+\eps
\left(\smallmatrix
   0    &    0    &   16 \\
   0    &   16    &   12x \\
  16    &   12x   &    8y
\endsmallmatrix\right);                                 \eqno(\eqVI)
$$
\item"(vi${}''$)" $g$ is given by (\eqVI) with $\eps=1$.
\endroster
\endproclaim


\newpage

\head\sectDOP. Solutions of SDOP problem bounded by tangent developables
\endhead

\subhead\sectDOP.1. Bounded solutions
\endsubhead

\proclaim{ Theorem \thB }
Up to affine linear change of coordinates,
the following is a complete list of solutions $(\Omega,g,\rho)$
of DOP problem in $\bR^3$ such that $\Omega$ is a bounded domain
whose boundary $\partial\Omega$ contains a piece of a tangent developable surface.
In each case $g$ is as in the corresponding case of Lemmas~\lemCubicI, \lemQuartic,
\lemQuiSext\ (sometimes with additional restrictions on the parameters)
and $\Omega$ is the only bounded component of the complement of $\{\det g=0\}$;
see Figures~\figAB--\figABD\ and comments on them in \S\sectFig.
\roster
\item"(i${}_4$)" $\rho=\Gamma_4^{p-1} z^{q-1} (1-x)^{r-1}$,
                 $6p>1$, $q>0$, $r>0$, $2p+q>1$;

\item"(i${}_5$)" $\rho=\Gamma_4^{p-1} P(1)^{q-1} P(-1)^{r-1}$,
                 $6p>1$, $q>0$, $r>0$, $2p+q>1$, $2p+r>1$;

\item"(i${}_6$)" $\alpha>0$ and $e=-1$,
       $\rho=\Gamma_4^{p-1}\Gamma_2^{q-1}$, $6p>1$, $q>0$ (see Remark~\remCubicI);

\item"(iii${}_2$)" $\rho=\Gamma_5^{p-1}\Gamma_1^{q-1}$, $4p>1$, $q>0$, $2p+q>1$;

\item"(iii${}_3$)" $c=-1$, $\rho=\Gamma_5^{p-1}(1-x)^{q-1}$, $4p>1$, $q>0$;

\item"(v,vi)" $\rho=(\det g)^{p-1}$, $4p>1$;
\endroster
\endproclaim

\midinsert
\centerline{
   \epsfxsize=42mm\epsfbox{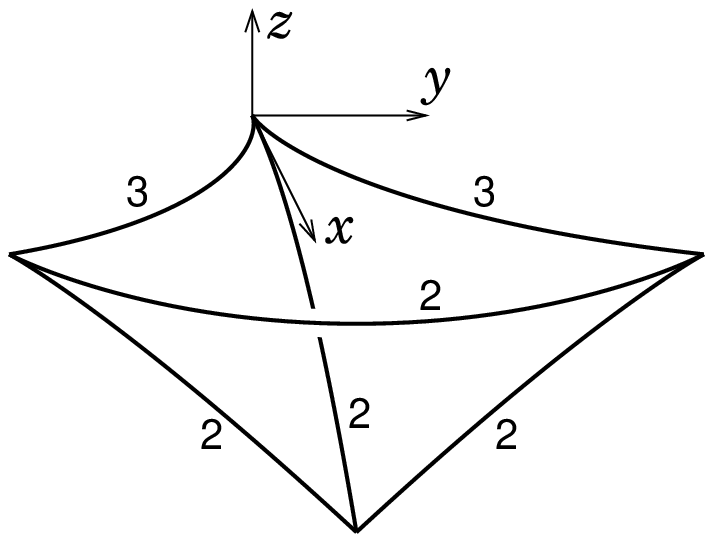}
   \hskip 10 mm
   \epsfxsize=32mm\epsfbox{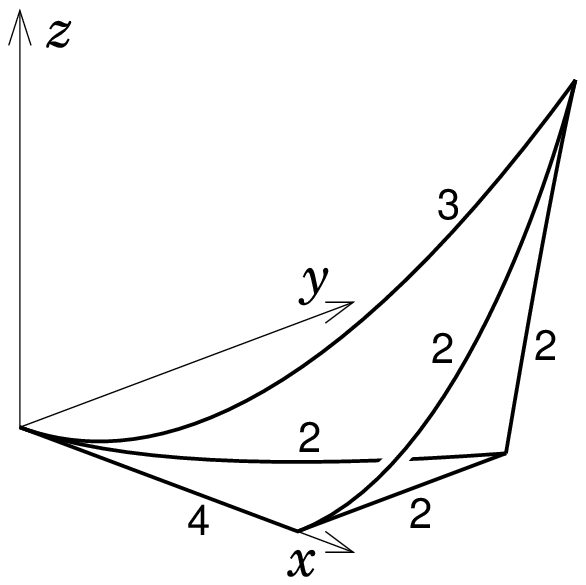}
}
\botcaption{ Figure \figAB }
(iii${}_3$) and (i${}_4$): the quotients of $\bS^3$
by the reflection groups $A_1+A_3$ (the truncated swallow tail) and $A_1+B_3$.
\endcaption
\endinsert

\midinsert
\centerline{
   \lower-4mm\hbox{\epsfxsize=38mm\epsfbox{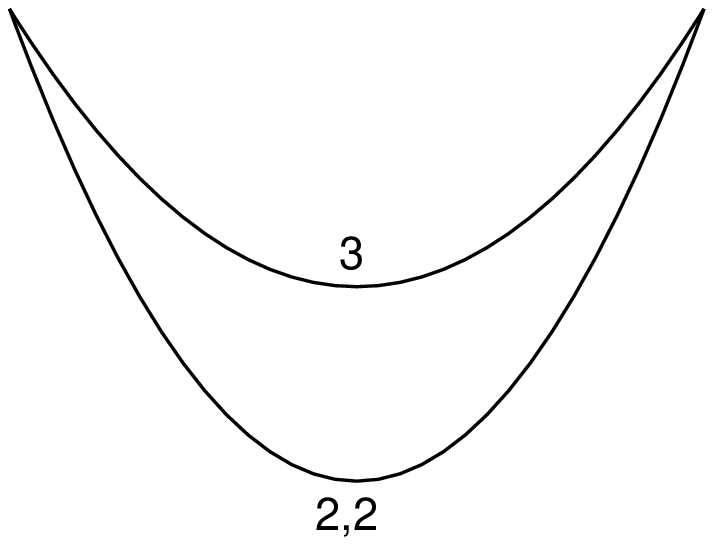}}
   \hskip 12mm
   \epsfxsize=44mm\epsfbox{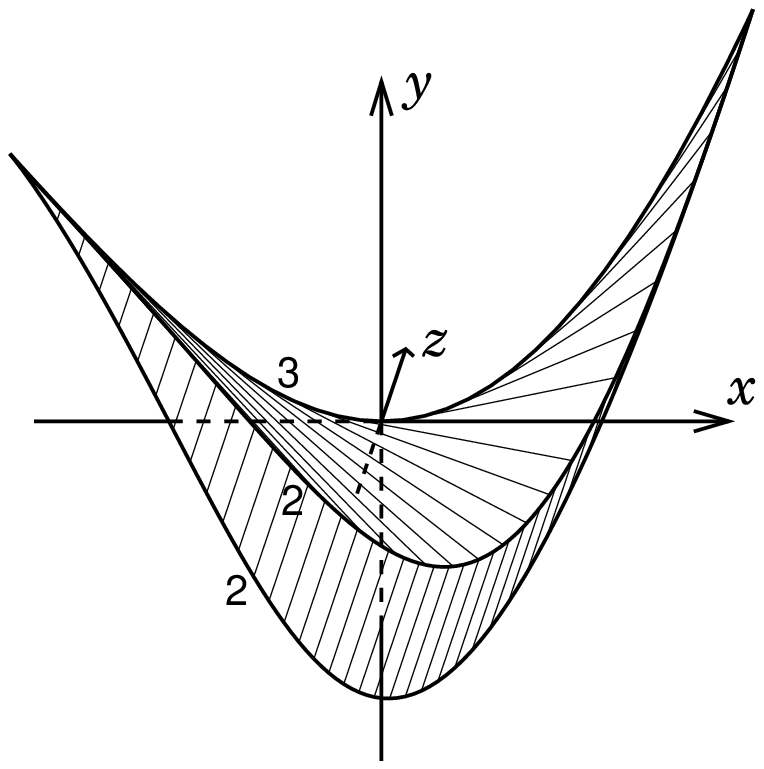}
}
\botcaption{ Figure \figAiAii } (i${}_6$):
 the quotient of $\bS^3$ by the reflection group $A_1+A_2$ \hbox to 5mm{}
 (the projection on the $xy$-plane is on the left hand side).
\endcaption
\endinsert

\midinsert
\centerline{\hskip0pt
   \lower-11mm\hbox{\epsfxsize=10mm\epsfbox{a3aff1.eps}}
   \hskip 15mm
   \lower-4mm\hbox{\epsfxsize=31mm\epsfbox{a3aff0.eps}}
   \hskip 15mm
   \epsfxsize=36mm\epsfbox{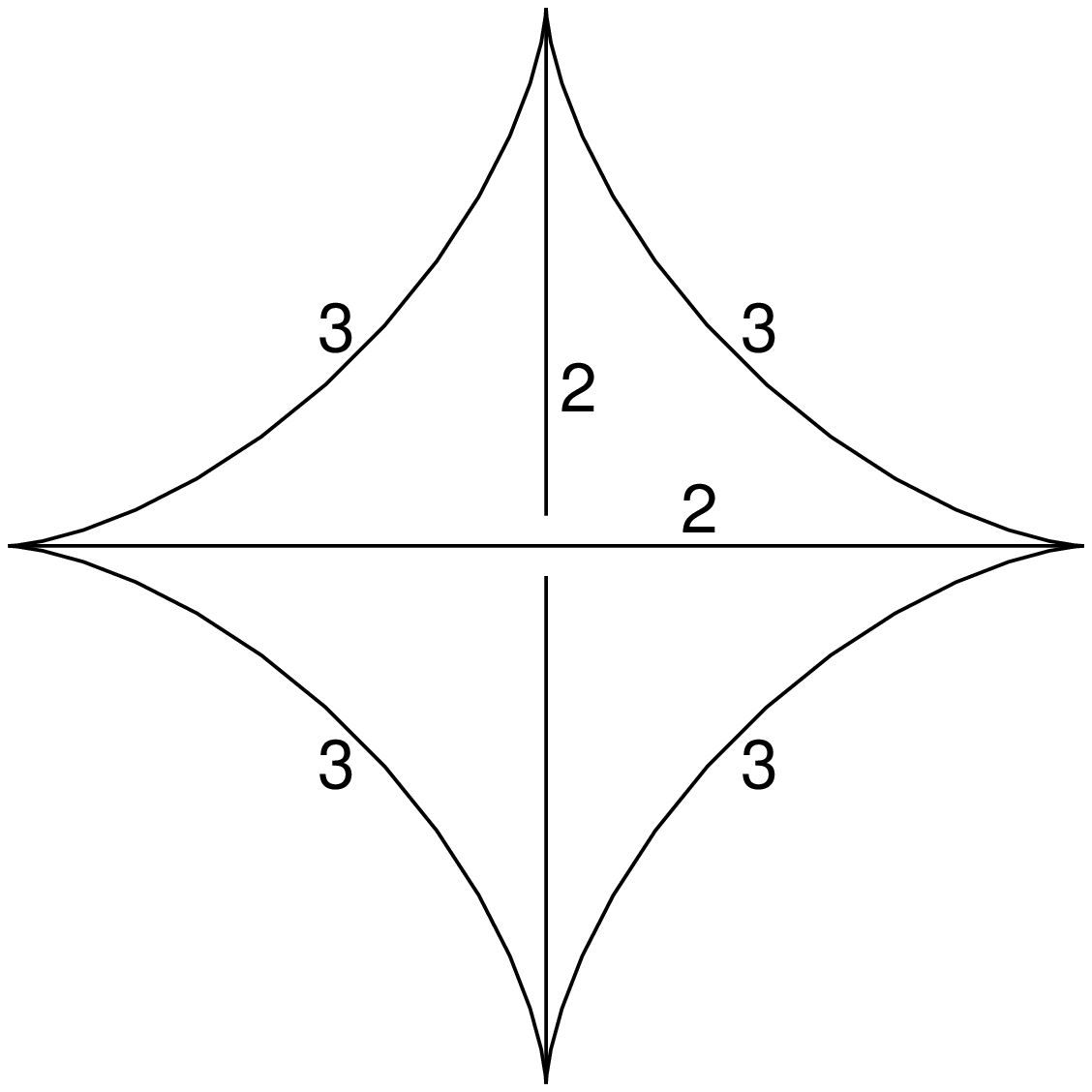}
}
\botcaption{ Figure \figAaff }
 (vi): the quotient of $\bR^3$ by the affine reflection group $\widetilde A_3$.
\endcaption
\endinsert

\midinsert
\centerline{\hskip0pt
   \lower-4mm\hbox{\epsfxsize=6mm\epsfbox{c3aff1.eps}}
   \hskip 12mm
   \lower-3mm\hbox{\epsfxsize=30mm\epsfbox{c3aff0.eps}}
   \hskip 15mm
   \epsfxsize=46mm\epsfbox{c3aff.eps}
}
\botcaption{ Figure \figCaff} (i${}_5$):
 the quotient of $\bR^3$ by the affine reflection group $\widetilde C_3$.
 The faces $ABC$ and $BCD$ are on the osculating planes at $A$ and $D$
resp.
\endcaption
\endinsert

\midinsert
\centerline{
   \lower-8mm\hbox{\epsfxsize=44mm\epsfbox{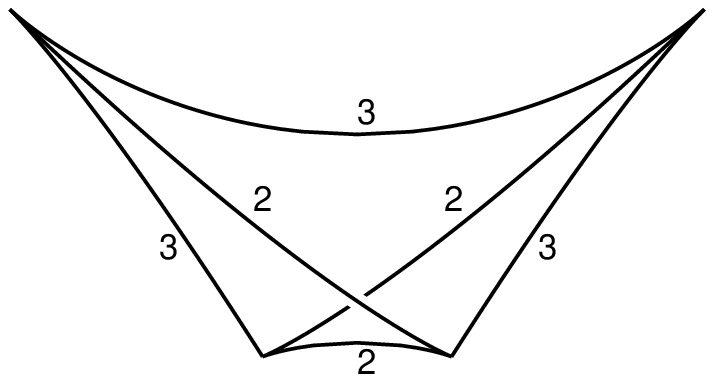}}%
   \hskip 2mm
   \lower-2mm\hbox{\epsfxsize=38mm\epsfbox{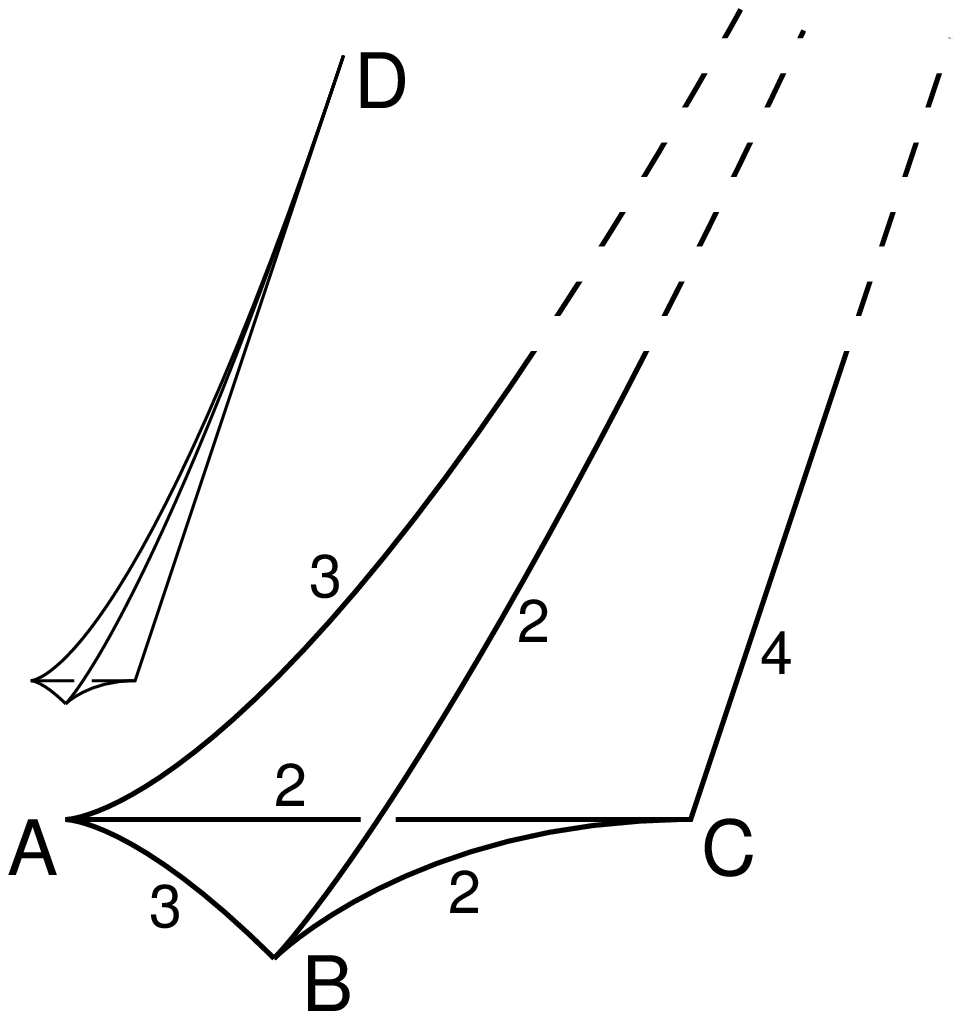}}%
   \epsfxsize=40mm\epsfbox{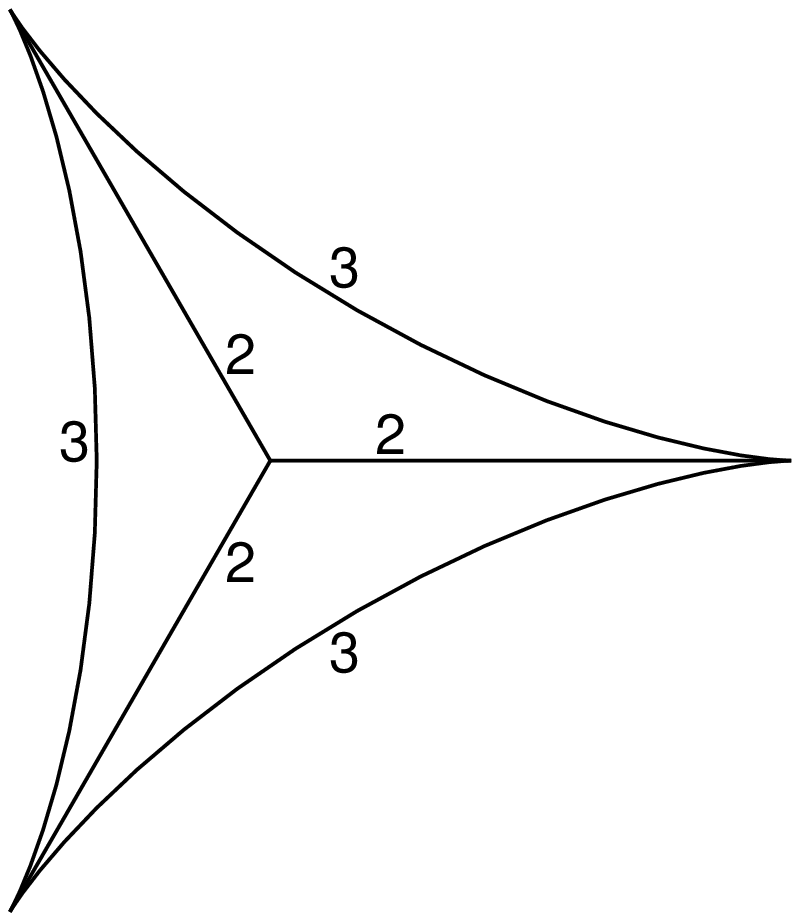}
}
\botcaption{ Figure \figABD } (v), (iii${}_2$), \S\sectD:
the quotients of $\bS^3$ by the reflection groups
 $A_4$, $B_4$, $D_4$.
The face $BCD$ belongs to the osculating plane at $D$.
\endcaption
\endinsert

\demo{ Proof }
{\it Boundedness of $\Omega$.}
Let us show that $\bR^3\setminus\Sigma$ where $\Sigma=\{\Gamma=0\}$
does not have any bounded component
in all other cases of Lemmas~\lemCubicI--\lemQuiSext.
For (i${}_1$), (i${}_2$), (i${}_3$), (iii${}_1$) this fact is evident
because $\Gamma$ is quasihomogeneous.
In other cases we consider the projection $\pi:(x,y,z)\mapsto(x,y)$ and find
the regions on the $xy$-plane over which $\Sigma$
is a disjoint union of graphs of smooth functions
(i.e. over which $\pi|_\Sigma$ is a covering).
This is the complement of the real curve
$R=\{D_z(x,y)C_z(x,y)=0\}$ where
$D_z$ is the discriminant of $\Gamma$ with respect to $z$,
and $C_z$ is the coefficient of the highest power of $z$ in $\Gamma$.
This curve is depicted in Figure~\figIrrel\ in the respective cases.
The dashed line represents $\pi(B_-)$ where $B_-$ is the part of the curve
of self-intersection of (the complexification of) $\Sigma$ such that
two non-real local branches of $\Sigma$ cross at points of $B_-$
(that is $\Sigma$ has the equation $u^2+v^2=0$ in some local real analytic coordinates $(u,v,w)$ near each point of $B_-$).
We see in Figure~\figIrrel\ that
all components of $\bR^2\setminus(R\setminus\pi(B_-))$ are unbounded.
Hence so are all components of $\bR^3\setminus\Sigma$.
It remains to exclude Case (i${}_6$) when $\alpha<0$ or $e=1$.
In this case we have $D_z=y-x^2$ and $C_z=(\alpha+1)x^2-y+e\alpha$, $e=\pm1$.
If $e=1$, then all components of $\bR^2\setminus R$ are unbounded.
If $e=-1$ and $\alpha<0$, then there is a bounded component $\Omega$, but
$\pi^{-1}(\Omega)\cap\Sigma$ is empty.

\midinsert
\centerline{%
   \epsfxsize=17mm\epsfbox{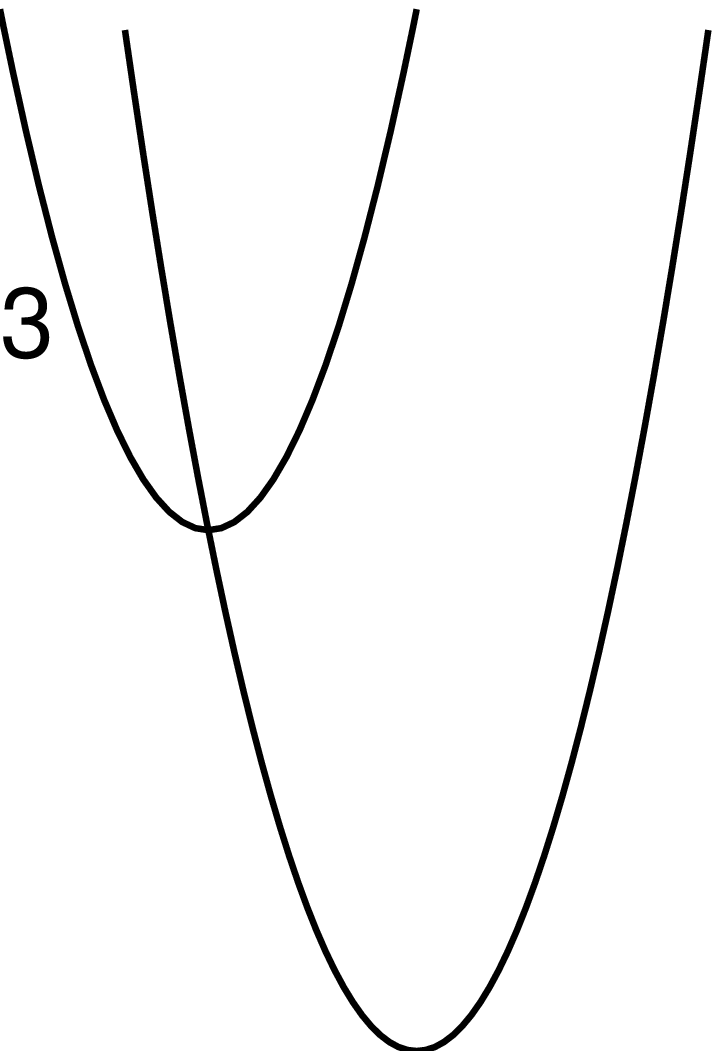}
   \hskip 3mm
   \epsfxsize=23mm\epsfbox{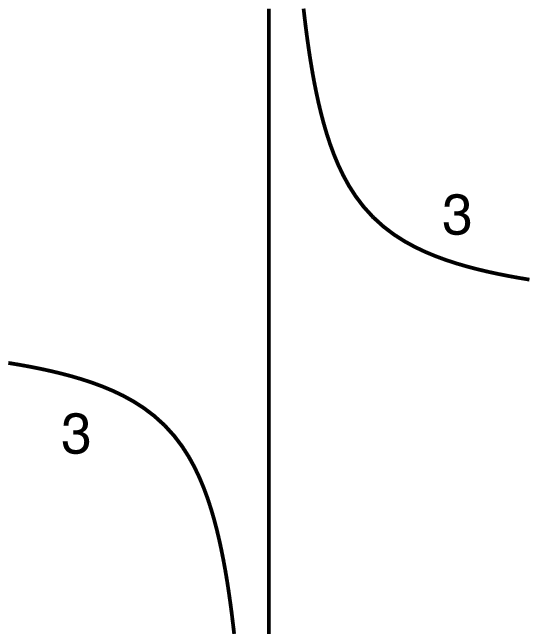}
   \hskip 3mm
   \epsfxsize=45mm\epsfbox{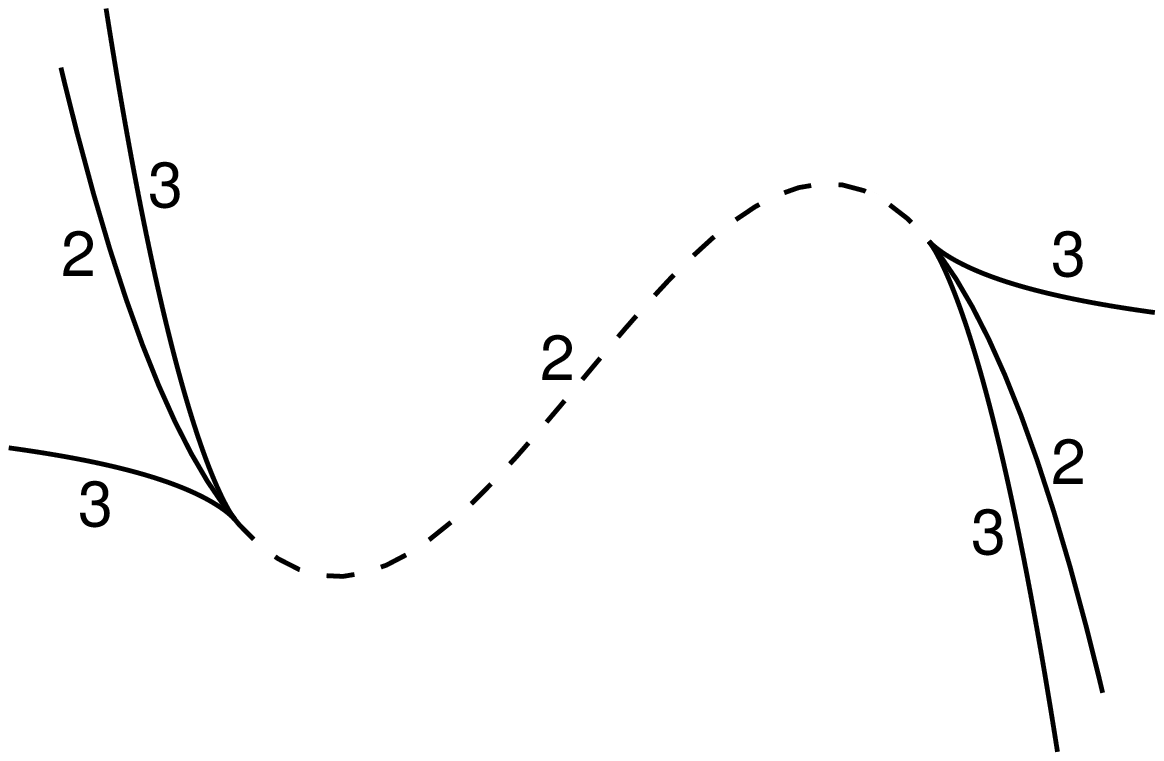}
   \hskip 3mm
   \epsfxsize=22mm\epsfbox{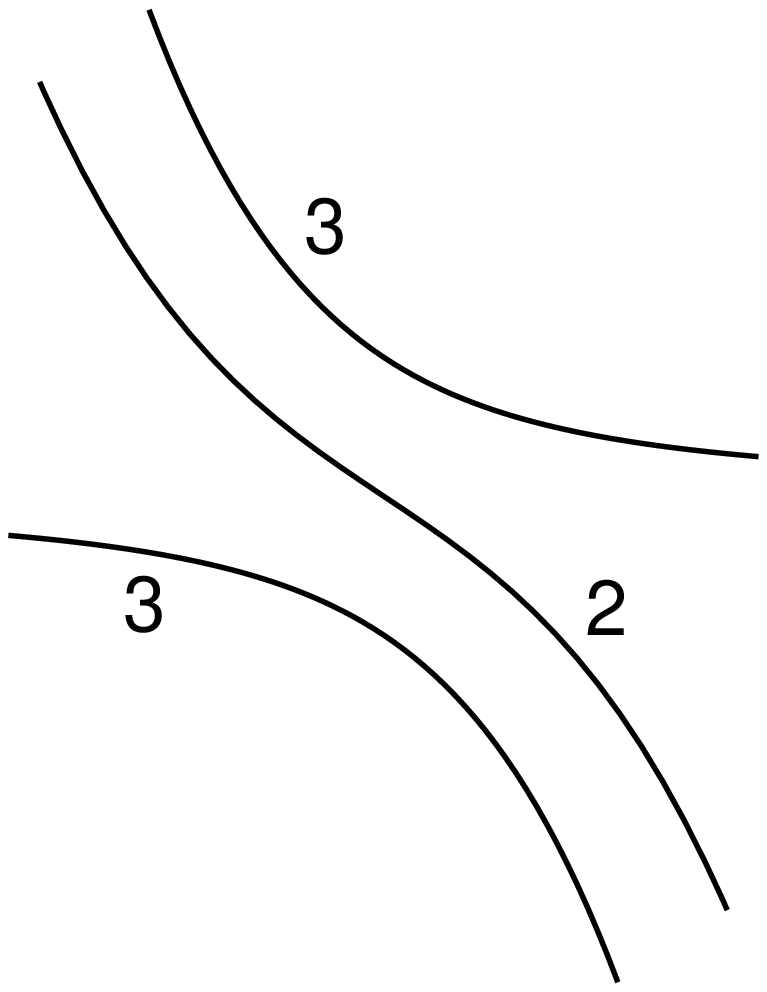}
}
\vskip-15pt
\centerline{%
\hskip 3mm
\hbox to 22mm{(i${}_7$)}
\hbox to 26mm{(ii)}
\hbox to 48mm{(iv)}
\hbox to 22mm{(iv${}'$)}
}
\medskip
\centerline{%
   \epsfxsize=25mm\epsfbox{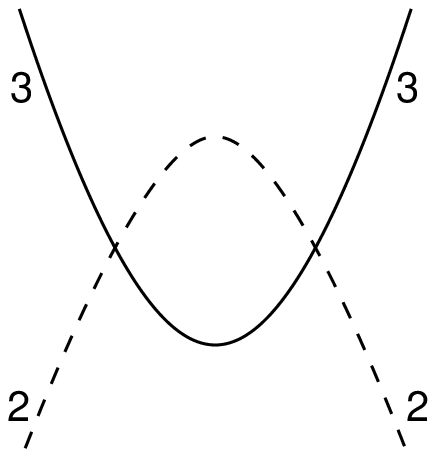}
   \hskip5mm
   \epsfxsize=42mm\epsfbox{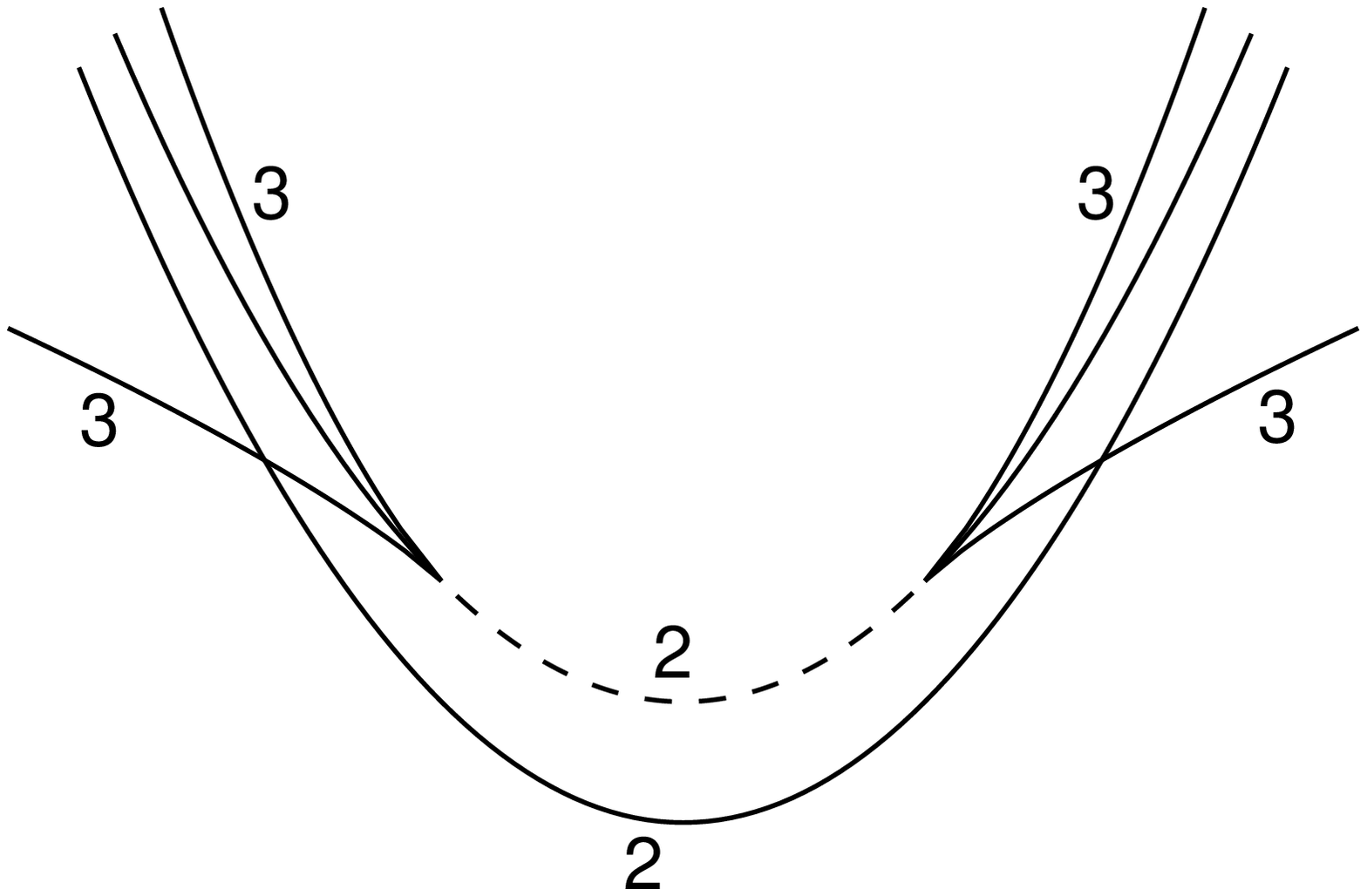}
   \hskip5mm
   \epsfxsize=24mm\epsfbox{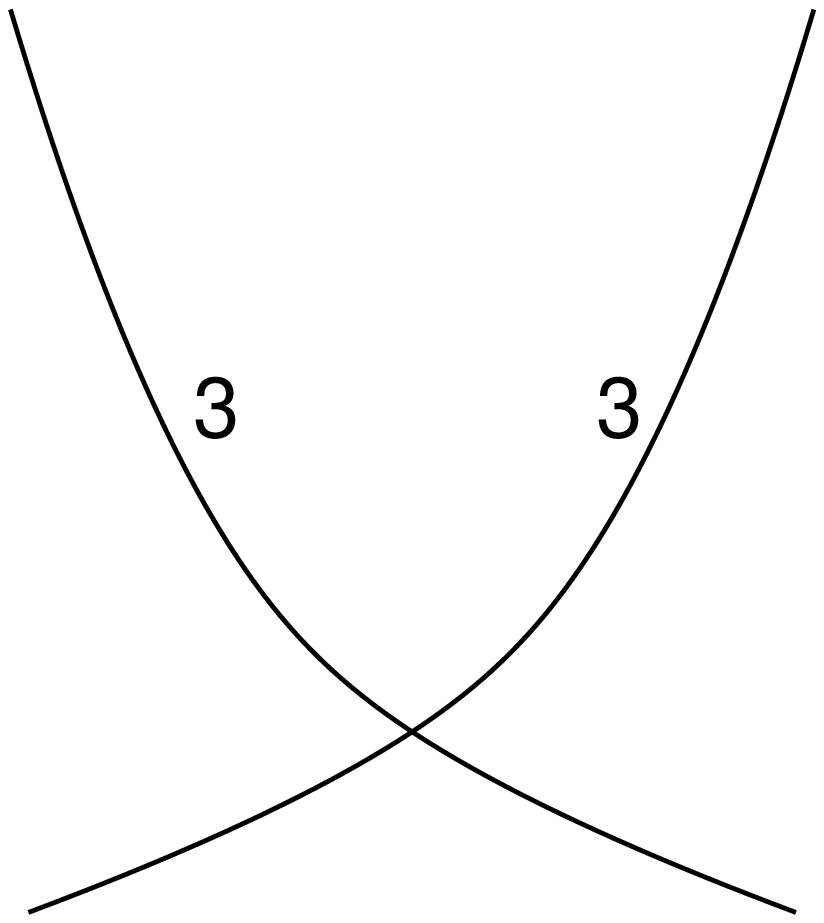}
}
\vskip-8pt
\centerline{%
\hskip 20pt
\hbox to 30mm{(v${}'$)}
\hbox to 47mm{(vi${}'$)}
\hbox to 24mm{(vi${}''$)}
}
\botcaption{ Figure \figIrrel }
Irrelevant solutions of the AlgDOP problem.
\endcaption
\endinsert

\smallskip
{\it Integrability of $\rho$.} In Case (i${}_4$), the integrability
conditions at the origin, at cuspidal edge, and at the $x$-axes (the line of
tangency) are, respectively, $2p+q>1$, $p>1/6$, and $p+q>1/2$ but the last
condition follows from the first one because $q>0$.
Let us prove the first condition. Let $\Omega_1=\{(y,z)\mid (1,y,z)\in\Omega\}$.
We have $\Gamma_4(x,y,z)=x^6\Gamma_4(1,y/x^2,z/x^3)$,
hence, using the variable change $y=x^2\eta$, $z=x^3\zeta$, we obtain
$$
   \int_{\Omega\cap\{x<\eps\}} z^{q-1}\Gamma_4^{p-1}\,dx\,dy\,dz
   = \int_0^\eps dx\int_{\Omega_1}(x^3\zeta)^{q-1}
        \big(x^6\Gamma_4(1,\eta,\zeta)\big)^{p-1} x^5 d\eta\,d\zeta,
$$
which is finite if and only if $6(p-1)+3(q-1)+5>-1$ (i.e., $2p+q>1$) and
$\int_{\Omega_1}\zeta^{q-1}\Gamma_4(1,\eta,\zeta)^{p-1}d\eta\,d\zeta$ is finite.
The integrability conditions in dimension 2
are obtained in the same way (in [\refBOZ, Remark~2.28] they are stated as an
evident fact). In our case they are $6p>1$ and $2p+2q>1$.
The same or similar arguments work in Cases (i${}_5$), (i${}_6$), (iii${}_3$) as well.
In the remaining three cases
the surface is not quasihomogeneous, however, one can show that
the restrictions are the same as in the quasihomogeneous case (we omit the proof).
\qed\enddemo

\subhead\sectDOP.2. Unbounded solutions
\endsubhead

\proclaim{ Theorem \thU }
Up to affine linear change of coordinates,
the following is a complete list of solutions $(\Omega,g,\rho)$
of SDOP problem in $\bR^3$ such that $\Omega$ is an unbounded domain
whose boundary $\partial\Omega$ contains a piece of a tangent developable surface.
In each case $g$ is as in the corresponding case of Lemmas~\lemCubicI, \lemQuartic\
with additional restrictions on the parameters.
\roster
\item"(i${}_1$)" $(a,\dots,f)=(2\alpha,1,0,0,0,0)$, $\alpha>0$;
                   $\Omega$ is the component of $\bR^3\setminus\{\det g=0\}$
                   containing $(0,-1,0)$ (i.e., the domain in Figure~\figAiAii\
                   is $\Omega\cap\{y\ge(\alpha+1)x^2-\alpha\}$),
                 $\rho=\Gamma_4^{p-1}\exp\big(\lambda y-\lambda(1+\alpha)x^2\big)$,
                 $p>1/6$, $\lambda>0$;

\item"(i${}_3$)" $(a,\dots,f)=(0,0,0,1,0,0)$;
                   $\Omega$ is the only component of $\bR^3\setminus\{\det g=0\}$
                   such that $\Omega\cap\{x=1\}$ is bounded
                   (i.e., the right domain in Figure~\figAB\ is $\Omega\cap\{x\le 1\}$),
                   $\rho=\Gamma_4^{p-1}z^{q-1}\exp(-\lambda x)$,
                   $p>1/6$, $q>0$, $2p+q>1$, $\lambda>0$;

\item"(iii${}_1$)" $(a,b,c)=(1,0,0)$; 
                   $\Omega$ is the only component of $\bR^3\setminus\{\det g=0\}$
                   such that $\Omega\cap\{x=1\}$ is bounded
                   (i.e., the left domain in Figure~\figAB\ is $\Omega\cap\{x\le 1\}$),
                   $\rho=\Gamma_5^{p-1}\exp(-\lambda x)$, $p>1/4$,
                   $\lambda>0$.
\endroster
\endproclaim

\demo{ Proof }
Let $\Gamma$ be the minimal polynomial vanishing on $\partial\Omega$. Then
$(g,\Gamma)$ is a solution of the AlgDOP problem. Set $\Delta=\det g$ and
$\Sigma=\{\Gamma=0\}$.
If $\Omega$ is unbounded and $\Delta$ does not have multiple components, then
$\deg\Delta<6$. This fact excludes all the cases
of Lemmas~\lemCubicI--\lemQuiSext\ for $(g,\Gamma)$
except those considered below.

\smallskip
(i${}_1$). Then $\partial\Omega\subset\Sigma_4=\{\Gamma_4=0\}$ (recall that $\Gamma_4$
is given in (\eqTCubic)). Let $\pi$ be the projection $\bR^3\to\bR^2$,
$(x,y,z)\mapsto(x,y)$. Then $\Sigma_4$ cuts $\bR^3$ into two unbounded components
$\Omega_+$, $\Omega_-$. One of them (let it be $\Omega_+$) is projected by $\pi$
onto the non-convex set $\{y<x^2\}$, and $\pi^{-1}(\pi(p))$ is a finite interval
for any $p\in\Omega_+$ (see Figure~\figAiAii).
Since $\Omega$ is one of $\Omega_+$, $\Omega_-$,
we have $\partial\Omega=\Sigma_4$.

Let us show that any affine plane $P$ intersects each of $\Omega_+$, $\Omega_-$.
The curve $C$ (whose tangent developable is $\Sigma_4$) has only one
point at infinity: the infinite point of the $z$-axis. If the projective closure of
$P$ does not pass through this point, then $P$ cuts $C$ in some finite
real point because the degree of $C$ is odd.
Otherwise $P=\pi^{-1}(L)$ where $L$ is a line in $\bR^2$,
hence $P$ cuts $\Sigma$ because $\pi(\Sigma)=\{y<x^2\}$.
In both cases $P$ cuts each of $\Omega_\pm$.
This fact implies that $\Delta=\Gamma_4$ (up to a scalar factor).
Indeed, recall that either $\Delta$ has a multiple component, or $\deg\Delta<6$.
Thus, if $\deg\Delta>4$, then in both cases $\Delta$
would vanishes on some plane $P$. This is impossible because $\Delta|_\Omega\ne0$
and $P\cap\Omega\ne\varnothing$.

By solving the system of linear equations (\eqParam),
we obtain the required form of $\rho$,
maybe, multiplied by $e^{\lambda_1 x}$, however, this factor can be killed by
the transformation $\varphi_\mu$ with a suitable $\mu$; see (\eqPhiMu).
We have $\Omega=\Omega_+$ because $\Omega_-$ contains cylinders parallel to the $z$-axis,
which contradicts the integrability condition for the measure of this form.
The positive definiteness of $g$ in $\Omega$ implies that $a>0$ and $b>0$. Then we
may set $b=1$, $a=2\alpha$, $\alpha>0$. The integrability conditions
near $C$ and at the infinity are, respectively, $p>1/6$ (see [\refBOZ, Remark~2.28]) 
and $\lambda>0$.

\smallskip
(i${}_2$). Then $\deg\Delta=6$, and $\Delta$ has multiple factors of $\Delta$
if and only if and $b=c=d=f=0$, $ae\ne0$. In this case $\Delta=x^2\Gamma_4$.
No exponential factor of $\rho$.

\smallskip
(i${}_3$), $a=b=c=e=f=0$, $d\ne0$. The solution is antisymmetric under the rotation
$(x,y,z)\mapsto(-x,y,-z)$, hence we may set $d=1$.
Then $g$ is positive definite only in the indicated domain. Solving the linear
equations, we obtain that the measure is of the required form.
The integrability condition at the infinity is $\lambda>0$. The others are
the same as in Theorem~\thB(i${}_4$).

\smallskip
(i${}_3$), $a=b=c=d=0$, $ef\ne0$. Solving the linear equations, we obtain that
$\rho$ has an exponential factor only when $e=f$, and it is $\exp(\lambda y/z)$.
Since $\Gamma_4(\lambda x,\lambda^2 y,\lambda^3 z)=\lambda^6 \Gamma_4(x,y,z)$, using the variable change $y_1=y/x^2$, $z_1=z/x^3$, one can easily show
that the integrability condition fails
for any choice of $\Omega$.

\smallskip
(i${}_6$), $\alpha=-1$. No exponential factor of $\rho$.

\smallskip
(iii${}_1$), $b=c=0$. Straightforward; see the bound for $p$ in
Theorem~\thB(iii${}_3$).
\qed\enddemo

\medskip\noindent{\bf Remark \remCubicI.}
The solutions (i${}_1$) and (i${}_6$) with different values of the parameter $\alpha$
(and in the latter case even the underlying domains)
cannot be transformed to each other by any affine linear transformation.
However, these solutions are also solutions of the weighted DOP problem
(see [\refBB], [\refOr]) with weights $(1,2,3)$ and the $(1,2,3)$-admissible
change of variables (see the definition in [\refOr, \S2.2])
\if01{
$$
  (x,y,z)\mapsto\big(\,x,\;\lambda^2 y + (1-\lambda^2)x^2,\;
      \lambda^3 z + 3(\lambda^2-\lambda^3)xy + (1-3\lambda^2+2\lambda^3)x^3\,\big)
$$}\fi
$$
   (x,y,z)\mapsto\big(\,x,\, (x^2-y)/\alpha,\, (2x^3 - 3xy + z)/(2\alpha^{3/2})\,\big)
$$
transforms (i${}_1$) and (i${}_6$) into, respectively,
\roster
\item"(i${}_1^*$)" $\Omega=\{y^3>z^2\}$, $g=g_0$ (see below), and
    $\rho=(y^3-z^2)^{p-1}e^{-\lambda\alpha(x^2+y)}$;
\item"(i${}_6^*$)" $\Omega$ is the only bounded component of
$\bR^3\setminus\{(y^3-z^2)(1-x^2-y)=0\}$, $g=g_0-g_1$,
$\rho=(y^3-z^2)^{p-1}(1-x^2-y)^{q-1}$, where
$$
    g_0 = \left(\matrix
      1 &  0 &  0  \\
      0 & 4y & 6z  \\
      0 & 6z & 9y^2
    \endmatrix\right), \qquad
    g_1 = \left(\matrix
      x^2 & 2xy  & 3xz  \\
      2xy & 4y^2 & 6yz  \\
      3xz & 6yz  & 9z^2
    \endmatrix\right) \qquad
$$
\endroster
($g_1$ is the coefficient of $e$ in the matrix in Lemma~\lemCubicI).
Thus we have one-parameter families of pairwise non-equivalent solutions
of the DOP problem such that the members of each family become equivalent to each other
when they are considered as solutions of the weighted DOP problem with suitable weights. 
The same phenomenon was observed in dimension 2 in [\refBOZ, \S4.5],
[\refOr, Remark~6.5].
\medskip

\subhead\sectFig. Comments on the figures
\endsubhead
In Figures \figAB--\figABD\ we show the bounded domains appearing in Theorem~\thB\
and the domain discussed in \S\sectD. All the domains are curvilinear polyhedra
(all but one being tetrahedra), so we present them by the planar
projections of their edges.
When the axes are not shown, the projection is assumed to be $(x,y,z)\mapsto(x,y)$
(or $(x,y,z)\mapsto(x,z)$ in Figure~\figABD\ on the right). The number $n$
near an edge means that the surface is given by the
equation $v^2=u^n$ in some local curvilinear coordinates $(u,v,w)$ in a neighbourhood
of this edge. 
In all the cases the metric $(g_{ij})=g^{-1}$ is of constant non-negative curvature
and the boundary of $\Omega$ is totally geodesic
(which well agrees with Soukhanov's results [\refSouArx], [\refSouAFST]).
Thus $\Omega$ can be identified with
the quotient of $\bR^3$ or $\bS^3$ by a group generated by reflections
(a Coxeter group) and $\LL$ is the image of the Laplace operator.
The types of the groups are indicated in the figure captions. The explicit formulas
for these identifications are given in \S\sectQuo. Note that if an edge of
$\Omega$ is marked by a number $n$, then the angle at the corresponding edge of
the fundamental polyhedron in $\bR^3$ or $\bS^3$ is $\pi/n$.
For affine Coxeter groups (Figures~\figAaff\ and \figCaff) we also present
the fundamental tetrahedra and the corresponding Coxeter graphs.
Notice also that the curves on the right hand sides of Figures~\figAaff\ and \figABD\
are $(3,1)$- and $(4,1)$-hypocycloids.


\head\sectQuo. Higher dimensional solutions of the DOP problem on the 
      quotients of $\bS^n$ or $\bR^n$ by Coxeter groups
\endhead

Using the approach from [\refBOZ,~\S4], in this section we realize
each solution from Theorem~\thB\ as an image of the Laplace operator on
$\bS^3$ or $\bR^3$ through the quotient by a discrete group generated by
reflections (a Coxeter group). Moreover, we include each of these solutions
into an infinite series of solutions in all dimensions.

\subhead\sectGa.  Generalities 
\endsubhead
With a second order differential operator $\LL$ with no 0-order term
on a manifold $M$, 
is associated the operator ``carr\'e du champ''
$$
   \Ga_{\LL}(f_1,f_2)=\frac12\Big(\LL(f_1f_2)-f_1\LL(f_2)-f_2\LL(f_1)\Big).
$$
(see [\refBGL]). Notice that the operator $\Ga_{\De}$ (for the Laplace
operator on $\bR^n$) plays a key r\^ole in [\refSYS] where it is denoted
by $\langle df_1,df_2\rangle$.
If $\LL$ is given by (\eqLL) in some coordinates $(x_1,\dots,x_n)$, then
$g^{ij}=\Gamma_{\LL}(x_i,x_j)$ and $b^i=\LL(x_i)$.
Let $\ff:M\to\bR^n$, $p\mapsto(f_1(p),\dots,f_n(p))$
be a mapping such that $\Gamma_{\LL}(f_i,f_j)=G^{ij}\circ\ff$ and
$\LL(f_i)=B^i\circ\ff$ for some functions $G^{ij}$ and $B^i$ defined on $\ff(M)$.
Then a direct computation shows that the operator
$$
   \ff_*(\LL)=\sum_{i,j}G^{ij}\partial_{ij} + \sum_i B^i\partial_i
$$
is such that $\ff_*(\LL)(\varphi)=\LL(\varphi\circ\ff)$ for any smooth
$\varphi:\ff(M)\to\bR$. We say that $\LL_*(\ff)$ is the
{\it image of $\LL$ through} $\ff$.

Let $G$ be a discrete group generated by orthogonal reflections acting on $\bR^n$
(see [\refBou] for a general introduction to the subject).
We discuss here only bounded solutions of the DOP problem. Therefore,
when $G$ is finite (a {\it spherical Coxeter group} or just {\it Coxeter group}),
we assume that the origin is a fixed point and we
restrict the action form $\bR^n$ to the unit sphere $\bS^{n-1}$.
If $G$ is infinite (an {\it affine Coxeter group}),
we assume that it contains a full rank subgroup of translations.
So, in both cases the orbit space $M/G$ is compact ($M$ is $\bR^n$ or $\bS^{n-1}$).

If $G$ is finite, it is known (see [\refCox], [\refBou, Ch.~V, \S\S5--6]) that
the ring of invariant polynomials is freely generated by some invariant homogeneous
forms $I_1,\dots,I_n$. The choice of the invariants $I_j$'s
is not unique (see, e.g., [\refMe], [\refSYS] for different concrete choices)
but their degrees $d_1,\dots,d_n$ are uniquely determined.
These numbers (called {\it exponents} in [\refBou])
for each Coxeter group can be found in Tables (Planches) I--X in [\refBou]. 
One of the basic invariants is (if the action is irreducible)
or can be chosen to be $x_1^2+\dots+x_n^2$.
Let it be $I_1$. Then [\refBOZ,~Eq.~(4.5)] implies that the image of the Laplace
operator $\De_{\bS^{n-1}}$ for $\ff:\bS^{n-1}\to\bR^{n-1}$,
$p\mapsto(I_2(p),\dots,I_n(p))$, is a solution
of the weighted DOP problem (see [\refBB], [\refOr], [\refSouAFST] for the definition)
with weights $(d_2,\dots,d_n)$ on $\ff(\bS^{n-1})$.
However, in \S\S\sectA--\sectAB\ we show that for the Coxeter groups of types
$A_n$, $B_n$, and their direct products, as well as for $D_4$, one can choose
the basic invariants so that the image of the Laplace operator is a solution
of the DOP problem (with weights $(1,\dots,1)$).

Consider now the case when $G$ is an affine Coxeter group acting on $E=\bR^n$.
It is shown in [\refBou, Ch.~VI, \S3.4] that the ring of invariant Fourier
polynomials is freely generated by certain elements $f_1,\dots,f_n$ which are
explicitly described via the fundamental weights $\omega_1,\dots,\omega_n$
corresponding to some Weyl chamber $C$. One can check that the image of $\De_E$ through
$p\mapsto(f_1(p),\dots,f_n(p))$ is a solution of the weighted DOP problem
with the weights $\alpha(\omega_1),\dots,\alpha(\omega_n)$ where $\alpha$
is any linear function positive on $C$. In \S\S\sectCt--\sectAt\ we show that
these are also solutions of the DOP problem (with weights $(1,\dots,1)$) for
the affine Coxeter groups of types $\tilde A_n$ and $\tilde C_n$.
It seems plausible that the quotients by other spherical or affine
Coxeter groups never give a solution of the DOP problem.
In dimension $2$ this fact follows from the classification in [\refBOZ].

For each solution $(\Omega,g,\rho)$ obtained as the image of a Laplace operator through
the quotient by a Coxeter group,
$\LL$ is the Laplace-Beltrami operator for the metric $g^{-1}$,
hence $\rho=(\det g)^{-1/2}$.


\subhead\sectA. Quotient of $\bS^{n-2}$ by the Coxeter group $A_{n-1}$
\endsubhead

Let $E=\bR^n$ with coordinates $x_1,\dots,x_n$, let $H\subset E$ be the
hyperplane $x_1+\dots+x_n=0$, and $\bS^{n-2}$ be the unit sphere in $H$.
The Coxeter group $A_{n-1}$ acting on $\bS^{n-2}$ is generated by the
orthogonal reflections in the hyperplanes $x_i=x_j$. The ring of
invariants is freely generated by the
elementary symmetric polynomials $s_2,\dots,s_n$.
So, we consider the mapping $\Phi:\bS^{n-2}\to\bR^{n-2}$,
$(x_1,\dots,x_n)\mapsto(s_3,\dots,s_n)$ where
$P(u)=(u+x_1)\dots(u+x_n)=\sum_{k=0}^n s_k u^{n-k}$.
Note that $(s_0,s_1,s_2)|_{\bS^{n-2}}=(1,0,-\frac12)$ and we set $s_k=0$ for
$k\not\in[0,n]$. Then $\Omega=\Phi(\bS^{n-2})$ is bounded by
the hypersurface
$$
   \text{discr}_u(u^n-\tfrac12u^{n-2}+X_{3}u^{n-3}+\dots+X_{n-1}u+X_n)=0.
$$
(cf.~Thm.~\thB(v)). Here $(X_3,\dots,X_n)$ are the coordinates
in the target space $\bR^{n-2}$.
Let $\De=\De_{\bS^{n-2}}$ and let $\Ga$ be the corresponding
carr\'e du champ. We are going to check that
$\Phi_*(\De)$ is a Laplace-Beltrami solution of the DOP problem on $\Omega$.
We have $\Ga(s_k,s_m) = \Ga_H(s_k,s_m) - kms_ks_m$ (see [\refBOZ, Eq.~(4.5)])
and $\Ga_H(s_k,s_m)$ is the coefficient of $u^{n-k}v^{n-m}$ in $\Ga_H(P(u),P(v))$.
We have
$$
     \De_E=\De_H+\tfrac{1}{n}\partial_0^2
     \quad\text{ where }\quad \partial_0=\sum\partial_i.       \eqno(\eqModelAnI)
$$
Hence $\Ga_H(f_1,f_2) = \Ga_E(f_1,f_2)-\frac{1}{n}(\partial_0f_1)(\partial_0f_2)$.
It is clear that
$$
   \partial_0 P(u) = \sum_i\frac{P(u)}{u+x_i} = P'(u),          \eqno(\eqModelAnII)
$$
thus $\Ga_H(P(u),P(v))=\Ga_E(P(u),P(v))-\frac{1}{n}P'(u)P'(v)$.
Finally, by [\refBOZ, p.~1033],
$$
\split
  &\Ga_E(P(u),P(v))
 =\sum_{i,j}\big(\partial_i P(u)\big)\big(\partial_j P(v)\big)\Ga_E(x_i,x_j)
   = \sum_i\big(\partial_i P(u)\big)\big(\partial_i P(v)\big)
\\&   =\sum_i\frac{P(u)P(v)}{(u+x_i)(v+x_i)}
     =\frac{P(u)P(v)}{v-u}\sum_i\Big(\frac{1}{u+x_i}-\frac{1}{v+x_i}\Big)
\\&
   \overset{\text{(\eqModelAnII)}}\to=\frac{P'(u)P(v)-P'(v)P(u)}{v-u}
  =\sum_{k,m}(n-k)s_ks_m\frac{u^{n-k-1}v^{n-m}-v^{n-k-1}u^{n-m}}{v-u}
\if01{
\\& =
  \sum_{k=0}^{n-1}\,\sum_{m=0}^k\!(n-k)s_ks_m\!\sum_{l=0}^{k-m}\!u^{n-k+l-1}v^{n-m-l-1}
\\&\qquad
  -\!\sum_{k=0}^{n-1}\,\sum_{m=k+2}^n\!(n-k)s_ks_m\!
                       \sum_{l=1}^{m-k-1}\!u^{n-k-l-1}v^{n-m+l-1}
}\fi%
\\& =
  \sum_{k,m}\,\!(n-k)s_ks_m\left(
   \sum_{l=0}^{k-m}\!u^{n-k+l-1}v^{n-m-l-1}
                   -\!\!\sum_{l=1}^{m-k-1}\!\!u^{n-k-l-1}v^{n-m+l-1}\right),
\endsplit
$$
hence for $a\le b$ we have
$$
  \Ga(s_a,s_b) = 
    (a-1)(1-\tfrac{b-1}{n})s_{a-1}s_{b-1} - abs_as_b
      + \sum_{l\ge 1}(a-b-2l)s_{a-l-1}s_{b+l-1}.
$$
Thus the coefficients $g^{ab}$, $a\le b$, of $\Phi_*(\De)$ are given by
the same expression where $s_0,s_1,\dots,s_n$ are replaced by
$1,0,-\tfrac12,X_3,\dots,X_n$ and $s_k$ is set to zero when $k\not\in[0,n]$.
Here $X_3,\dots,X_n$ are coordinates in the target space $\bR^{n-2}$.

By [\refBOZ, Eq.~(4.5)] we have $\De(s_a)=\De_H(s_a)-a(n+a-3)s_a$.
Counting the number of monomials, one obtains
$$
    \partial_0(s_a) = (n-a+1)s_{a-1}.              \eqno(\eqModelAnIII)
$$
These formulae combined with (\eqModelAnI) and with $\De_E(s_a)=0$ yield
$$
  \De(s_a)=-\tfrac1n(n-a+1)(n-a+2)s_{a-2} - a(n+a-3)s_a.  \eqno(\eqModelAnIV)
$$
\if01{

(***   A_{n-1}   ***)
n=5;
PX=Sum[s[k]X^(n-k),{k,0,n}]; dPX=D[PX,X];
PY=Sum[s[k]Y^(n-k),{k,0,n}]; dPY=D[PY,Y];
EPP=Factor[(dPX*PY-dPY*PX)/(Y-X)]
HPP=EPP-dPX*dPY/n
G=Table[Coefficient[Coefficient[HPP,X,n-k],Y,n-m]-k*m*s[k]*s[m],{k,3,n},{m,3,n}]/.
{s[0]->1,s[1]->0,s[2]->-1/2}
F=Factor[Det[G]]
Factor[(G.Table[D[F,s[i]],{i,3,n}])/F]
Factor[F/Discriminant[T^n - 1/2*T^(n-2) + Sum[s[k]T^(n-k),{k,3,n}],T]]

GG=Table[Coefficient[Coefficient[HPP,X,n-k],Y,n-m]-k*m*s[k]*s[m],{k,3,n},{m,3,n}]
GE=Table[Coefficient[Coefficient[EPP,X,n-k],Y,n-m],{k,1,n},{m,1,n}]
GS=Table[Coefficient[Coefficient[EPP,X,n-k],Y,n-m]-k*m*s[k]*s[m],{k,1,n},{m,1,n}]/.s[0]->1
Factor[Det[GS]/Discriminant[T^n + Sum[s[k]T^(n-k),{k,1,n}],T]]
(**  GS is the quotient by  A0+A1+A_n  **)

EPP1 = Sum[(n-k)s[k]s[m]X^(n-k+l-1)Y^(n-m-l-1),{k,0,n-1},{m,0,k},  {l,0,k-m}]
EPP2 = Sum[(n-k)s[k]s[m]X^(n-k-l-1)Y^(n-m+l-1),{k,0,n-1},{m,k+2,n},{l,1,m-k-1}]

GE0 =
Table[(n-b+1)s[b-1]s[a-1] + Sum[(a-b-2l)s[a-l-1]s[b+l-1],{l,n}],
  {a,n},{b,n}]/.
  Join@@Table[{s[-i]->0,s[n+i]->0},{i,n}]
Do[GE0[[b,a]]=GE0[[a,b]],{a,n},{b,a,n}]

G0 =
Table[(n-b+1)s[b-1]s[a-1] - a*b*s[a]*s[b]
       - (n-a+1)(n-b+1)s[a-1]s[b-1]/n + Sum[(a-b-2l)s[a-l-1]s[b+l-1],{l,n}],
  {a,3,n},{b,3,n}]/.
  Join@@Table[{s[-i]->0,s[n+i]->0},{i,n}]
Do[G0[[b,a]]=G0[[a,b]],{a,n-2},{b,a,n-2}];
G0=G0/.{s[0]->1,s[1]->0,s[2]->-1/2}

G00 =
Table[(a-1)(n-b+1)s[b-1]s[a-1]/n - a*b*s[a]*s[b] + Sum[(a-b-2l)s[a-l-1]s[b+l-1],{l,n}],
  {a,3,n},{b,3,n}]/.
  Join@@Table[{s[-i]->0,s[n+i]->0},{i,n}]
Do[G00[[b,a]]=G00[[a,b]],{a,n-2},{b,a,n-2}];
G00=G00/.{s[0]->1,s[1]->0,s[2]->-1/2}
F=Factor[Det[G00]];
Factor[F/Discriminant[T^n - 1/2*T^(n-2) + Sum[s[k]T^(n-k),{k,3,n}],T]]
(* De *)
rho=F^(-1/2)
De=Factor[Sum[D[G00[[i-2]],s[i]],{i,3,n}]+G00.Table[D[Log[rho],s[i]],{i,3,n}]]
Table[-a(n+a-3)s[a] - (n-a+1)(n-a+2)/n*s[a-2], {a,3,n}]/.
{s[0]->1,s[1]->0,s[2]->-1/2}
Factor[

(*   A_{n-1} + A_1    *)

n=6
G0 =
Table[(a-1)(n-b+1)s[b-1]s[a-1]/n - a*b*s[a]*s[b] + Sum[(a-b-2l)s[a-l-1]s[b+l-1],{l,n}],
  {a,2,n},{b,2,n}]/.
  Join@@Table[{s[-i]->0,s[n+i]->0},{i,n}];
Do[G0[[b,a]]=G0[[a,b]],{a,n-1},{b,a,n-1}];
G0=G0/.{s[0]->1,s[1]->0};
F0=Factor[Det[G0]];
Factor[F0/Discriminant[T^n + Sum[s[k]T^(n-k),{k,2,n}],T]]
rho=F0^(-1/2);
De=Factor[Sum[D[G0[[i-1]],s[i]],{i,2,n}]+G0.Table[D[Log[rho],s[i]],{i,2,n}]]
De1=Table[-(n-a+1)(n-a+2)/n*s[a-2]-a(n+a-3)s[a]-a*s[a],{a,2,n}]/.{s[0]->1,s[1]->0}
Factor[De/De1]

(*   A_{n-1} + A_1 + A_0   *)

n=5
G0 =
Table[(a-1)(n-b+1)s[b-1]s[a-1]/n - a*b*s[a]*s[b] + Sum[(a-b-2l)s[a-l-1]s[b+l-1],{l,n}],
  {a,n},{b,n}]/.
  Join@@Table[{s[-i]->0,s[n+i]->0},{i,n}];
G0=G0/.{s[0]->1,s[1]->0};
G0[[1,1]]=1-s[1]^2;
Do[G0[[1,k]]=-k*s[1]s[k],{k,2,n}];
Do[G0[[b,a]]=G0[[a,b]],{a,n},{b,a,n}];
F0=Factor[Det[G0]];
Factor[(G0.Table[D[F0,s[i]],{i,n}])/F0]
Factor[F0/Discriminant[T^n + Sum[s[k]T^(n-k),{k,2,n}],T]]
rho=F0^(-1/2);
De=Factor[Sum[D[G0[[i]],s[i]],{i,n}]+G0.Table[D[Log[rho],s[i]],{i,n}]]
De1=Table[-(n-a+1)(n-a+2)/n*s[a-2]-a(n+a-3)s[a]-2a*s[a],{a,n}]/.{s[0]->1,s[1]->0}
De1[[1]]=-n*s[1];
Factor[De/De1]

(*   A_{n-1} + A_1 + A_0,  another projection (the most natural)  *)
n=5
G0 =
Table[(n-b+1)s[b-1]s[a-1] - a*b*s[a]*s[b] + Sum[(a-b-2l)s[a-l-1]s[b+l-1],{l,n}],
  {a,n},{b,n}]/.
  Join@@Table[{s[-i]->0,s[n+i]->0},{i,n}];
G0=G0/.{s[0]->1};
Do[G0[[b,a]]=G0[[a,b]],{a,n},{b,a,n}];
F0=Factor[Det[G0]];
Factor[(G0.Table[D[F0,s[i]],{i,n}])/F0]
Factor[F0/Discriminant[T^n + Sum[s[k]T^(n-k),{k,1,n}],T]]
rho=F0^(-1/2);
De=Factor[Sum[D[G0[[i]],s[i]],{i,n}]+G0.Table[D[Log[rho],s[i]],{i,n}]]
De1=Table[-a(n+a-1)s[a],{a,n}]
Factor[De/De1]

(**  B_n  **)

(**  first projection:  t_i = x_i^2 **)

n=5;
G0 = 4Table[-a*b*s[a]s[b]+Sum[(b-a+2l-1)s[a-l]s[b+l-1],{l,n}], {a,2,n},{b,2,n}]/.
   Join@@Table[{s[-i]->0,s[n+i]->0},{i,n}];
Do[G0[[b,a]]=G0[[a,b]],{a,n-1},{b,a,n-1}];
G0=G0/.{s[0]->1,s[1]->1}
F=Factor[Det[G0]]
Factor[(G0.Table[D[F,s[i]],{i,n}])/F]
Factor[ F/Discriminant[T^n + T^(n-1) + Sum[s[k]T^(n-k),{k,2,n}],T] ]  (* s[n] *)
(* De *)
rho=F^(-1/2)
De=Factor[Sum[D[G0[[i-1]],s[i]],{i,2,n}]+G0.Table[D[Log[rho],s[i]],{i,2,n}]]
Table[-2a(n+2a-2)s[a] + 2(n-a+1)s[a-1], {a,2,n}]/.
{s[0]->1,s[1]->1}
Factor[

n=5;         (***  2nd proj, first version:  t_i = x_i^2 - 1/n  ***)
G0 = Table[
     ((n-b+1)s[a-1]s[b-1] + Sum[(a-b-2l)s[a-l-1]s[b+l-1],{l,n}])/n
     + Sum[(b-a+2l-1)s[a-l]s[b+l-1],{l,n}]
     - (a*s[a]+(n-a+1)/n*s[a-1])(b*s[b]+(n-b+1)/n*s[b-1]), {a,2,n},{b,2,n}]/.
   Join@@Table[{s[-i]->0,s[n+i]->0},{i,n}];
Do[G0[[b,a]]=G0[[a,b]],{a,n-1},{b,a,n-1}];
G0=G0/.{s[0]->1,s[1]->0}
F=Factor[Det[G0]]
Factor[(G0.Table[D[F,s[i]],{i,2,n}])/F]
Factor[ F/Discriminant[T^n + Sum[s[k]T^(n-k),{k,2,n}],T] ]

n=5;         (***  2nd proj, second version:  t_i = n x_i^2 - 1  ***)
G0 = Table[
     (n-b+1)s[a-1]s[b-1] + Sum[(a-b-2l)s[a-l-1]s[b+l-1],{l,n}]
     + Sum[(b-a+2l-1)s[a-l]s[b+l-1],{l,n}]
     - (a*s[a]+(n-a+1)s[a-1])(b*s[b]+(n-b+1)s[b-1])/n, {a,2,n},{b,2,n}]/.
   Join@@Table[{s[-i]->0,s[n+i]->0},{i,n}];
Do[G0[[b,a]]=G0[[a,b]],{a,n-1},{b,a,n-1}];
G0=G0/.{s[0]->1,s[1]->0}
F=Factor[Det[G0]]
Factor[(G0.Table[D[F,s[i]],{i,2,n}])/F]
Factor[ F/Discriminant[T^n + Sum[s[k]T^(n-k),{k,2,n}],T] ]
(* De *)
rho=F^(-1/2)
De=Factor[4n*(Sum[D[G0[[i-1]],s[i]],{i,2,n}]+G0.Table[D[Log[rho],s[i]],{i,2,n}])]
Table[2(2a-2a^2-a*n)s[a] + (n-a+1)(8-8a)s[a-1] - 4(n-a+1)(n-a+2)s[a-2], {a,2,n}]/.
{s[0]->1,s[1]->0}
Factor[

(**  Bn + A1   first projection:  t_i = x_i^2 **)

n=5;
G0 = 4Table[-a*b*s[a]s[b]+Sum[(b-a+2l-1)s[a-l]s[b+l-1],{l,n}], {a,n},{b,n}]/.
   Join@@Table[{s[-i]->0,s[n+i]->0},{i,n}];
Do[G0[[b,a]]=G0[[a,b]],{a,n},{b,a,n}];
G0=G0/.{s[0]->1};
F=Factor[Det[G0]];
Factor[(G0.Table[D[F,s[i]],{i,n}])/F]
Factor[ F/Discriminant[T^n + Sum[s[k]T^(n-k),{k,n}],T] ]  (***  (s[1]-1)s[n]  ***)
rho=F^(-1/2);
De=Factor[Sum[D[G0[[i]],s[i]],{i,n}]+G0.Table[D[Log[rho],s[i]],{i,n}]]
De1=Table[-2a(n+2a-2)s[a] + 2(n-a+1)s[a-1] -2a*s[a], {a,n}]/.{s[0]->1}
Factor[De1/De]

(**  Bn + A1   second projection:  t_i = n x_i^2 - 1 **)

n=3;
G0 = Table[
     (n-b+1)s[a-1]s[b-1] + Sum[(a-b-2l)s[a-l-1]s[b+l-1],{l,n}]
     + Sum[(b-a+2l-1)s[a-l]s[b+l-1],{l,n}]
     - (a*s[a]+(n-a+1)s[a-1])(b*s[b]+(n-b+1)s[b-1])/n, {a,n},{b,n}]/.
   Join@@Table[{s[-i]->0,s[n+i]->0},{i,n}];
Do[G0[[b,a]]=G0[[a,b]],{a,n},{b,a,n}];
G0=G0/.{s[0]->1}
F=Factor[Det[G0]]
Factor[(G0.Table[D[F,s[i]],{i,n}])/F]
Factor[ F/Discriminant[T^n + Sum[s[k]T^(n-k),{k,n}],T] ]
}\fi

\subhead\sectB. Quotient of $\bS^{n-1}$ by the Coxeter group $B_n$
\endsubhead

Let $E$ be $\bR^n$ with coordinates $x_1,\dots,x_n$ and $\bS^{n-1}$ be the unit sphere
in $E$. The Coxeter group $B_n$ acting on $E$ is generated by the reflections in
the hyperplanes $x_i=x_j$ and $x_i=0$. The ring of polynomial invariants is
generated by the elementary symmetric polynomials in $x_i^2$.
We consider the mapping $\Phi:\bS^{n-1}\to\bR^{n-1}$,
$(x_1,\dots,x_n)\mapsto(s_2,\dots,s_n)$ where
$P(u)=(u+t_1)\dots(u+t_n)=\sum_{k=0}^ns_ku^{n-k}$, $t_i=x_i^2$. We have
$(s_0,s_1)|_{\bS^{n-1}}=(1,1)$ and we set $s_k=0$ for $k\not\in[0,n]$.
Then $\Phi(\bS^{n-1})$ is bounded by the hypersurface
$$
   X_n\, \text{discr}_u(u^n + u^{n-1} + X_2 u^{n-2} + \dots + X_{n-1} u + X_n)=0
$$
(cf.~Thm.~\thB(iii${}_2$) and Figure~\figABD).
Its component $X_n=0$ is the image of the hyperplanes $x_i=0$ and the other
component is the image of the planes $x_i=x_j$.

Let $\De=\De_{\bS^{n-1}}$ and let $\Ga$ be the corresponding carr\'e du champ.
$\Ga(s_k,s_m)$ is the coefficient of $u^{n-k}v^{n-m}$ in $\Ga(P(u),P(v))$.
The $s_k$ are homogeneous of degree $2k$, hence (see [\refBOZ, Eq.~(4.5)])
$\Ga(s_k,s_m)=\Ga_E(s_k,s_m)-4kms_ks_m$ and
$$
 \Ga_E(t_i,t_j)=\Ga_E(x_i^2,x_j^2)=4x_ix_j\Ga_E(x_i,x_j)=4t_i\delta_{ij}.
$$
Then (cf.~[\refSYS,~Prop.~2.2.2])
$$
  \frac14\Ga_E\big(P(u),P(v)\big)=\frac14\sum_{i,j}
   \big(\partial_{t_i}P(u)\big)
   \big(\partial_{t_j}P(v)\big)\Ga_E(t_i,t_j)
   =\sum_i\frac{t_iP(u)P(v)}{(u+t_i)(v+t_i)}
$$
$$
   =\frac{P(u)P(v)}{u-v}\sum_i\Big(\frac{u}{u+t_i}-\frac{v}{v+t_i}\Big)
   =\frac{uP'(u)P(v) - vP'(v)P(u)}{u-v}
$$
$$
  =\sum_{k,m}(n-k)s_ks_m\frac{u^{n-k}v^{n-m}-v^{n-k}u^{n-m}}{u-v}
$$
$$
  =\sum_{k,m}\,\!(n-k)s_k s_m \left(
                        \sum_{l=1}^{m-k}\!\!u^{n-k-l}v^{n-m+l-1}
                   -\!\!\sum_{l=1}^{k-m}\!u^{n-k+l-1}v^{n-m-l}    \right)
$$
Hence for $a\le b$ we have
$$
\Ga(s_a,s_b)=-4ab\, s_a s_b +
  \sum_{l\ge 1}4(b-a+2l-1)s_{a-l}s_{b+l-1}.           \eqno(\eqModelBn)
$$
The coefficients $g^{ab}$, $a\le b$, of $\Phi_*(\De)$ are given by
the same expression where $s_0,s_1,\dots,s_n$ are replaced by
$1,1,X_2,\dots,X_n$ and $s_k$ is set to zero when $k\not\in[0,n]$.
\if01{
n=7;
A=Sum[(n-k)s[k]s[m]Factor[(u^(n-k)v^(n-m)-v^(n-k)u^(n-m))/(u-v)],{k,0,n},{m,0,n}]
B=Sum[(n-k)s[k]s[m]
    Factor[Sum[u^(n-k-l)v^(n-m+l-1),{l,m-k}]-Sum[u^(n-k+l-1)v^(n-m-l),{l,k-m}]],
{k,0,n},{m,0,n}];
Factor[A-B]
CC=Sum[u^(n-a)v^(n-b)Sum[(b-a+2l-1)s[a-l]s[b+l-1],{l,n}],{a,n},{b,n}]/.
   Flatten[{Table[{s[n+k]->0,s[-k]->0},{k,n}]}]
}\fi

We have $\De(s_a)=\De_E(s_a) - 2a(n+2a-2)s_a$ (see [\refBOZ, Eq.~(4.5)]).
Similarly to (\eqModelAnIII) one obtains $\De_E(s_a)=2(n-a+1)s_{a-1}$,
hence
$$
     \De(s_a)=2(n-a+1)s_{a-1} - 2a(n+2a-2)s_a.           \eqno(\eqModelBnDe)
$$
\subhead\sectBb. Quotient of $\bS^{n-1}$ by the Coxeter group $B_n$
                 (another mapping)
\endsubhead

In this subsection we compute $\Phi_*(\De_{\bS^{n-1}})$ for another polynomial
mapping $\Phi$ invariant under the action of $B_n$.
This time $\Omega=\Phi(\bS^{n-1})$ is bounded by
$$
  \{X=(X_2,\dots,X_n)\mid P_X(1)\,\text{discr}_u P_X(u)=0\},\quad
  P_X(u)=u^n + \sum_{k=2}^n X_k u^{n-k}.
$$
Up to rescaling of the coordinates,
we obtain the solution in Thm.~\thB(iii${}_2$) (see Figure~\figABD) when $n=4$,
and the solution in [\refBOZ, \S4.10] when $n=3$.

The mapping $\Phi:\bS^{n-1}\to\bR^{n-1}$ is given by
$(x_1,\dots,x_n)\mapsto(s_2,\dots,s_n)$ where
$P(u)=(u+t_1)\dots(u+t_n)=\sum_{k=0}^ns_ku^{n-k}$, $t_i=nx_i^2-1$.
It factors through $\Phi_1:(x_2,\dots,x_n)\mapsto(t_2,\dots,t_n)$
and $\Phi_1(\bS^{n-1})$ is the $(n-1)$-simplex $\sigma$
given by $\sum t_i=0$, $t_i\ge-1$.
The image of $\partial\sigma$ is the hyperplane $P_X(1)=0$, and the image
of the $(n-2)$-planes $\sigma\cap\{t_i=t_j\}$ is the discriminantal hypersurface.
This solution is obtained from the one in \S\sectB\ by the change of variable
$u\mapsto u-\frac1n$ which corresponds to an evident affine linear transformation
in the coefficient space. It seems, however, that it is easier to recompute
$\Phi_*(\De)$ rather than to perform this change of variables. Let us do it.
By linearity, 
$\Ga(s_k,s_m)$ is the coefficient of $u^{n-k}v^{n-m}$ in $\Ga(P(u),P(v))$.
We have (see [\refBOZ, Eq.~(4.2)])
$$
 \Ga(t_i,t_j)=\Ga(nx_i^2,nx_j^2)=4n^2x_ix_j(\delta_{ij}-x_ix_j)
     =4n^2\delta_{ij}x_i^2-4n^2x_i^2x_j^2.
$$
Hence
$$
  \frac{\Ga(P(u),P(v))}{4n^2}=\sum_{i,j}
   \big(\partial_{t_i}P(u)\big)
   \big(\partial_{t_j}P(v)\big)\frac{\Ga(t_i,t_j)}{4n^2}
$$
$$
  =\sum_{i,j}\frac{P(u)P(v)\Ga(t_i,t_j)}{4n^2(u+t_i)(v+t_j)}
   =\sum_i\frac{P(u)P(v)x_i^2}{(u+t_i)(v+t_i)}
   -\sum_{i,j}\frac{P(u)P(v)x_i^2x_j^2}{(u+t_i)(v+t_j)}
$$
$$
   =\frac{P(u)P(v)x_i^2}{v-u}\sum_i\Big(\frac{1}{u+t_i}-\frac{1}{v+t_i}\Big)
   -\sum_{i,j}\frac{P(u)x_i^2P(v)x_j^2}{(u+t_i)(v+t_j)}
$$
$$
   =\frac{Q(u)P(v)-Q(v)P(u)}{v-u}-Q(u)Q(v)
$$
where
$$
\split
   Q(u)&=\sum_i\frac{P(u)x_i^2}{u+t_i}
       =\frac1{n}\sum_i\frac{P(u)(t_i+1)}{u+t_i}
       =\frac1n\sum_i\Big(P(u)-\frac{(u-1)P(u)}{u+t_i}\Big)
\\&
       =P(u)-\frac1n(u-1)P'(u)
       =\frac1n\sum_k s_k\Big(ku^{n-k}+(n-k)u^{n-k-1}\Big)
\endsplit
$$
Thus $\big(Q(u)P(v)-Q(v)P(u)\big)/(v-u)$ is equal to
$$
\split
   & \sum_{k,m}\!\frac{s_ks_m}n\!\left(
    k\,\frac{u^{n-k}v^{n-m}-v^{n-k}u^{n-m}}{v-u}
    +(n-k)\frac{u^{n-k-1}v^{n-m}-v^{n-k-1}u^{n-m}}{v-u}\right)
\\&\qquad =
  \sum_{k,m}\,\!\frac{s_k s_m}n\Bigg\{k\left(
   \sum_{l=1}^{k-m}\!u^{n-k+l-1}v^{n-m-l}
                   -\!\!\sum_{l=1}^{m-k}\!\!u^{n-k-l}v^{n-m+l-1}\right)
\\& \qquad\qquad\qquad+(n-k)\left(
   \sum_{l=0}^{k-m}\!u^{n-k+l-1}v^{n-m-l-1}
                   -\!\!\sum_{l=1}^{m-k-1}\!\!u^{n-k-l-1}v^{n-m+l-1}\right)\Bigg\}
\endsplit
$$
If $a\le b$, then
$$
\split
   \tfrac14\Gamma&(s_a,s_b)=
 \sum_{l\ge1}n(a-b-2l)s_{a-l-1}s_{b+l-1}
     + \sum_{\l\ge 1}n(b-a+2l-1)s_{a-l}s_{b+l-1}
\\&+n(n-b+1)s_{a-1}s_{b-1}
 - \Big(as_a+(n-a+1)s_{a-1}\Big)
      \Big(bs_b+(n-b+1)s_{b-1}\Big).
\endsplit
$$
The coefficients $g^{ab}$, $a\le b$, of $\Phi_*(\De)$ are given by
the same expression where $s_0,s_1,\dots,s_n$ are replaced by
$1,0,X_2,\dots,X_n$ and $s_k$ is set to zero when $k\not\in[0,n]$.

We have $\De=\De_E-(r\partial_r)^2-(n-2)r\partial_r$ where
$r\partial_r=\sum x_i\partial_{x_i}$
(see [\refBOZ, Eq.~(4.4)])
and $\De_E P = 2nP'$, $r\partial_r P=2nP+2(1-u)P'$, hence
$$
    (r\partial_r)^2 P=4n^2P + 4(2n-1)(1-u)P' + 4(1-u)^2 P'',
$$
$$
   \De P = 2nP' - 2(3n^2-n)P - 2(5n-3)(1-u)P' - 4(1-u)^2 P'',
$$
and we obtain
$$
  \De s_a = 2a(2-2a-n)s_a + 8(n-a+1)(1-a)s_{a-1} - 4(n-a+1)(n-a+2)s_{a-2}.
$$


\subhead\sectD.
           Quotient of $\bS^{n-1}$ by the Coxeter group $D_n$
\endsubhead
Let the notation be as in \S\sectB.
The Coxeter group $D_n$ acting on $E$ is generated by the reflections in
the hyperplanes $x_i\pm x_j=0$. The ring of polynomial invariants is generated
by $s_1,\dots,s_{n-1}$ and $\hat s_n=\sqrt{s_n}=x_1\dots x_n$.
The values of $\Ga(s_a,s_b)$ and $\De(s_a)$
are already computed in \S\sectB, and we have (recall that
$(s_0,s_1)|_{\bS^{n-1}}=(1,1)$)
$$
\split
    \Ga(s_a,\hat s_n) &= \Ga(s_a,s_n^{1/2}) = \tfrac12 s_n^{-1/2}\Ga(s_a,s_n)
    \overset{\text{(\eqModelBn)}}\to=
     -2ans_a\hat s_n + 2(n-a+1)s_{a-1}\hat s_n,
\\
    \Ga(\hat s_n,\hat s_n) &= \Ga(s_n^{1/2},s_n^{1/2}) = \tfrac14 s_n^{-1}\Ga(s_n,s_n)
    \overset{\text{(\eqModelBn)}}\to=
     -n^2 \hat s_n^2 + s_{n-1},
\\
    \De(\hat s_n) &= \De_E(\hat s_n) - 2n(n-1)\hat s_n = -2n(n-1)\hat s_n
                \qquad\text{(by [\refBOZ, Eq.~(4.5)]).}
\endsplit
$$
Thus, for a given $n$, the image of $\De$ is a solution of the DOP problem
if and only if, for any $a,b<n$, $\Ga(s_a,s_b)$ does not contain any monomial
of the form $s_ks_n$ with $2\le k\le n$. This is the case for $n=4$.
The corresponding matrix $g$ is
$$
   \left(\matrix
     -16x^2 + 4x + 12y    & -24xy  + 8y  + 16z^2  & -16xz  + 6z   \\
     -24xy  + 8y + 16z^2  & -36y^2 + 4xy + 12z^2  & -24yz  + 4xz  \\
     -16xz  + 6z          & -24yz  + 4xz          & -16z^2 + y
   \endmatrix\right)
$$
The projection of the cuspidal edge of the surface $\deg g=0$
onto the $xz$-plane is the deltoid
(see Figure~\figABD\ and [\refBOZ, \S4.12]) up to affine
transformation of $\bR^2$.

If $n\ge 5$, then $\Phi_*(\De)$ is not a solution of the DOP problem
because, for example,
$$
   \Ga(s_4,s_{n-1})=-16(n-1)s_4 s_{n-1} + 4(n-4)s_3 s_{n-1} + 4(n-2) s_2\hat s_n^2
$$
has a monomial of degree $3$.
However, for any $n$, it is, evidently, a solution of the
weighted DOP problem with weights $(1,\dots,1,\tfrac12)$.


\subhead\sectJS. Direct products of Coxeter groups
\endsubhead

Let finite Coxeter groups $G_\alpha$,
$\alpha=1,\dots,m$,
act on vector spaces $E_\alpha$, $\dim E_\alpha=n_\alpha$.
We assume that these representations are irreducible or trivial.
Consider the diagonal action of $G=G_1\times\dots\times G_m$
on $E=\bigoplus_\alpha E_\alpha$. Let $n=\sum_\alpha n_\alpha$.
We denote the Laplace operator and the corresponding ``carr\'e du champ" on the
unit sphere in $E$ (resp. in $E_\alpha$) by $\De$ and $\Ga$
(resp. by $\De_\alpha$ and $\Ga_\alpha$).
Let $I_{\alpha,k}$, $k=1,\dots,n_\alpha$, be sets of basic invariant
homogeneous polynomials for the respective group actions,
$d_{\alpha,k}=\deg I_{\alpha,k}$.
We assume that $d_{\alpha,1}$ is minimal among the $d_{\alpha,k}$'s.
Then $d_{\alpha,1}=2$ unless $G_\alpha$ is trivial.

Let $g_\alpha^{ij}\big(x_{\alpha,1},\dots,x_{\alpha,{n_\alpha}}\big)$
and $b_\alpha^{i}\big(x_{\alpha,1},\dots,x_{\alpha,{n_\alpha}}\big)$ be the polynomials
such that
$$
   \Ga_\alpha(I_{\alpha,i},I_{\alpha,j})=
     g_\alpha^{ij}\big(I_{\alpha,1} ,\dots,I_{\alpha,n_\alpha}\big),
\qquad
   \De_\alpha(I_{\alpha,i})=
     b_\alpha^{i}\big(I_{\alpha,1} ,\dots,I_{\alpha,n_\alpha}\big).
$$
We assume that $\deg g^{ij}_\alpha\le 2$ and $b^i_\alpha\le 1$
for any $i,j,\alpha\ge 1$
(notice that this condition is fulfilled for
$A_n$ and $B_n$, but not for $D_4$; see \S\S\sectA--\sectD). 

\smallskip
{\it First construction.}
Suppose that $d_{1,1}=\dots=d_{m,1}=2$,
i.e. all the $I_{\alpha,1}$ are positive definite quadratic forms.
Let $\bS^{n-1}$ be the sphere in $E$ given by the equation
$\sum_\alpha I_{\alpha,1}=0$ and let
$\Phi:\bS^{n-1}\to\bR^{n-1}$ be the mapping defined by
$p\mapsto\big(\tilde I_1(p),I_2(p),\dots,I_m(p)\big)$, where
$I_\alpha=(I_{\alpha,1},\dots,I_{\alpha,n_\alpha})$ and
$\tilde I_1=(I_{1,2},\dots,I_{1,n_1})$.
It is easy to see that the image of
$\De$ through $\Phi$ is a solution of the DOP problem.
Denote the coordinates in the target space $\bR^{n-1}$ by
$(\tilde\xx_1,\xx_2,\dots,\xx_m)$ where
$\xx_\alpha=(x_{\alpha,1},\dots,x_{\alpha,{n_\alpha}})$ and
$\tilde\xx_\alpha=(x_{\alpha,2},\dots,x_{\alpha,{n_\alpha}})$.
Then the corresponding matrix
$g$ is the block matrix $(g_{\alpha\beta})_{\alpha,\beta=1}^m$ with the
block dimensions $(n_1-1,n_2,\dots,n_m)$ and the blocks
$g_{\alpha\beta}=\big(g_{\alpha\beta}^{ij}(\tilde\xx_1,\xx_2,\dots,\xx_m)\big)_{i,j}$
defined by
$$
  g_{\alpha\beta}^{ij}=\cases
     g_1^{ij}(1-x_{2,1}-\dots-x_{m,1},\tilde\xx_1), & \alpha=\beta=1,\\
     g_\alpha^{ij}(\xx_\alpha),   & \alpha=\beta\ge 2,\\
     -d_{\alpha,i}\,d_{\beta,j}\,x_{\alpha,i}\,x_{\alpha,j}, & \alpha\ne\beta.
  \endcases
$$
Up to affine linear change of coordinates, $\Phi_*(\De)$ does not depend on
the order of the summands. For example, 
if we exchange $E_1$ and $E_2$, then
the resulting solution is obtained from the initial one by
the affine linear change of coordinates 
$$
  (\tilde\xx_1,\xx_2,\dots,\xx_m)\mapsto
  (\tilde\xx_2,\,1-x_{2,1}-\dots-x_{m,1},\tilde\xx_1,\,\xx_3,\dots,\xx_m).
$$

\smallskip
{\it Second construction.}
Suppose now that $G_1$ is trivial, $d_{2,1}=\dots=d_{m,1}=2$, and $\dim E_m=1$.
Then $I_{1,1},\dots,I_{1,n_1}$ are just linear coordinates on $E_1$ and
$d_{1,1}=\dots=d_{1,n_1}=1$.
Let $\bS^{n-1}$ be the sphere in $E$ given by the equation
$\sum_i I_{1,i}^2+\sum_{\alpha\ge 2} I_{\alpha,1}=0$ and let
$\Phi:\bS^{n-1}\to\bR^{n-1}$ be the mapping defined by
$p\mapsto\big(I_1(p),\dots,I_{m-1}(p)\big)$, where
$I_\alpha=(I_{\alpha,1},\dots,I_{\alpha,n_\alpha})$.
Then $\Phi_*(\De)$ is a solution of the DOP problem.
Denote the coordinates in the target space $\bR^{n-1}$ by
$(\xx_1,\dots,\xx_{m-1})$ where
$\xx_\alpha=(x_{\alpha,1},\dots,x_{\alpha,{n_\alpha}})$.
Then the corresponding matrix
$g$ is the block matrix $\big(g(\alpha,\beta)\big)_{\alpha,\beta=1}^m$ with the
block dimensions $(n_1,n_2,\dots,n_{m-1})$ and with the blocks
$g(\alpha,\beta)=\big(g_{\alpha\beta}^{ij}\big)_{i,j}$
defined by
$$
  g_{\alpha\beta}^{ij}=\cases
     \delta^{ij}-x_{1,i}\, x_{1,j}, & \alpha=\beta=1,\\
     g_\alpha^{ij}(\xx_\alpha),   & 2\le\alpha=\beta\le m-1,\\
     -d_{\alpha,i}\,d_{\beta,j}\,x_{\alpha,i}\, x_{\alpha,j}, & \alpha\ne\beta.
  \endcases
$$


\subhead\sectAA. Quotient of $\bS^{n-1}$ and $\bS^n$ by the Coxeter group $A_1+A_{n-1}$
\endsubhead

Let the notation be as in \S\sectA. Let $\bS^{n-1}$ be the unit sphere in
$H_+=\bR\oplus H\subset\bR\oplus E$ (we denote the coordinate on $\bR$ by $x_0$).
Let $\De_+$ be the Laplace operator on $\bS^{n-1}$.
Consider the product $G$ of the Coxeter groups $A_1$ and $A_{n-1}$
diagonally acting on $H_+$.
%
According to \S\sectJS\ (first construction), the image of $\De_+$
through the mapping $\Phi_+:\bS^{n-1}\to\bR^{n-1}$,
$(x_0,x_1,\dots,x_n)\mapsto(s_2,\dots,s_n)$, provides a solution of the DOP problem
on the domain $\Phi_+(\bS^{n-1})$,
which is bounded by the hypersurface
$$
    (1+2X_2)F=0,
   \qquad F=\text{discr}_u\big(u^n+X_2 u^{n-2}+\dots+X_{n-1}u+X_n\big)=0.
                                                                      \eqno(\eqModelAA)
$$
Its component $1+2X_2=0$ is the image of $H\cap\bS^{n-1}$.
For $n=4$, this is the solution in Thm.~\thB(iii${}_3$)
(see Figure~\figAB) up to rescaling of the coordinates.
The entries of the matrix $g$ are given by the formulas in \S\sectB\ with
$s_0,s_1,\dots,s_n$ replaced
by $1,0,X_2,\dots,X_n$. We have $\De_+(s_a)=\De(s_a)-as_a$ with $\De(s_a)$ as in
(\eqModelAnIV).

\medskip
Let $\bS^n$ be the unit sphere in $\bR\oplus H_+$.
We denote the newly added coordinate by $\hat x_0$. Extend the above action
of $G$ to $\bR\oplus H_+$ assuming that it acts trivially on the first component.
Consider the image of $\De_{\bS^n}$ through
$(\hat x_0,x_0,x_1,\dots,x_n)\mapsto(\hat x_0,s_2,\dots,s_n)$.
According to \S\sectJS\ (second construction), it gives a solution of the DOP problem
in the domain in $\bR^n$ with coordinates $(X_1,\dots,X_n)$ bounded by the
hypersurface $(1+2X_2-X_1^2)F=0$; see (\eqModelAA).
We have $g^{1b}=\delta^{1b}-bX_1X_b$, $g^{ab}$ for $2\le a\le b$ are as above,
$\De_{\bS^n}(s_a)=\De(s_a)-2as_a$ ($\De(s_a)$ is as in (\eqModelAnIV)), and
$\De_{\bS^n}(\hat x_0)=-n\hat x_0$
For $n=3$ we obtain (i${}_6^*$) in Remark~\remCubicI\
up to rescaling.

The solution (i${}_6$) in Theorem~\thB\ and its generalization for higher dimensions
can be obtained as the image of $\De_{\bS^n}$ through a quotient by $A_1+A_{n-1}$
using a more direct (and somewhat more natural) construction as follows.
Let the notation still be as in \S\sectA.
Let $\bS^n$ be the unit sphere in $\bR\oplus E$.
Consider the mapping $\bS^n\to\bR^n$, $(x_0,x_1,\dots,x_n)\mapsto(s_1,\dots,s_n)$.
Its image is bounded by the hypersurface
$$
    (1+2X_2-X_1^2)
    \,\text{discr}_u\big(u^n+X_1 u^{n-1}+\dots+X_{n-1}u+X_n\big)=0.
$$
Using the computations in \S\sectA, for  $1\le a\le b\le n$, we obtain
$$
   g^{ab} = (n-b+1)s_{a-1}s_{b-1} - ab s_a s_b
            + \sum_{l\ge 1}(a-b-2l)s_{a-l-1} s_{b+l-1}
$$
with $s_0,\dots,s_n$ replaced by $1,X_1,\dots,X_n$
and $s_k=0$ for $k\not\in[0,n]$. We have 
$\Delta_{\bS^n}(s_a)=-a(n+a-1)s_a$.
When $n=3$,
we obtain the solution in Theorem~\thB(i${}_6$)
with $\alpha=1/2$ (see Figure~\figAiAii)
after rescaling $(x,y,z)=(3^{-1/2}X_1,X_2,3^{3/2}X_3)$.

\subhead\sectAB. Quotient of $\bS^n$ by the Coxeter group $A_1+B_n$
\endsubhead

Let the notation be as in \S\sectB. Let $\bS^n$ be the unit sphere in
$E_+=\bR\oplus E$
and let $\De_+$ be the Laplace operator on $\bS^n$.
 We denote the coordinate on $\bR$ by $x_0$
(recall that the coordinates on $E$ are $x_1,\dots,x_n$).
Consider the product of the Coxeter groups $A_1$ and $B_n$
diagonally acting on $E_+$. It is generated by the reflections in the hyperplanes
$x_i=0$ ($0\le i\le n$) and $x_i=x_j$ ($1\le i<j\le n$).

According to \S\sectJS\ (first construction), the image of $\De_+$
through the mapping $\Phi_+:\bS^n\to\bR^n$,
$(x_0,x_1,\dots,x_n)\mapsto(s_1,\dots,s_n)$, provides a solution of the DOP problem
on the domain $\Phi_+(\bS^n)$,
which is bounded by the hypersurface
$$
  X_n(1-X_1)\,\text{discr}_u\big(u^n+X_1 u^{n-1}+\dots+X_{n-1}u+X_n\big)=0.
$$
Its component $X_1=1$ is the image of $E\cap\bS^n$.
The other components are
as in \S\sectB. In the case $n=3$, this is the solution in Thm.~\thB(i${}_4$)
(see Figure~\figAB) up to rescaling of the coordinates.
The entries of $g$ are as in \S\sectB, but with
$s_0,s_1,\dots,s_n$ replaced
by $1,X_1,\dots,X_n$; $\De_+(s_a)=\De(s_a)-2a s_a$
(with $\De(s_a)$ as in (\eqModelBnDe)).


\subhead\sectCt. Quotient of $\bR^n$ by the affine Coxeter group $\widetilde C_n$
\endsubhead

Let $E=\bR^n$ with coordinates $\theta_1,\dots,\theta_n$.
The ring of invariant Fourier polynomials for the affine Coxeter group 
$\widetilde B_{n-1}$ is freely generated by $s_1,\dots,s_n$ where
$P(u)=(u+t_1)\dots(u+t_n)=\sum_{k=0}^ns_ku^{n-k}$, $t_i=\cos\theta_i$.
We consider the mapping $\Phi:\bR^n\to\bR^n$,
$(\theta_1,\dots,\theta_n)\mapsto(s_1,\dots,s_n)$. Its image $\Omega$
is the set of all $n$-tuples $X=(X_1,\dots,X_n)$ such that all roots of
the polynomial $P_X(u)=u^n+\sum_{k=1}^n X_k u^{n-k}$ are real and belong to
the interval $[-1,1]$. Therefore $\Omega$ is bounded by the union of the hypersurface
$\{X\mid\text{discr}_u P_X(u)=0\}$ and two hyperplanes $\{X\mid P_X(\pm1)=0\}$.
When the point $X$ moves from $\Omega$ crossing the discriminantal hypersurface,
two real roots disappears. When it crosses the hyperplane $X=\pm 1$, one of the roots
gets out from the interval $[-1,1]$.
One easily checks that $\Omega$ is the only bounded component of the complement
(cf.~Thm.~\thB(i${}_5$), Figure~\figCaff, \S\sectFig; for $n=2$ see [\refBOZ, \S4.7]).

By linearity, $\Ga(s_k,s_m)$ is the coefficient of $u^{n-k}v^{n-m}$ in $\Ga(P(u),P(v))$.
We have 
$$
 \Ga(t_i,t_j)=\Ga(\cos\theta_i,\cos\theta_j)=\delta_{ij}\sin\theta_i\sin\theta_j
    =\delta_{ij}(1-t_i^2).
$$
Hence
$$
   \Ga\big(P(u),P(v)\big)=\sum_{i,j}
   \big(\partial_{t_i} P(u)\big)\big(\partial_{t_j} P(v)\big)\Ga(t_i,t_j)
   =\sum_i\frac{P(u)P(v)(1-t_i^2)}{(u+t_i)(v+t_i)}
$$
$$
   =\frac{P(u)P(v)}{v-u}\sum_i\Big(\frac{1+ut_i}{u+t_i}-\frac{1+vt_i}{v+t_i}\Big)
   =\frac{Q(u)P(v)-Q(v)P(u)}{v-u}
$$
where
$$
\split
   Q(u)&=\sum_i\frac{P(u)(1+ut_i)}{u+t_i}
       =\sum_i P(u)\Big(u+\frac{1-u^2}{u+t_i}\Big)
 \\& = nuP(u) + (1-u^2)P'(u)
     = \sum_k s_k\Big( ku^{n-k+1} + (n-k)u^{n-k-1} \Big),
\endsplit
$$
thus $\Ga\big(P(u),P(v)\big)$ is equal to
$$
\split
 &\sum_{k,m} s_k s_m\! \left(
       k\frac{u^{n-k+1}v^{n-m} - v^{n-k+1}u^{n-m}}{v-u}
      + (n-k)\frac{u^{n-k-1}v^{n-m} - v^{n-k-1}u^{n-m}}{v-u} \right)
\\&\qquad=\sum_{k,m} s_k s_m \Bigg\{
   k\left(\sum_{l=1}^{k-m-1}u^{n-k+l}v^{n-m-l}
         -\sum_{l=0}^{m-k}  u^{n-k-l}v^{n-m+l}\right)
\\&\qquad\qquad\quad
     +(n-k)\left(\sum_{l=0}^{k-m}  u^{n-k+l-1}v^{n-m-l-1}
                -\sum_{l=1}^{m-k-1}u^{n-k-l-1}v^{n-m+l-1}\right)\Bigg\}.
\endsplit
$$
Hence, for $a\le b$, we have
$$
   \Ga(s_a,s_b) = (n-b+1)s_{a-1}s_{b-1} - as_as_b
   + \sum_{l\ge 1}(b-a+2l)(s_{a-l}s_{b+l} - s_{a-l-1}s_{b+l-1})
$$
It is easy to see that $\De(s_a) = -a s_a$.


\subhead\sectBDt.
           Quotients of $\bR^n$ by the affine Coxeter groups $\widetilde B_n$
           and  $\widetilde D_n$
\endsubhead

Let the notation be as in \S\sectCt.
The ring of invariant Fourier polynomials for the affine Coxeter group 
$\widetilde B_n$ is freely generated by $s_1,\dots,s_{n-1}$
and 
$$
  \hat s_n
    = \prod_{i=1}^n \sqrt{2}\,\cos(\theta_i/2) = \prod_{i=1}^n\sqrt{1+\cos\theta_i}
    = P(1)^{1/2}.
$$
$\Gamma(s_a,s_b)$ for $a\le b\le n-1$ are as in \S\sectCt\ with the substitution
$s_n=\hat s_n^2-\sum_{k=0}^{n-1}s_k$ (recall that $s_0=1$).
Using the computations in \S\sectCt\ 
we obtain
$$
   2\Ga(P,\hat s_n) = \frac{\Ga(P(u),P(1))}{P(1)^{1/2}}
   =\frac{Q(u)P(1) - Q(1)P(u)}{P(1)^{1/2}(1-u)}
   = \hat s_n\big((1+u)P'-nP\big),
$$
hence
$2\Ga(s_a,\hat s_n) = ((n-a+1)s_{a-1} - a s_a)\hat s_n$ and
$$
\split
   4\Ga(\hat s_n,\hat s_n)&=\frac{\Ga(P(1),P(1))}{P(1)}
   =\frac{1}{P(1)}\sum_i\frac{P(1)^2(1-t_i^2)}{(1+t_i)^2}
   =\sum_i\left(\frac{2P(1)}{1+t_i}-P(1)\right)
   \\&=2P'(1)-nP(1)
   =-n\hat s_n^2+\sum_{k=0}^{n-1}(n-k)s_k. 
\endsplit
$$
We have $\De(s_a)=-as_a$ and $\De(\hat s_n)=-\tfrac14 n\hat s_n$.

The image of $\De$ through $\Phi:\bR^n\to\bR^n$,
$(\theta_1,\dots,\theta_n)\mapsto(s_1,\dots,s_{n-1},\hat s_n)$
is not a solution of the DOP problem when $n\ge 3$
(and $\widetilde B_2$ is the same as $\widetilde C_2$).
Indeed, $\Ga(s_2,s_{n-1})$ has monomial $(n-1)s_1\hat s_n^2$ of degree $3$.
However, $\Phi_*(\De)$ is a solution of the weighted DOP problem with
weights $(1,\dots,1,\frac12)$.


\smallskip
For the affine group $\widetilde D_n$,
all the computations are almost the same and we omit the details.
The ring of invariant Fourier polynomials
is generated by $s_1,\dots,s_{n-2}$, $\hat s_n$,
and $\hat s_{n-1}=\prod_{i=1}^n\sqrt2\,\sin(\theta/2)=\sqrt{(-1)^nP(-1)}$.
For $a\le b\le n-2$, $\;\Ga(s_a,s_b)$, $\Ga(s_a,\hat s_n)$, and $\Ga(\hat s_n,\hat s_n)$
are the same as above but with the substitutions
$$
     s_n=\tfrac12\big(\hat s_n^2+(-1)^n\hat s_{n-1}^2\big) - \sum_{k\ge 1}s_{n-2k},
 \quad
 s_{n-1}=\tfrac12\big(\hat s_n^2+(-1)^n\hat s_{n-1}^2\big) - \sum_{k\ge 1}s_{n-2k-1}.
$$
\if01{
$$
     s_n=\frac{\hat s_n^2+(-1)^n\hat s_{n-1}^2}2 - \sum_{k\ge 1}s_{n-2k},
 \quad
 s_{n-1}=\frac{\hat s_n^2+(-1)^n\hat s_{n-1}^2}2 - \sum_{k\ge 1}s_{n-2k-1}.
$$
}\fi
\if01{
We have
$2\Gamma( P, \hat s_{n-1})=\hat s_{n-1}((1-u)P'+nP)$.
Finally,
$$
\split
  2\Gamma(s_a,\hat s_{n-1})&=((n-a+1)s_{a-1}+a s_a)\hat s_{n-1},\quad a=1,\dots,n-2,\\
  2\Gamma(s_a,\hat s_{n})  &=((n-a+1)s_{a-1}-a s_a)\hat s_{n-1},\quad a=1,\dots,n-2,\\
  4\Gamma(\hat s_{n-1},\hat s_{n-1})&= \\
  4\Gamma(\hat s_{n-1},\hat s_{n})  &= \\
  4\Gamma(\hat s_{n},  \hat s_{n})  &= 
\endsplit
$$
}\fi
The values of $\Ga(\hat s_{n-1},*)$ are computed similarly
and we arrive to the same conclusion as above:
the quotient by $\widetilde D_n$, $n\ge 4$,
does not provide a solution of the DOP problem, but it provides a solution
of the weighted DOP problem with weights $(1,\dots,1,\frac12,\frac12)$.


\subhead\sectAt.
           Quotient of $\bR^{n-1}$ by the affine Coxeter group $\widetilde A_{n-1}$
\endsubhead

Let $E$ be $\bR^n$ with coordinates $\theta_1,\dots,\theta_n$ and
$H=\{\theta_1+\dots+\theta_n=0\}$.
The affine Coxeter group $\widetilde A_{n-1}$ acting on $H$ is generated by the
orthogonal reflections in the hyperplanes $x_i=x_j$ and a suitable translation.
The ring of invariant Fourier polynomials is freely generated by
$s_1,\dots,s_{n-1}$
where $P(u)=(u+t_1)\dots(u+t_n)=\sum_{k=0}^ns_ku^{n-k}$, $t_i=\exp(\ii\,\theta_i)$,
$\ii=\sqrt{-1}$.
Notice that $\bar s_k|_H=s_{n-k}|_H$, in particular $s_n|_H=1$ and $s_{n/2}|_H$ is real
when $n$ is even.
We consider the mapping $\Phi:H\to\bR^{n-1}$,
$(\theta_1,\dots,\theta_n)\mapsto s=(s_1,\dots,s_{\lfloor n/2\rfloor})$
where we identify $\bR^{n-1}$ with $\bC^{(n-1)/2}$ (which we define as
$\bC^{(n-2)/2}\times\bR$ when $n$ is even).
Then $\Phi(H)$ is bounded by the hypersurface
$\text{discr}_u P(u,Z)=0$, $Z=(Z_1,\dots,Z_{\lfloor n/2\rfloor})\in\bC^{(n-1)/2}$,
$$
  P(u,Z)=\cases
 \!u^{2k}+Z_1 u^{2k-1}+\dots+Z_k u^k + \bar Z_{k-1}u^{k-1}+\dots+\bar Z_1 u+1, &\!n=2k,\\
 \!u^{2k+1}+Z_1 u^{2k}+\dots+Z_k u^k + \bar Z_k u^{k+1}+\dots+\bar Z_1 u+1, &\!n=2k+1
  \endcases
$$
\if01{$$
  P(u,z)=\cases
 \!u^{2k}+z_1 u^{2k-1}+\dots+z_k u^k + \bar z_{k-1}u^{k-1}+\dots+\bar z_1 u+1, &\!n=2k,\\
 \!u^{2k+1}+z_1 u^{2k}+\dots+z_k u^k + \bar z_k u^{k+1}+\dots+\bar z_1 u+1, &\!n=2k+1
  \endcases
$$}\fi
(cf.~Thm.~\thB(vi) and Figure~\figAaff; for $n=2$ see [\refBOZ, \S 4.12]).
Let $\De=\De_H$ and let $\Ga$ be the associated carr\'e du champ.
Then $\Ga(s_k,s_m)$ is the coefficient of $u^{n-k}v^{n-m}$
in $\Ga(P(u),P(v))$.
For any functions $f$, $g$ we have
$\Ga(f,g)=\Ga_E(f,g)-\tfrac1n(\partial_0f)(\partial_0g)$ where
$\partial_0 = \sum_{i=1}^n \tfrac{\partial}{\partial\theta_i}$ (see (\eqModelAnI)).
Denote $\tfrac\partial{\partial t_i}$ by $\partial_i$.
We have
$$
\split
   \partial_0 P(u)&=\sum_i\ii\, t_i\partial_i P = \ii\sum_i\frac{t_iP}{u+t_i}
    =\ii\sum_i\Big(1-\frac{u}{u+t_i}\Big)P
    = \ii\big(nP(u)-uP'(u)\big),  
\endsplit
$$
thus
$
  \Ga\big(P(u),P(v)\big) = \Ga_E\big(P(u),P(v)\big)
    + \tfrac1n\big(nP(u)-uP'(u)\big)\big(nP(v)-vP'(v)\big).
$
We also have
$\Ga_E(t_i,t_j)=(dt_i/d\theta_i)(dt_j/d\theta_j)\Gamma_E(\theta_i,\theta_j)
=-\delta_{ij}t_i^2$.
Then
$$
   -\Ga_E\big(P(u),P(v)\big)
   =-\sum_{i,j}\big(\partial_i P(u)\big)\big(\partial_j P(v)\big)\Ga_E(t_i,t_j)
   =\sum_i\frac{P(u)P(v)t_i^2}{(u+t_i)(v+t_i)}
$$
$$
   =\frac{P(u)P(v)}{v-u}\sum_i\Big(v-u+\frac{u^2}{u+t_i}-\frac{v^2}{v+t_i}\Big)
   =nP(u)P(v) + \frac{u^2P'(u)P(v) - v^2P'(v)P(u)}{v-u}
$$
$$
   =nP(u)P(v)+\sum_{k,m}(n-k)s_k s_m
            \frac{u^{n-k+1}v^{n-m} - v^{n-k+1}v^{n-m}}{v-u}
$$
$$
  =nP(u)P(v)+\sum_{k,m}(n-k)s_k s_m
        \left( \sum_{l=1}^{k-m-1}u^{n-k+l}v^{n-m-l}
         -\sum_{l=0}^{m-k}  u^{n-k-l}v^{n-m+l}\right)
$$
Hence, for $a\le b$, we have
$$
   \Ga(s_a,s_b) = \frac{a(b-n)}{n}s_a s_b + \sum_{l\ge 1}(b-a+2l)s_{a-l}s_{b+l}.
$$
Setting $s_k=x_k+\ii y_k$, $\Ga(s_a,s_b) = A+\ii B$,
and $\Ga(s_a,\bar s_b) = C+\ii D$, we obtain
for $a\le b\le n/2$:
$$
  2\Ga(x_a,x_b)=A+C,\;
  2\Ga(x_a,y_b)=B-D,\;
  2\Ga(y_a,x_b)=B+D,\;
  2\Ga(y_a,y_b)=C-A,
$$
$$
\split&
   A
    =\frac{a(b-n)}{n}(x_ax_b-y_ay_b)
                     +\sum_{l\ge1}(b-a+2l)(x_{a-l}x_{b+l}-y_{a-l}y_{b+l}),
\\&
   B=\frac{a(b-n)}{n}(x_ay_b+y_ax_b)
                     +\sum_{l\ge1}(b-a+2l)(x_{a-l}y_{b+l}+y_{a-l}x_{b+l}),
\\&
   C=-\frac{ab}{n}(x_ax_b+y_ay_b)
                     +\sum_{l\ge1}(n-a-b+2l)(x_{a-l}x_{b-l}+y_{a-l}y_{b-l}),
\\&
   D=\frac{ab}{n}(x_ay_b - y_ax_b)
                    -\sum_{l\ge1}(n-a-b+2l)(x_{a-l}y_{b-l}-y_{a-l}x_{b-l}).
\endsplit
$$
For $a\le n/2$ we have $\De(x_a) = \lambda_a x_a$,
    $\De(y_a) = \lambda_a y_a$, $\lambda_a = a(a-n)/n$.

\if01{ (* D4 *)
G={{ -16*X^2 + 4*X + 12*Y,    -24*X*Y + 8*Y   + 16*Z^2,  -16*X*Z + 6*Z    },
   { -24*X*Y + 8*Y + 16*Z^2,  -36*Y^2 + 4*X*Y + 12*Z^2,  -24*Y*Z + 4*X*Z  },
   { -16*X*Z + 6*Z,           -24*Y*Z + 4*X*Z,           -16*Z^2 + Y      }}
F=Factor[Det[G]]
Factor[(G.{D[F,X],D[F,Y],D[F,Z]})/F]

v[1]=X; v[2]=Y; v[3]=Z;
rho=F^(-1/2)
De=Factor[(Sum[D[G[[i]],v[i]],{i,3}]+G.Table[D[Log[rho],v[i]],{i,3}])]

}\fi

\if01{

(**  tilde A3  **)
n=7
Pu=Sum[s[k]u^(n-k),{k,0,n}]; dPu=D[Pu,u];
Pv=Sum[s[k]v^(n-k),{k,0,n}]; dPv=D[Pv,v];

GaEPP = Factor[ -n*Pu*Pv - (u^2*dPu*Pv - v^2*dPv*Pu)/(v-u) ];
GaPP = Factor[ GaEPP + (n*Pu-u*dPu)(n*Pv-v*dPv)/n ]/.{s[0]->1,s[n]->1};

Factor[GaEPP + n*Pu*Pv
  + Sum[(n-k)s[k]s[m]( Sum[u^(n-k+l)v^(n-m-l),{l,1,k-m-1}]
                      -Sum[u^(n-k-l)v^(n-m+l),{l,0,m-k}] ),
    {k,0,n},{m,0,n}]
]

su0=Join@@Table[{s[-k]->0,s[n+k]->0},{k,n}];
Table[ Coefficient[Coefficient[GaEPP,u,n-a],v,n-b] +
        a*s[a]s[b] + Sum[(a-b-2l)s[a-l]s[b+l],{l,1,n}]/.su0,
{a,n-1},{b,a,n-1} ]

Expand[Table[ -Coefficient[Coefficient[GaPP,u,n-a],v,n-b] +
        a*(b-n)/n*s[a]s[b] + Sum[(b-a+2l)s[a-l]s[b+l],{l,1,n}]/.su0,
{a,n-1},{b,a,n-1} ]/.{s[0]->1,s[n]->1}]

G=Table[Coefficient[Coefficient[GaPP,u,n-a],v,n-b]n,{a,n-1},{b,n-1}]
F=Factor[Det[G]]
Factor[(G.Table[D[F,s[i]],{i,n-1}])/F]

(* Ga(xa,xb), a<=b *)
Expand[ (a(n-b)sa*sb + a*b*sa*sbb + a*b*saa*sb + a(n-b)saa*sbb)/.
    {sa->xa+ya*I,saa->xa-ya*I,sb->xb+yb*I,sbb->xb-yb*I}]
Expand[ ( (b-a+2l)(sa*sb+saa*sbb) + (n-a-b+2l)(sa*sBB + saa*sB) )/.
    {sa->xa+ya*I,saa->xa-ya*I,sb->xb+yb*I,sbb->xb-yb*I,sB->xB+yB*I,sBB->xB-yB*I}]

(* Ga(ya,ab), a<=b *)
Expand[ -(a(n-b)sa*sb - a*b*sa*sbb - a*b*saa*sb + a(n-b)saa*sbb)/.
    {sa->xa+ya*I,saa->xa-ya*I,sb->xb+yb*I,sbb->xb-yb*I}]
Expand[ -( (b-a+2l)(sa*sb+saa*sbb) - (n-a-b+2l)(sa*sbb + saa*sb) )/.
    {sa->xa+ya*I,saa->xa-ya*I,sb->xb+yb*I,sbb->xb-yb*I}]

}\fi

\if01{
(* D4 *)
n=4;
P = Product[u+x[i]^2,{i,n}]
tn = Product[x[i],{i,n}]
Do[s[k]=Coefficient[P,u,n-k],{k,n}]
va=Table[x[i],{i,n}];
FunGa=Function[{f,g},Sum[D[f,v]D[g,v],{v,va}]]
Expand[FunGa[s[3],s[3]]-4s[2]s[3]-12s[1]s[4]]
Expand[FunGa[s[3],s[2]]-8s[3]s[1]-16s[4]]
Expand[FunGa[s[2],s[2]]-4s[2]s[1]-12s[3]]
Expand[FunGa[s[3],tn]-4s[2]tn]
Expand[FunGa[s[2],tn]-6s[1]tn]
Expand[FunGa[tn,tn]-s[3]]

n=3;
P = Product[u+x[i],{i,n}]
Do[s[k]=Coefficient[P,u,n-k],{k,n}]
va=Table[x[i],{i,n}];
FunGa=Function[{f,g},Sum[D[f,v]D[g,v],{v,va}]]
Expand[FunGa[s[1],s[1]]]

n=3;
P = Product[u+x[i]^2,{i,n}]
Do[s[k]=Coefficient[P,u,n-k],{k,n}]
va=Table[x[i],{i,n}];
FunGa=Function[{f,g},Sum[D[f,v]D[g,v],{v,va}]]
Expand[FunGa[s[1],s[1]]]

(* D5 *)
n=5;
P = Product[u+x[i]^2,{i,n}]
tn = Product[x[i],{i,n}]
Do[s[k]=Coefficient[P,u,n-k],{k,n}]
va=Table[x[i],{i,n}];
FunGa=Function[{f,g},Sum[D[f,v]D[g,v],{v,va}]]
Expand[FunGa[s[4],s[4]]]
Expand[FunGa[s[3],s[3]]-4s[2]s[3]-12s[1]s[4]-20s[5]]
Expand[FunGa[s[3],s[2]]-8s[3]s[1]-16s[4]]
Expand[FunGa[s[2],s[2]]-4s[2]s[1]-12s[3]]
Expand[FunGa[s[3],tn]-4s[2]tn]
Expand[FunGa[s[2],tn]-6s[1]tn]
Expand[FunGa[tn,tn]-s[3]]

(* B4 *)
n=4;
P = Product[u+x[i]^2,{i,n}]
Do[s[k]=Coefficient[P,u,n-k],{k,n}]
va=Table[x[i],{i,n}];
FunGa=Function[{f,g},Sum[D[f,v]D[g,v],{v,va}]]
Expand[FunGa[s[4],s[4]]-4s[3]s[4]]
Expand[FunGa[s[4],s[3]]-8s[2]s[4]]
Expand[FunGa[s[4],s[2]]-12s[1]s[4]]
Expand[FunGa[s[3],s[3]]-4s[2]s[3]-12s[1]s[4]]
Expand[FunGa[s[3],s[2]]-8s[3]s[1]-16s[4]]
Expand[FunGa[s[2],s[2]]-4s[2]s[1]-12s[3]]
G={{4X+12Y-4X*4X,    8Y+16Z-4X*6Y,  12Z-4X*8Z},
   {8Y+16Z-4X*6Y,  4X*Y+12Z-6Y*6Y, 8X*Z-6Y*8Z},
   {   12Z-4X*8Z,      8X*Z-6Y*8Z, 4Y*Z-8Z*8Z}}
F=Factor[Det[G]]
Factor[(G.{D[F,X],D[F,Y],D[F,Z]})/F]
Factor[Discriminant[F,Z]]

Factor[F/Discriminant[T^4-T^3+X*T^2-Y*T+Z,T]]  (** 64 Z **)

(* tilde B3 *)
n=3
P = Product[u+(x[i]^2+1/x[i]^2)/2,{i,n}]
Do[s[k]=Coefficient[P,u,n-k],{k,n}]
s[n] = Product[x[i]+1/x[i],{i,n}]
va=Table[x[i],{i,n}];
FunGa=Function[{f,g},Sum[-v^2 D[f,v]D[g,v],{v,va}]]

Expand[FunGa[s[1],s[1]] + 4*s[1]^2 - 8*s[2] - 12]
Expand[FunGa[s[1],s[2]] + 4*s[1]s[2] - 3/2*s[3]^2 + 12*s[2] + 4*s[1] + 12]
Expand[FunGa[s[1],s[3]] + 2*s[1]s[3] - 6*s[3]]
Expand[FunGa[s[2],s[2]] + 8*s[2]^2 - s[1]*s[3]^2 + 8*s[1]s[2] + 8*s[2] + 8*s[1]]
Expand[FunGa[s[2],s[3]] + 4*s[2]s[3] - 4*s[1]s[3]]
Expand[FunGa[s[3],s[3]] + 3*s[3]^2 - 16*s[2] - 32*s[1] - 48]
g11 = -4X^2+8Y+12; g12 = -4X*Y+3/2Z^2-12Y-4X-12; g13 = -2X*Z+6Z;
                   g22 = -8Y^2+X*Z^2-8X*Y-8Y-8X; g23 = -4Y*Z+4X*Z;
                                                 g33 = -3Z^2+16Y+32X+48;
G={{g11,g12,g13},{g12,g22,g23},{g13,g23,g33}}
F=Factor[Det[G]]
Factor[(G.{D[F,X],D[F,Y],D[F,Z]})/F]

(* tilde D4 *)
n=4
P = Product[u+(x[i]^2+1/x[i]^2)/2,{i,n}]
Do[s[k]=Coefficient[P,u,n-k],{k,n}]
s[n-1] = Product[x[i]-1/x[i],{i,n}]
s[n]   = Product[x[i]+1/x[i],{i,n}]
va=Table[x[i],{i,n}];
FunGa=Function[{f,g},Sum[-v^2 D[f,v]D[g,v],{v,va}]]
sp=Expand[(s[4]+s[3])/2]
sm=Expand[(s[4]-s[3])/2]
Expand[FunGa[s[2],s[2]]+8s[2]^2-sm*sp*s[1]-4s[1]^2-sp^2-sm^2+24s[2]+16]

(* tilde F4 *)

(* Bn (Delta in the 2nd proj) *)

n=6
Do[t[i]=n*x[i]^2-1,{i,n}];
rdr = Function[F, Expand[Sum[x[i]D[F,x[i]],{i,n}]] ];
rdr2 = Function[F,rdr[rdr[F]]];
De = Function[F, Expand[Sum[D[F,x[i],x[i]],{i,n}]]];
P=Product[u + t[i], {i,n}];
dP=D[P,u];
ddP=D[dP,u];

Expand[ 2n*dP - De[P] ]
Expand[ 2n*P+2(1-u)dP - rdr[P] ]

Expand[ 2n*(2n*P+2(1-u)dP)+2(1-u)(2(n-1)*dP+2(1-u)ddP) - rdr2[P] ]  (* 0 *)
Expand[ (4n^2*P+4n(1-u)dP)+(4(n-1)(1-u)dP+4(1-u)^2ddP) - rdr2[P] ]  (* 0 *)

Expand[ 4n^2*P + 4(2n-1)(1-u)dP + 4(1-u)^2ddP - rdr2[P] ]  (* 0 *)

Expand[ 4n^2*P + 4(2n-1)(1-u)dP + 4(1-u)^2ddP + 2n(n-1)P + 2(n-1)(1-u)dP
   - rdr2[P] - (n-1)rdr[P] ]                                                (** 0 **)

Expand[2(3n^2-n)*P + 2(5n-3)(1-u)dP + 4(1-u)^2*ddP - rdr2[P] - (n-1)rdr[P] ] (** 0 **)

Expand[2n*dP - ( 2(3n^2-n)*P + 2(5n-3)(1-u)dP + 4(1-u)^2*ddP ) -
   (De[P] - rdr2[P] - (n-1)rdr[P]) ]                               (** 0 **)

Expand[ - ( 2(3n^2-n)*P + 2((5n-3)(1-u)-n)dP + 4(1-u)^2*ddP ) -
   (De[P] - rdr2[P] - (n-1)rdr[P]) ]                               (** 0 **)

Expand[ - ( 2(3n^2-n)*P + 2(4n-3)dP + 2(-5n+3)u*dP + 4(1-u)^2*ddP ) -
   (De[P] - rdr2[P] - (n-1)rdr[P]) ]                               (** 0 **)

Do[s[k] = Coefficient[P,u,n-k], {k,-2,n+2}]

Factor[ Sum[(2(-k*n-2k^2+k)s[k] + 2(n-k+1)(3-4k)s[k-1]
            - 4(n-k+1)(n-k+2)s[k-2])u^(n-k), {k,0,n} ] -
   (De[P] - rdr2[P] - (n-1)rdr[P]) ]                               (** 0 **)

Factor[ P*Sum[x[i]^2/(u+t[i]),{i,n}] - ((1-u)/n*dP + P) ] (** 0 **)
}\fi


\head\sectCone. Conical surfaces 
\endhead


\proclaim{ Theorem \thC }
Let $(\Omega,g,\rho)$ be a solution of the SDOP problem in $\bC^3$ such that
$\partial\Omega$ contains a relatively open subset
of an irreducible conical surface $\Sigma$,
i.e. of a surface $\Sigma=\{\Gamma(x,y,z)=0\}$
where $\Gamma$ is an irreducible homogeneous polynomial.
Then $\deg\Gamma\le 2$.
\endproclaim

\noindent
{\bf Remark \remC.}
There exist solutions of the DOP problem in bounded domains
whose boundaries contain a piece of the quadratic cone. For example,
the solutions (1b), (1e), (3e), (3i), (5g), (6f) in [\refBOZ, \S7.2]
(see Figure~\figCone).
\medskip

\midinsert
\centerline{%
\noindent\hskip0pt
\epsfysize=23mm\epsfbox{1b.eps}
\epsfysize=23mm\epsfbox{1e.eps}
\epsfysize=23mm\epsfbox{3e.eps}
\epsfysize=23mm\epsfbox{3i.eps}
\epsfysize=23mm\epsfbox{5g.eps}
\epsfysize=23mm\epsfbox{6f.eps}}
\vskip-7pt
\centerline{
(1b)\hskip 20mm
(1e)\hskip 18mm
(3e)\hskip 12mm
(3i)\hskip 9mm
(5g)\hskip 7mm
(6f)\hskip 10mm}
\botcaption{ Figure~\figCone } Bounded domains $\Omega$ from [\refBOZ, \S7.2]
                            admitting solutions of the DOP problem, such that
                            $\partial\Omega$ contains a piece of the quadratic cone.
\endcaption
\endinsert

\proclaim{ Proposition \propC }
Let $(g,\Gamma)$ be a solution of the AlgDOP problem in $\bC^3$ such that
$\Gamma$ is an irreducible homogeneous polynomial and
$\det g$ is not a homogeneous polynomial of degree $6$.
Then $\deg\Gamma\le 2$.
\endproclaim

One can easily derive Theorem~\thC\ from Proposition~\propC.
Indeed, if $\det g$ (in the setting of Theorem~\thC) were
a homogeneous polynomial of degree $6$, then
affine coordinates $(x,y,z)$ could be chosen so that $\deg_x\det g=6$ and
$\Omega$ contains a half-cylinder $\{x>0, y^2+z^2<1\}$, which
contradicts [\refBOZ, Cor.~2.19].

\if01{
eq={g+n+2a4+3a6+4a8+5a10+6a12==(d-1)(d-2)/2,                                       
    dd==d(d-1)-2n-5a4-7a6-9a8-11a10-13a12,                                       
    2-2g==2dd-d-a4-a6-a8-a10-a12-2p}
Solve[eq/.{d->5,p->0},{dd,n,g}]
(***
{{dd -> 5,
   n -> (15 - 11*a10 - 13*a12 - 5*a4 - 7*a6 - 9*a8)/2, 
   g -> (-3 + a10 + a12 + a4 + a6 + a8)/2 }}
***)
Solve[eq/.{d->4,a8->0,a10->0,a12->0},{dd,n,g}]
(***
{{dd -> (2*(6 + p))/3,
   n -> (24 - 15*a4 - 21*a6 - 2*p)/6, 
   g -> (-6 + 3*a4 + 3*a6 + 2*p)/6}}
***)
Solve[eq/.{d->3,a8->0,a10->0,a12->0},{dd,n,g}]
}\fi

The rest of this section is devoted to the proof of Proposition~\propC.
Let $\Gamma$ be as in Proposition~\propC.
Let $\Sigma$ be the surface in $\bC^3$ defined by the equation $\Gamma=0$, and
let $C$ be the curve in $\bP^3=\bP(\bC^3)$ defined by the same equation.
Any local branch $\gamma$ of $C$ has a parametrization of the form
$t\mapsto(t^p,t^q+o(t^q))$, $1\le p<q$, in some affine coordinates.
We then say that $\gamma$ is of type $(p,q)$. 

\proclaim{ Lemma \lemC }
Any local branch of $C$ is of type $(1,2)$ or $(2,4)$.
\endproclaim

\demo{Proof}
The arguments are as in \S\sectLoc\ but simpler.
Let $\pi:\bC^3\setminus\{(0,0,0)\}\to\bP^2$ be the quotient map
(then $\Sigma=\pi^{-1}(C)$).
Let $\gamma$ be a local branch of $C$ at $p\in\CP^2$ parametrized by
$t\mapsto\gamma(t)=(\xi_1(t):\xi_2(t):\xi_3(t))$.
Then $\Sigma$ near the line $\pi^{-1}(p)$ is parametrized by
$(t,u)\mapsto(u\xi_1(t),u\xi_2(t),u\xi_3(t))$.
Similarly to \S\sectLoc, we rewrite the equations (\eqParam)
in the form $E_1=E_2=E_3=0$ where
$$
   E_i 
    = u\sum_{j=1}^3 (\dot\xi_{j+1}\xi_{j-1} - \dot\xi_{j-1}\xi_{j+1})
         g^{ij}(u\xi_1,u\xi_2,u\xi_3)
    = \sum_{\alpha=0}^\infty t^\alpha\sum_{\beta=1}^3 E_{\alpha,\beta,i}\,u^\beta
$$
($i$ and $j$ are considered mod $3$)
and the $E_{\alpha,\beta,i}$ are linear forms in $g_{klm}^{ij}$ whose coefficients are
polynomial functions of the coefficients of the $\xi_i$'s.

We have $\deg C\le 5$ because otherwise $\det g$ would be homogeneous of degree $6$.
Hence $C$ may have only local branches of type $(p,q)$ with $q\le 5$.
For each pair $(p,q)$, $1\le p<q\le 5$, except $(1,2)$ and $(2,4)$
(thus for $8$ pairs) we consider a branch $\gamma$ of type $(p,q)$ of the form
$t\mapsto(1:t^p:t^q+\sum_{k>q}a_kt^k)$ with indeterminate coefficients $a_k$
and solve the maximal triangular subsystem of the system of
equations $E_{\alpha,\beta,i}=0$ for the unknowns $g_{klm}^{ij}$.
This means that we find an equation implying that some unknown is zero,
replace this unknown by zero in all other equations,
and repeat this process as long as we can do.
%
For all pairs $(p,q)$ except $(1,4)$, $(1,5)$ we obtain that
$\deg g$ is homogeneous of degree $6$. In the two exceptional cases
we obtain that $z^2$ divides $\det g$. Since $q\le\deg\Gamma$, this implies
$\deg(z^2\Gamma)\ge 6$, hence $\det g=z^2\Gamma$ up to a scalar factor.
This means that $\det g$ is homogeneous of degree $6$. 
\qed\enddemo

Let $d$ and $\gg$ be the degree and the genus of $C$ respectively.
Let $a_{2k}$, $k\ge 2$, be the number of local branches of type $(2,4)$
which admit a parametrization $t\mapsto(t^2,t^{2k+1})$ in some local curvilinear
coordinates (the $A_{2k}$-singularity).
Let $\check d$ be the degree of the projectively dual curve $\check C$.
Let
$
  \nn = \sum_{\gamma}\delta(\gamma) + \sum_{\gamma_1,\gamma_2}(\gamma_1\cdot\gamma_2)
$
where $\gamma$ runs over all local branches of $C$, $\delta(\gamma)$ is the
delta-invariant of $\gamma$, and $(\gamma_1,\gamma_2)$ runs over
all unordered pairs of local branches  (see [\refBOZ, \S3.2] for more details).
Due to Lemma~\lemC, the Pl\"ucker-like equations [\refBOZ, Eqs.~(3.13)--(3.15)]
take the form
$$
\split
          &\gg + \nn + \sum k a_{2k} = (d-1)(d-2)/2, \\
          &\check d = d(d-1) - 2\nn - \sum(2k+1)a_{2k}, \\
          &2 - 2\gg = 2\check d - d - \sum a_{2k}
\endsplit
$$
(all the summations run over $k\ge 2$).
One easily checks that these equations do not have any integer
non-negative solution with $3\le d\le 5$. Proposition~\propC\
is proven.

\enddocument